\documentclass{amsart}
\usepackage{amssymb}

\usepackage{amscd}
\usepackage{thmdefs}
\usepackage[dvips]{epsfig}

\input{tcilatex}
\theoremstyle{definition}
\theoremstyle{remark}
\numberwithin{equation}{section}

\input tcilatex

\begin{document}
\title[Fractaloids and Labeling Operators]{Applications of Automata and Graphs: Labeling-Operators in Hilbert Space 
\textrm{I}}
\author{Ilwoo Cho and Palle E. T. Jorgensen}
\address{Saint Ambrose Univ., Dep. of Math., 421 Ambrose Hall, 518 W. Locust St.,
Davenport, Iowa, 52803, U. S. A. / Univ. of Iowa, Dep. of Math, 25McLean
Hall, Iowa City, Iowa, 52242, U. S. A. }
\email{chowoo@sau.edu\\
Jorgense@math.uiowa.edu}
\thanks{}
\date{Mar., 2008.}
\dedicatory{}
\thanks{}
\maketitle

\begin{abstract}
We show that certain representations of graphs by operators on Hilbert space
have uses in signal processing and in symbolic dynamics. Our main result is
that graphs built on automata have fractal characteristics. We make this
precise with the use of Representation Theory and of Spectral Theory of a
certain family of Hecke operators. Let $G$ be a directed graph. We begin by
building the graph groupoid $\Bbb{G}$ induced by $G$, and representations of 
$\Bbb{G}.$ Our main application is to the groupoids defined from automata.
By assigning weights to the edges of a fixed graph $G,$ we give conditions
for $\Bbb{G}$ to acquire fractal-like properties, and hence we can have
fractaloids or $G$-fractals. Our standing assumption on $G$ is that it is
locally finite and connected, and our labeling of $G$ is determined by the
``out-degrees of vertices''. From our labeling, we arrive at a family of
Hecke-type operators whose spectrum is computed. As applications, we are
able to build representations by operators on Hilbert spaces (including the
Hecke operators); and we further show that automata built on a finite
alphabet generate fractaloids. Our Hecke-type operators, or labeling
operators, come from an amalgamated free probability construction, and we
compute the corresponding amalgamated free moments. We show that the free
moments are completely determined by certain scalar-valued functions.
\end{abstract}

\strut

A number of recent papers have addressed an intriguing interplay between
Discrete Potential Theory\ \strut \strut on the one hand and Harmonic
Analysis / Spectral Theory on an associated family of fractals on the other,
with the Sierpinski gasket serving as a preferred model.

\strut

This is the first of two papers studying representations of graphs by
operators in Hilbert space, and their applications. While graph theory is
traditionally considered part of discrete mathematics, in this paper we show
that applications of tools from automata and operators on Hilbert spaces
yield global results for representations of a class of infinite graphs, as
well as spin-off applications. We begin with an outline of the use of
automata, and more generally, of finite state models (FSMs) in the
processing of numbers, or more importantly in sampling and in quantization
of digitized information such as speech signals and digital images. In these
models, the finite input states of a particular FSM might be frequency-bands
(for example a prescribed pair of high-pass and low-pass digital filters),
or a choice of subdivision filters; where the subdivision refers to data
sets with self-similarity; such as is typically seen in fractals. Hence,
these applications make connections to discrete wavelet algorithms as used
in signal and image processing, as well as in science and engineering. If
the input-options for a particular FSM are chosen from a prescribed system
of low-pass and high-pass filters, the resulting discrete model can then be
realized by operators on Hilbert spaces. Similarly, images are digitized
into matrix shapes used in computer programs for compression of images.

\strut

In a general operator theoretic setting, this paper introduces the relevant
representations of the generators of graphs and automata. Hence data from
FSMs and graphs are represented with the use of Hilbert space geometry.
Recall that digital images are typically given by matrices of pixels, and
that spectral analysis and matrix operations can be done with operators in
Hilbert space. However, the Hilbert spaces needed in a particular
application are typically not immediately apparent from the particular
engineering problem under discussion. This paper focused on making the link
between graphs and automata on the discrete side to Hilbert space operators
and representations on the spectral side.

\strut \strut \strut \strut \strut \strut

The potential-theoretic part of the subject may be understood as a
mathematical idealization of electrical networks on infinite graphs $G$ (See
[41], [42] and [43]). We address two issues from Analysis: Find
representations (by operators on Hilbert spaces) of the graph systems $G,$
and identify a class of operators whose spectral theory captures significant
information about $G.$ We focus on the graphs themselves, and our motivation
derives in part from work by Strichartz and Kigami (See [44], [45] and [46])
and others for restricted classes of fractals. The focus there is the
adaptation of a rescaling and energy-renormalization on graphs and an
adaptation to fractal models $X$; for example, to Brownian motion on $X,$ or
to a version of differential operators on spaces of functions on $X.$ Here,
our focus is on Operator Theory needed for analysis of graphs and automata
(finite-state machines: e.g., [52]) such as are used in for example signal
processing algorithm, e.g., those based on a discrete multiresolution (e.g.,
see [50]).

\strut

We begin our paper with an outline of three trends: (i) Hilbert space and
Spectral Analysis on the graphs, (ii) Analysis on associated fractals
arising from automata, and finally (iii) the interplay between (i) and (ii).
We futher stress some of the differences between the two.

\strut

\strut The main purpose of this paper is to introduce a new algebraic
structures having certain fractal property. In [10], [11], [13], [14] and
[15], we introduced graph groupoids induced by countable directed graphs. A
graph groupoid is a categorial grouopoid having as a base the set of all
vertices of the given graph. We know that every groupoid having only one
base element is a group. So, if $G$ is a finite directed graph with graph
groupoid $\Bbb{G},$ and if the vertex set $V(G)$ $\subset $ $\Bbb{G}$
consists of only one element, then the graph groupoid $\Bbb{G}$ of this
one-vertex-$n$-loop-edge graph $G$, for $n$ $=$ $\left| E(G)\right| $ $\in $ 
$\Bbb{N},$ is a group; futhermore, it is group-isomorphic to the free group $%
F_{n}$ with $n$-generators (See [10] and [11]). Here, $E(G)$ denotes the
edge set of $G.$ We show that the automaton $\mathcal{A}_{G}$ $=$ $<X,$ $%
\Bbb{F}^{+}(G^{\symbol{94}}),$ $\varphi ,$ $\psi >$, having its alphabet $X$ 
$=$ $\{1,$ ..., $2n\},$ generates a fractal group $\Gamma $ $\overset{\text{%
Group}}{=}$ $F_{n}$. How about the general case where the vertex set
contains more than one element? This question is the motivation for our
paper. Also, if such generalization is determined suitably, then what are
the properties of the corresponding Hecke-type operator (with respect to the
natural representation of $\Bbb{G}$)? We will provide the answers of these
questions in the paper.

\strut

We realize that, by giving weights on the edges of a given graph $G$, which
are determined by the out-degrees of vertices, we can get the corresponding
automaton $\mathcal{A}_{G}$. By observing the properties of $\mathcal{A}%
_{G}, $ we can determine the groupoidal version of a fractal property,
relative to a fractal group. We will say that a graph groupoid $\Bbb{G}$ is
a fractaloid, if $\mathcal{A}_{G}$ acts ``fully'' on a $\left| V(G^{\symbol{%
94}})\right| $-copies of regular growing directed trees. For example, if a
directed graph $G$ is a one-flow circulant graph, then the graph groupoid $%
\Bbb{G}$ is a fractaloid.

\strut \strut \strut

Similar to the construction in [10], we can have a right graph von Neumann
algebra $\Bbb{M}_{G}$ $=$ $\overline{\Bbb{C}[\beta (\Bbb{G})]}^{w},$ as a $%
W^{*}$-subalgebra of the operator algebra $B(H_{G}),$ consisting of all
bounded linear operators on $H_{G},$ where $H_{G}$ is a generalized Fock
space, called the graph Hilbert space induced by $\Bbb{G}.$ We are
interested in the case where $\Bbb{G}$ is a fractaloid. Then, similar to the
classical case, we can define the Hecke-type operator $\tau $ $\in $ $\Bbb{M}%
_{G}.$ Instead of considering the pure operator-theoretical data of $\tau ,$
we will observe the amalgamated free distributional data of $\tau .$ Since $%
\tau $ is self-adjoint, the amalgamated free moments of it contains the
operator-valued spectral measure theoretical information. This means that
the free moments of $\tau $ will contain the operator-theoretical data of $%
\tau .$ Consequently, in Section 5, we can see that the free moments of $%
\tau $ is totally depending on the scalar-values, recursively, whenever $%
\Bbb{G}$ is a fractaloid.

\strut \strut

From the theory of algebras of operators on Hilbert space, we will need von
Neumann algebra constructions (e.g., [59]), free probability, in particular,
amalgamated free products (e.g., [5] and [21]), groupoids and groupoid
actions (e.g., [17] and [58]), and Hecke operators (e.g., [25]). In
technical discussions, we will use ``von Neumann algebras'' and ``$W^{*}$%
-algebras'' synonynously. If $H$ is a given Hilbert space, the algebra of
all bounded operators on $H$ is denoted by $B(H).$ From Graph Theory, we use
such notations as ``sets of vertices'', ``sets of edges'', ``loops'' and
``degrees'', etc. If $G$ is a given directed graph, we introduce a
``shadow'' construction from a reversal of edges, denoted by $G^{-1}$; see
details below. From symbolic dynamics, we shall use fundamentals of
``automata'', as well as free constructions, such as the free group with
multi-generators. Hence, our next section will contain a number of
definitions that will needed later.

\strut

Since our paper is interdisciplinary and directed at several audiences, we
have included details from one area of mathematics which might not be
familiar to readers from another.

\strut \strut

To get to our main theorems, a number of technical terms must be introduced.
We have collected them in Section 2 after Introduction. The essential ones
are: graphs, groupoids, and automata. In Section 2.2, we explain the
interplay between automata and the graphs, and we introduce ``canonical
weighted graphs'', ``shadowed graphs'', ``fractaloids'', and ``labeling
operators''. The first of the results on automata are in Section 3.3. From
Section 2.3 to Section 2.5, we review facts from Free Probability and von
Neumann algebras we will need. Section 5.5 introduces labeling operators,
and in Corollary 5.2, we show that these operators are generated by edges of
shadowed graphs.

\strut

\strut

\textbf{Table of Contents}

\strut

$
\begin{array}{lll}
1.\ \text{{\small Introduction}} &  & \;\text{{\small 6}} \\ 
2.\ \text{{\small Preliminaries}} &  & \ \text{{\small 9}} \\ 
\quad 2.1.\ \text{{\small Groupoids and Groupoid Actions}} &  & \ \text{%
{\small 9}} \\ 
\quad 2.2.\ \text{{\small Automata and Fractal Groups}} &  & \text{{\small 11%
}} \\ 
\quad 2.3.\ \text{{\small Free Probability}} &  & \text{{\small 15}} \\ 
\quad 2.4.\ \text{{\small Graph Groupoids and Right Graph von Neumann
Algebras}} &  & \text{{\small 18}} \\ 
\quad 2.5.\ M\text{{\small -Valued Graph}}\ W^{*}\text{{\small -Probability
Spaces}} &  & \text{{\small 22}} \\ 
3.\ \text{{\small Labeled Graph Groupoids and Graph Automata}} &  & \text{%
{\small 27}} \\ 
\quad 3.1.\ \text{{\small Labeled Graph Groupoids}} &  & \text{{\small 28}}
\\ 
\quad 3.2.\ \text{{\small The Operation}}\ \theta \ \text{{\small and}}\
\omega _{+} &  & \text{{\small 34}} \\ 
\quad 3.3.\ \text{{\small Graph Automata}} &  & \text{{\small 35}} \\ 
4.\ \text{{\small Fractaloids}} &  & \text{{\small 38}} \\ 
\quad 4.1.\ \text{{\small Fractaloids}} &  & \text{{\small 38}} \\ 
\quad 4.2.\ \text{{\small Examples}} &  & \text{{\small 45}} \\ 
5.\ \text{{\small Labeling Operators of Fractaloids}} &  & \text{{\small 47}}
\\ 
\quad 5.1.\ \text{{\small Labeling Operators}} &  & \text{{\small 47}} \\ 
\quad 5.2.\ \text{{\small Amalgamated Free Distributional Data of Labeling
Operators}} &  & \text{{\small 50}} \\ 
\quad 5.3.\ \text{{\small Labeling Operators of Fractaloids}} &  & \text{%
{\small 51}} \\ 
\quad 5.4.\ \text{{\small Refinements of (5.3)}} &  & \text{{\small 55}} \\ 
\quad 5.5.\ \text{{\small An Example}} &  & \text{{\small 62}}
\end{array}
$

\strut

\strut

Our main results are Theorems 4.2, 4.4, and Corollary 4.5, 5.5, 5.11, 5.14,
5.15. A number of examples and applications are computed at the end of
Section 5. Theorem 4.2 shows that the pair consisting a weighted graph and a
labeled graph groupoid induce a tree (graph), and we display its properties.
Futher, we give a necessary and sufficient condition for this to be a
fractaloid. In Theorem 4.4 and Corollary 4.5, we show that the automata
actions coming from fractaloids act fully on a certain tree. Moreover, they
are identified in terms of roots and out-degrees. Theorem 5.7 derives the
essential properties of the labeling operators. In Section 5.4, for a fixed
graph $G$, we compute the amalgamated free moments of the labeling operator $%
T_{G}$; see Corollaries 5.9, 5.11 5.14, and Theorem 5.13.

\strut \strut \strut \strut \strut \strut

A \emph{graph} is a set of objects called \emph{vertices} (or points or
nodes) connected by links called \emph{edges} (or lines). In a \emph{%
directed graph}, the two directions are counted as being distinct directed
edges (or arcs). A graph is depicted in a diagrammatic form as a set of dots
(for vertices), jointed by curves (for edges). Similarly, a directed graph
is depicted in a diagrammatic form as a set of dots jointed by arrowed
curves, where the arrows point the direction of the directed edges.

\strut

Throughout this paper, every graph is a locally finite countably directed
graph. Recall that we say that a countably directed graph $G$ is \emph{%
locally finite} if each vertex of $G$ has only finitely many incident edges.
Equivalently, the degree of $v$ is finitely determined. Also, recall that
the degree $\deg (v)$ of a vertex $v$ is defined to be the sum of the
out-degree $\deg _{out}(v)$ and the in-degree $\deg _{in}(v),$ where

\strut

\begin{center}
$\deg _{out}(v)$ $\overset{def}{=}$ $\left| \{e\in E(G):e\text{ has its
initial vertex }v\}\right| $
\end{center}

and

\begin{center}
$\deg _{in}(v)$ $\overset{def}{=}$ $\left| \{e\in E(G):e\text{ has its
terminal vertex }v\}\right| .$
\end{center}

\strut

Let $G^{\symbol{94}}$ be the shadowed graph of $G$, in the sense of [10].
Then we can consider the degree of each vertex of $G^{\symbol{94}},$ too,
since $G^{\symbol{94}}$ is also a locally finite countable directed graph.
Assume that

$\strut $

\begin{center}
$N$ $=$ $\max $ $\{\deg _{out}(v)$ $:$ $v$ $\in $ $V(G^{\symbol{94}})$ $=$ $%
V(G)\}.$
\end{center}

\strut

Then we can define the labeling set $X$ $=$ $\{1,$ $2,$ ..., $N\}.$ We
assign the weights $\{1,$ ..., $N\}$ to all elements in the edge set $E(G^{%
\symbol{94}})$ of the shadowed graph $G^{\symbol{94}}$ of $G.$ (This
weighting provides the weights of all elements in $\Bbb{G}$, which are the
sequences contained in $X_{0}^{\,\infty },$ where $X_{0}$ $=$ $\{0\}$ $\cup $
$X.$) We will call this process placing the weights onto all elements of $%
\Bbb{G},$ the labeling process. This labeling process lets us construct the
automaton $\mathcal{A}_{G}$ $=$ $<X_{0},$ $\Bbb{G},$ $\varphi ,$ $\psi $ $>$
induced by the graph groupoid $\Bbb{G}.$ If the automaton $\mathcal{A}_{G}$
satisfies certain \emph{fractal} property; we will call the graph groupoid $%
\Bbb{G}$ a fractaloid. Clearly, the word ``fractaloid'' hints at ``fractal
(graph) groupoid''.

\strut \strut \strut

As in [10] and [11], we will fix a representation $(H_{G},$ $\beta )$ of a
graph groupoid $\Bbb{G},$ where $H_{G}$ is the graph Hilbert space and $%
\beta $ is a certain groupoid action of $\Bbb{G}.$ Let $X$ $=$ $\{x_{1},$ $%
x_{2},$ ..., $x_{N}\}$ be the labeling set. Then we can define the operator $%
\tau _{j}$ $\in $ $B(H_{G})$ by

\strut

\begin{center}
$\tau _{j}(\xi _{w})$ $=$ $\xi _{w}$ $\xi _{e}$ $=$ $\xi _{w\,e},$ for all $%
\xi _{w}$ $\in $ $\mathcal{B}_{H_{G}},$
\end{center}

\strut

for each $j$ $\in $ $X,$ where $\mathcal{B}_{H_{G}}$ is the Hilbert basis of 
$H_{G},$ whenever an edge $e$ has its weight $x_{j},$ and $w$ $e$ $\neq $ $%
\emptyset $ in $\Bbb{G}.$ Then we can have the operator $\tau $ $\in $ $%
B(H_{G})$ defined by

\strut

\begin{center}
$\tau $ $=$ $\sum_{j=1}^{N}$ $\tau _{j}.$
\end{center}

\strut

This Hecke-type operator on $H_{G}$ is said to be the labeling operator of $%
\Bbb{G}$ on $H_{G}.$ We will consider the free distributional data of this
operator $\tau $ on $H_{G}.$ If $\Bbb{G}$ is a fractaloid, then $\tau $ is
self-adjoint. So, the free moments of it contain the spectral theoretical
properties of $\tau .$ In particular, we can show that the amalgamated free
moments $\left( E\left( \tau ^{n}\right) \right) _{n=1}^{\infty }$ of $\tau $
are determined by the cardinalities $\left( \eta _{n}\right) _{n=1}^{\infty
} $ of certain sets $\sum_{n}^{(N)},$ determined recursively for $n$ $\in $ $%
\Bbb{N}$. i.e., we can get that

\strut

\begin{center}
$E(\tau ^{n})$ $=$ $\eta _{n}$ $\cdot $ $1_{\Bbb{D}_{G}},$ for all $n$ $\in $
$\Bbb{N}.$
\end{center}

\strut

We provide the complete computation of $(\eta _{n})_{n=1}^{\infty }.$

\strut

\strut \strut \strut

\section{Introduction}

\strut

\strut \strut

\strut Recently, countable directed graphs have been studied in Pure and
Applied Mathematics, because not only are they important noncommutative
structures but they help us visualize such structures. Futhermore, the
visualization has nice matricial expressions, (sometimes, the
operator-valued matricial expressions dependent on) adjacency matrices or
incidence matrices of the given graph. In particular, the partial isometries
in an operator algebra can be expressed and visualized by directed graphs.
The main purpose of this paper is to introduce algebraic and
operator-algebraic structures induced by countable directed graphs. For
convenience, as we mentioned at the beginning of this paper, we will
restrict our interests to the case where the countable directed graphs are
locally finite. In [10] and [11], starting with a countable directed graph,
we assign certain algebraic elements gotten from the admissibility, and then
we assign partial isometries to those elements. From these operators, we
generated a von Neumann algebra and we then considered free probabilistic
properties of them (See [10]). In conclusion, we found the nice (free) block
structures of such von Neumann algebras and this provides the full
characterization of the von Neumann algebras (See [11]).

\strut

For a given countable directed graph $G,$ we can define the shadow $G^{-1}$
of $G,$ which is the oppositely directed graph of $G.$ Then we can define
the shadowed graph $G^{\symbol{94}}$ $=$ $G$ $\cup $ $G^{-1}$ of $G$ as a
directed graph with its vertex set $V(G^{\symbol{94}})$ $=$ $V(G)$ $=$ $%
V(G^{-1})$ and its edge set $E(G^{\symbol{94}})$ $=$ $E(G)$ $\cup $ $%
E(G^{-1})$. Construct the free semigroupoid $\Bbb{F}^{+}(G^{\symbol{94}})$
of the shadowed graph $G^{\symbol{94}},$ as a collection of all vertices and
finite paths of $G^{\symbol{94}}$ with its binary operation called the
admissibility, where the admissibility is nothing but the
direction-depending connectedness of elements in $\Bbb{F}^{+}(G^{\symbol{94}%
})$. Notice that all finite paths in $\Bbb{F}^{+}(G^{\symbol{94}})$ are the
words in $V(G^{\symbol{94}})$ $\cup $ $E(G^{\symbol{94}}).$ By defining the
reduction (RR) on $\Bbb{F}^{+}(G^{\symbol{94}}),$ we can construct the graph
groupoid $\Bbb{G}.$ i.e., the graph groupoid $\Bbb{G}$ is a set of all
``reduced'' words in $E(G^{\symbol{94}}),$ with the inherited admissibility
on $\Bbb{F}^{+}(G^{\symbol{94}}).$ In fact, this graph groupoid $\Bbb{G}$ is
a categorial groupoid with its base $V(G^{\symbol{94}})$.

\strut \strut \strut \strut \strut

Let $G$ be a countable directed graph and let $\Bbb{G}$ be the corresponding
graph groupoid of $G.$ Then it is a categorial groupoid (See Section 2.1).
It is well-known that every groupoid, having only one base element, is a
group. So, we can say that a graph groupoid $\Bbb{G}$ containing only one
vertex is a group. We can naturally expect that such graph groupoid $\Bbb{G}$
is induced by the one-vertex-multi-loop-edge graph. As we observed in [10]
and [11], such a graph groupoid is indeed a group, group-isomorphic to the
free group $F_{n}$ with $n$-generators, where $n$ is the cardinality of
loop-edges of $G.$ Notice that the free group $F_{n}$ acts fully on the $%
(2n) $-regular tree $\mathcal{T}_{2n},$ via the automata action of the
corresponding automaton $<X_{2n},$ $F_{n},$ $\varphi ,$ $\psi >,$ where $%
X_{n}$ $=$ $\{1,$ ..., $2n\}$ and $\varphi $ and $\psi $ are suitable maps
defined on $X_{2n}$ $\times $ $F_{n}.$ Moreover, the elements of $F_{n}$
acting on $\mathcal{T}_{2n}$ generate a fractal group. We will generalize
such fractal properties to the case where we have arbitrary groupoids.

\strut \strut \strut \strut \strut

For an arbitrary fixed von Neumann algebra $M$ in an operator algebra $B(K)$
of all bounded linear operators on a Hilbert space $K,$ we define a crossed
product algebra $\Bbb{M}_{G}$ $=$ $M$ $\times _{\beta }$ $\Bbb{G}$ of $M$
and $\Bbb{G}$ via a graph-representation $\beta $ $:$ $\Bbb{G}$ $\rightarrow 
$ $B(K$ $\otimes $ $H_{G}),$ where $H_{G}$ is the Hilbert space with its
Hilbert basis $\{\xi _{w}$ $:$ $w$ $\in $ $FP_{r}(G^{\symbol{94}})\},$ where

\strut

\begin{center}
$FP_{r}(G^{\symbol{94}})$ $\overset{def}{=}$ $\Bbb{G}$ $\setminus $ $\left(
V(G^{\symbol{94}})\text{ }\cup \text{ }\{\emptyset \}\right) .$
\end{center}

\strut

A graph-representation $\beta $ is a nonunital intertwined partial
representation or a groupoid action, determining the bounded operators $%
\beta _{w}$ on $K$ $\otimes $ $H_{G}$ satisfying that

\strut

\begin{center}
$\beta _{w}\left( m\right) $ $R_{w}$ $R_{w}^{*}=R_{w}^{*}$ $m$ $R_{w},$
\end{center}

and

\begin{center}
$\beta _{v}(m)$ $=$ $m,$
\end{center}

\strut \strut

for all \ $m$ $\in $ $M$,\ $w$ $\in $ $FP_{r}(G^{\symbol{94}})$ and $v$ $\in 
$ $V(G^{\symbol{94}}),$ where $R_{w}$ are the \emph{right} multiplication
operators on $H_{G}$ with their symbols $\xi _{w}$, for all $w$ $\in $ $\Bbb{%
G}.$ The adjoint $R_{w}^{*}$ of $R_{w}$ are defined to be $R_{w^{-1}},$ for
all $w$ $\in $ $\Bbb{G}.$

\strut

\begin{definition}
The crossed product algebra $\Bbb{M}_{G}$ $=$ $M$ $\times _{\beta }$ $\Bbb{G}
$ is the von Neumann algebra generated by $M$ and $\{R_{w}$ $:$ $w$ $\in $ $%
\Bbb{G}\},$ satisfying the above $\beta $-relations. We call $\Bbb{M}_{G}$ a
right graph von Neumann algebra induced by $G$ over $M$. A right graph von
Neumann algebra $\Bbb{M}_{G}$ has its canonical $W^{*}$-subalgebra $\Bbb{D}%
_{G}$ $\overset{def}{=}$ $\underset{v\in V(G^{\symbol{94}})}{\oplus }$ $(M$ $%
\cdot $ $R_{v}),$ called the $M$-diagonal subalgebra of $\Bbb{M}_{G}.$
\end{definition}

\strut

Recall that, in [10], we observed the (left) multiplication operators $L_{w}$%
's, for all $w$ $\in $ $\Bbb{G},$ instead of using right multiplication
operators $R_{w}$'s. Similar to [10], we construct an amalgamated $W^{*}$%
-probability space $(\Bbb{M}_{G},$ $E)$ over $\Bbb{D}_{G},$ where $E$ $:$ $%
\Bbb{M}_{G}$ $\rightarrow $ $\Bbb{D}_{G}$ is the canonical conditional
expectation. Under this setting, we realized that $\Bbb{M}_{G}$ is $*$%
-isomorphic to a $\Bbb{D}_{G}$-valued reduced free product $\underset{e\in
E(G)}{*_{\Bbb{D}_{G}}^{r}}$ $\Bbb{M}_{e}$ of the $\Bbb{D}_{G}$-free blocks $%
\Bbb{M}_{e}$ indexed by directed edges $e$ $\in $ $E(G),$ where

\strut

\begin{center}
$\Bbb{M}_{e}$ $\overset{def}{=}$ $vN\left( M\text{ }\times _{\beta }\text{ }%
\Bbb{G}_{e},\text{ }\Bbb{D}_{G}\right) ,$ for all $e$ $\in $ $E(G).$
\end{center}

\strut

Here, $\Bbb{G}_{e}$ is a substructure (or a subgroupoid) of $\Bbb{G},$
consisting of all reduced words in $\{e,$ $e^{-1}\},$ with the inherited
admissibility on $\Bbb{G},$ for all $e$ $\in $ $E(G).$ Notice that the
reduction of the free product ``$*_{\Bbb{D}_{G}}^{r}$'' is completely
dependent upon the admissibility on the graph groupoid $\Bbb{G}$.

\strut

In [11], we characterized $\Bbb{D}_{G}$-free blocks of graph von Neumann
algebras. Similar to [11], we can characterize the $\Bbb{D}_{G}$-free blocks 
$\Bbb{M}_{e}$'s of our right graph von Neumann algebra $\Bbb{M}_{G}.$
Because of the setting, in fact, the results are almost same. Especially, we
show that each $\Bbb{D}_{G}$-free block $\Bbb{M}_{e}$ of $\Bbb{M}_{G}$ is $*$%
-isomorphic to a certain von Neumann algebra contained in either $\mathcal{M}%
_{lp}^{\beta }$ or $\mathcal{M}_{non-lp}^{\beta },$ where

$\strut $

\begin{center}
$\mathcal{M}_{lp}^{\beta }$ $=$ $\{vN(M$ $\times _{\lambda _{e}}$ $\Bbb{Z},$ 
$\Bbb{D}_{G})$ $:$ $\lambda _{e}$ $=$ $\beta \mid _{\Bbb{S}_{e}},$ $e$ is a
loop edge$\}$
\end{center}

and

\begin{center}
$\mathcal{M}_{non-lp}^{\beta }$ $=$ $\{vN(M_{2}^{\beta _{e}}(M),$ $\Bbb{D}%
_{G})$ $:$ $e$ is a non-loop edge$\}$,
\end{center}

\strut

where $M$ $\times _{\lambda ^{(e)}}$ $\Bbb{Z}$ is a classical group crossed
products of $M$ and the infinite cyclic abelian group $\Bbb{Z}$ via a group
action $\lambda ^{(e)}$ satisfying that $\lambda _{(e)}$ $=$ $\beta $ $\mid
_{\Bbb{G}_{e}},$ and where $M_{2}^{\beta _{e}}(M)$ is a $W^{*}$-subalgebra
of $M_{2}(M)$ $=$ $M$ $\otimes _{\Bbb{C}}$ $M_{2}(\Bbb{C})$ satisfying the $%
\beta _{e}$-relation: $\beta _{e}(m)$ $R_{e}$ $R_{e}^{*}$ $=$ $R_{e}^{*}$ $m$
$R_{e}$, for all $m$ $\in $ $M$ and for $e$ $\in $ $E(G),$ where $M_{2}(\Bbb{%
C})$ is the matricial algebra generated by all $(2$ $\times $ $2)$-matrices.
In particular, if $M$ $=$ $\Bbb{C},$ we can conclude that each $\Bbb{D}_{G}$%
-free block $\Bbb{M}_{e}$ of $\Bbb{M}_{G}$ is $*$-isomorphic to either $%
L^{\infty }(\Bbb{T})$ or $M_{2}(\Bbb{C}),$ where $\Bbb{T}$ is the unit
circle in $\Bbb{C}$, for all $e$ $\in $ $E(G).$ This characterization says
that the study of right graph von Neumann algebras is the investigation of
graph groupoids and the above two types of von Neumann algebras.\strut

\strut

In this paper, we only consider the case where the fixed von Neumann algebra 
$M$ is $\Bbb{C}.$ Then the right graph von Neumann algebras $\Bbb{M}_{G}$ $=$
$\Bbb{C}$ $\times _{\beta }$ $\Bbb{G}$ are all $*$-isomorphic to $M_{G}$ $=$ 
$\overline{\Bbb{C}[\Bbb{G}]}^{w},$ as $W^{*}$-subalgebras of $B(H_{G}),$ for
all graph-representations $\beta ,$ by the linearity of $\beta $ on $\Bbb{C}$%
. We will say that the von Neumann algebra $M_{G}$ is \emph{the} right graph
von Neumann algebra of $G.$

\strut \strut

We will consider the labeling operator $\tau $ of $\Bbb{G}$, which is an
analogue of classical Hecke-type operators (of groups). As an element of the
right graph von Neumann algebra $M_{G}$, we can verify that $\tau $ has its
decomposition $\tau $ $=$ $\oplus _{j=1}^{N}$ $\tau _{j},$ where

\strut

\begin{center}
$N$ $\overset{def}{=}$ $\max \{\deg _{out}(v)$ $:$ $v$ $\in $ $V(G^{\symbol{%
94}})\}$
\end{center}

\strut and

\begin{center}
$\tau _{j}$ $(\xi _{w})$ $\overset{def}{=}$ $\left\{ 
\begin{array}{ll}
\xi _{w}\text{ }\xi _{e}=\xi _{w\,e} & 
\begin{array}{l}
\,\,\text{if }w\,e\neq \emptyset \text{ and} \\ 
e\text{ has its weight }j
\end{array}
\\ 
&  \\ 
\xi _{w}\text{ }\xi _{\emptyset }=\xi _{\emptyset }=0 & 
\begin{array}{l}
\text{if there is no edge }e\text{ such that} \\ 
e\text{ has its weight }j\text{ and }w\,e\neq \emptyset .
\end{array}
\end{array}
\right. $
\end{center}

\strut

The study of the labeling operator $\tau ,$ itself, is interesting in
Operator Theory and Quantum Physics, but we will concentrate on observing
its free probabilistic data, where $\Bbb{G}$ is a fractaloid. Consequently,
we can show that, if a graph groupoid $\Bbb{G}$ of a locally finite directed
graph $G$ is a fractaloid, then $\tau $ is self-adjoint and hence the free
moments of $\tau $ contain the Spectral Theoretical data of $\tau .$
Moreover, we can show that the free moments of $\tau $ are totally depending
on certain scalar-values.

\strut

\strut

\strut

\section{Preliminaries}

\strut

\strut

In this Section, we introduce the definitions and concepts which we will use
in this paper. We will review groupoids, groupoid actions, automata and
fractal groups. Also, we introduce and discuss right graph von Neumann
algebras.\strut

\strut

\strut

\subsection{Groupoids and Groupoid Actions}

\strut

\strut

\strut

While (categorial) groupoids and their actions are used in many areas of
mathematics (e.g., [17], [27], [34], and [58]), here we use them in
connections with graphs and representations, and we open with the necessary
definitions.

\strut

We say an algebraic structure $(\mathcal{X},$ $\mathcal{Y},$ $s,$ $r)$ is a
(categorial) groupoid if it satisfies that (i) $\mathcal{Y}$ $\subset $ $%
\mathcal{X},$ (ii) for all $x_{1},$ $x_{2}$ $\in $ $\mathcal{X},$ there
exists a partially-defined binary operation $(x_{1},$ $x_{2})$ $\mapsto $ $%
x_{1}$ $x_{2},$ for all $x_{1},$ $x_{2}$ $\in $ $\mathcal{X},$ depending on
the source map $s$ and the range map $r$ satisfying that:

\strut

(ii-1) $x_{1}$ $x_{2}$ is well-determined, whenever $r(x_{1})$ $=$ $s(x_{2})$
and in this case, $s(x_{1}$ $x_{2})$ $=$ $s(x_{1})$ and $r(x_{1}$ $x_{2})$ $%
= $ $r(x_{2}),$ for $x_{1},$ $x_{2}$ $\in $ $\mathcal{X},$

\strut

(ii-2) $(x_{1}$ $x_{2})$ $x_{3}$ $=$ $x_{1}$ $(x_{2}$ $x_{3})$, if they are
well-determined in the sense of (ii-1), for $x_{1},$ $x_{2},$ $x_{3}$ $\in $ 
$\mathcal{X},$

\strut

(ii-3) if $x$ $\in $ $\mathcal{X},$ then there exist $y,$ $y^{\prime }$ $\in 
$ $\mathcal{Y}$ such that $s(x)$ $=$ $y$ and $r(x)$ $=$ $y^{\prime },$
satisfying $x$ $=$ $y$ $x$ $y^{\prime }$ (Here, the elements $y$ and $%
y^{\prime }$ are not necessarily distinct),

\strut

(ii-4) if $x$ $\in $ $\mathcal{X},$ then there exists a unique element $%
x^{-1}$ for $x$ satisfying $x$ $x^{-1}$ $=$ $s(x)$ and $x^{-1}$ $x$ $=$ $%
r(x).$

\strut

Thus, every group is a groupoid $(\mathcal{X},$ $\mathcal{Y},$ $s,$ $r)$
with $\left| \mathcal{Y}\right| $ $=$ $1$ (and hence $s$ $=$ $r$ on $%
\mathcal{X}$). This subset $\mathcal{Y}$ of $\mathcal{X}$ is said to be the
base of $\mathcal{X}$. Remark that we can naturally assume that there exists
the empty element $\emptyset $ in a groupoid $\mathcal{X}.$ The empty
element $\emptyset $ means the products $x_{1}$ $x_{2}$ are not
well-defined, for some $x_{1},$ $x_{2}$ $\in $ $\mathcal{X}.$ Notice that if 
$\left| \mathcal{Y}\right| $ $=$ $1$ (equivalently, if $\mathcal{X}$ is a
group), then the empty word $\emptyset $ is not contained in the groupoid $%
\mathcal{X}.$ However, in general, whenever $\left| \mathcal{Y}\right| $ $%
\geq $ $2,$ a groupoid $\mathcal{X}$ always contain the empty word. So, if
there is no confusion, we will naturally assume that the empty element $%
\emptyset $ is contained in $\mathcal{X}.$

\strut

It is easily checked that our graph groupoid $\Bbb{G}$ of a finite directed
graph $G$ is indeed a groupoid with its base $V(G^{\symbol{94}}).$ i.e.,
every graph groupoid $\Bbb{G}$ of a countable directed graph $G$ is a
groupoid $(\Bbb{G},$ $V(G^{\symbol{94}}),$ $s$, $r)$, where $s(w)$ $=$ $s(v$ 
$w)$ $=$ $v$ and $r(w)$ $=$ $r(w$ $v^{\prime })$ $=$ $v^{\prime },$ for all $%
w$ $=$ $v$ $w$ $v^{\prime }$ $\in $ $\Bbb{G}$ with $v,$ $v^{\prime }$ $\in $ 
$V(G^{\symbol{94}}).$ i.e., the vertex set $V(G^{\symbol{94}})$ $=$ $V(G)$
is a base of $\Bbb{G}.$

\strut

Let $\mathcal{X}_{k}$ $=$ $(\mathcal{X}_{k},$ $\mathcal{Y}_{k},$ $s_{k},$ $%
r_{k})$ be groupoids, for $k$ $=$ $1,$ $2.$ We say that a map $f$ $:$ $%
\mathcal{X}_{1}$ $\rightarrow $ $\mathcal{X}_{2}$ is a groupoid morphism if
(i) $f$ is a function, (ii) $f(\mathcal{Y}_{1})$ $\subseteq $ $\mathcal{Y}%
_{2},$ (iii) $s_{2}\left( f(x)\right) $ $=$ $f\left( s_{1}(x)\right) $ in $%
\mathcal{X}_{2},$ for all $x$ $\in $ $\mathcal{X}_{1}$, and (iv) $%
r_{2}\left( f(x)\right) $ $=$ $f\left( r_{1}(x)\right) $ in $\mathcal{X}%
_{2}, $ for all $x$ $\in $ $\mathcal{X}_{1}.$ If a groupoid morphism $f$ is
bijective, then we say that $f$ is a groupoid-isomorphism, and the groupoids 
$\mathcal{X}_{1}$ and $\mathcal{X}_{2}$ are said to be groupoid-isomorphic.

\strut

Notice that, if two countable directed graphs $G_{1}$ and $G_{2}$ are
graph-isomorphic, via a graph-isomorphism $g$ $:$ $G_{1}$ $\rightarrow $ $%
G_{2},$ in the sense that (i) $g$ is bijective from $V(G_{1})$ onto $%
V(G_{2}),$ (ii) $g$ is bijective from $E(G_{1})$ onto $E(G_{2}),$ (iii) $%
g(e) $ $=$ $g(v_{1}$ $e$ $v_{2})$ $=$ $g(v_{1})$ $g(e)$ $g(v_{2})$ in $%
E(G_{2}),$ for all $e$ $=$ $v_{1}$ $e$ $v_{2}$ $\in $ $E(G_{1}),$ with $%
v_{1},$ $v_{2}$ $\in $ $V(G_{1}),$ then the graph groupoids $\Bbb{G}_{1}$
and $\Bbb{G}_{2}$ are groupoid-isomorphic. More generally, if $G_{1}$ and $%
G_{2}$ have graph-isomorphic shadowed graphs $G_{1}^{\symbol{94}}$ and $%
G_{2}^{\symbol{94}},$ then $\Bbb{G}_{1}$ and $\Bbb{G}_{2}$ are
groupoid-isomorphic.

\strut

\begin{proposition}
Let $G_{1}$ and $G_{2}$ be countable directed graph with their graph
groupoids $\Bbb{G}_{1}$ and $\Bbb{G}_{2},$ respectively. If the shadowed
graphs $G_{1}^{\symbol{94}}$ and $G_{2}^{\symbol{94}}$ are graph-isomorphic,
then $\Bbb{G}_{1}$ and $\Bbb{G}_{2}$ are groupoid-isomorphic.
\end{proposition}

\strut

\begin{proof}
Suppose $G_{1}^{\symbol{94}}$ and $G_{2}^{\symbol{94}}$ are
graph-isomorphic, via a graph-isomorphism $g$ $:$ $G_{1}^{\symbol{94}}$ $%
\rightarrow $ $G_{2}^{\symbol{94}}.$ Then we can define the morphism $%
\varphi $ $:$ $\Bbb{G}_{1}$ $\rightarrow $ $\Bbb{G}_{2}$, defined by

\strut

\begin{center}
$\varphi (w)$ $\overset{def}{=}$ $\left\{ 
\begin{array}{ll}
g(w) & \text{if }w\in V(G_{1}^{\symbol{94}})\cup E(G_{1}^{\symbol{94}}) \\ 
g(e_{1})...g(e_{n}) & 
\begin{array}{l}
\text{if }w=e_{1}...e_{n}\in FP_{r}(G^{\symbol{94}})\text{, } \\ 
\text{with }e_{1},\text{ ..., }e_{n}\in E(G^{\symbol{94}}),\text{ for }n>1
\end{array}
\\ 
\emptyset _{2} & \text{if }w=\emptyset _{1},
\end{array}
\right. $
\end{center}

\strut

where $\emptyset _{1}$ and $\emptyset _{2}$ are the empty word in $\Bbb{G}%
_{1}$ and $\Bbb{G}_{2},$ respectively. We can easily check that $\varphi $
is a groupoid-isomorphism, preserving the admissibility on $\Bbb{G}_{1}$ to
that on $\Bbb{G}_{2}.$
\end{proof}

\strut \strut \strut

Let $\mathcal{X}$ $=$ $(\mathcal{X},$ $\mathcal{Y},$ $s,$ $r)$ be a
groupoid. We say that this groupoid $\mathcal{X}$ acts on a set $Y$ if there
exists a groupoid action $\pi $ of $\mathcal{X}$ such that $\pi (x)$ $:$ $Y$ 
$\rightarrow $ $Y$ is a well-determined function, for all $x$ $\in $ $%
\mathcal{X}.$ Sometimes, we call the set $Y,$ a $\mathcal{X}$-set. We are
interested in the case where a $\mathcal{X}$-set $Y$ is a Hilbert space. The
nicest example of a groupoid action acting on a Hilbert space is a
graph-representation defined and observed in [10].

\strut

\begin{definition}
Let $\mathcal{X}_{1}$ $\subset $ $\mathcal{X}_{2}$ be a subset, where $%
\mathcal{X}_{2}$ $=$ $(\mathcal{X}_{2},$ $\mathcal{Y}_{2},$ $s,$ $r)$ is a
groupoid, and assume that $\mathcal{X}_{1}$ $=$ $(\mathcal{X}_{1},$ $%
\mathcal{Y}_{1},$ $s,$ $r),$ itself, is a groupoid, where $\mathcal{Y}_{1}$ $%
=$ $\mathcal{X}_{2}$ $\cap $ $\mathcal{Y}_{2}.$ Then we say that the
groupoid $\mathcal{X}_{1}$ is a subgroupoid of $\mathcal{X}_{2}.$
\end{definition}

\strut \strut

Recall that we say that a countable directed graph $G_{1}$ is a
full-subgraph of a countable directed graph $G_{2},$ if

\strut

\begin{center}
$E(G_{1})$ $\subseteq $ $E(G_{2})$
\end{center}

and

\begin{center}
$V(G_{1})$ $=$ $\{v$ $\in $ $V(G_{1})$ $:$ $e$ $=$ $v$ $e$ or $e$ $=$ $e$ $%
v, $ $\forall $ $e$ $\in $ $E(G_{1})\}.$
\end{center}

\strut

Remark the difference between full-subgraphs and subgraphs: We say that $%
G_{1}^{\prime }$ is a subgraph of $G_{2},$ if

\strut

\begin{center}
$V(G_{1}^{\prime })$ $\subseteq $ $V(G_{2})$
\end{center}

and

\begin{center}
$E(G_{1}^{\prime })$ $=$ $\{e$ $\in $ $E(G_{2})$ $:$ $e$ $=$ $v_{1}$ $e$ $%
v_{2},$ for $v_{1},$ $v_{2}$ $\in $ $V(G_{1}^{\prime })\}.$
\end{center}

\strut

We can see that the graph groupoid $\Bbb{G}_{1}$ of $G_{1}$ is a subgroupoid
of the graph groupoid $\Bbb{G}_{2}$ of $G_{2},$ whenever $G_{1}$ is a
full-subgraph of $G_{2}.$

\strut

\begin{proposition}
Let $G_{1}$ and $G_{2}$ be countable directed graphs with their graph
groupoids $\Bbb{G}_{1}$ and $\Bbb{G}_{2},$ respectively. If $G_{1}$ is a
full-subgraph of $G_{2},$ then $\Bbb{G}_{1}$ is a subgroupoid of $\Bbb{G}%
_{2} $
\end{proposition}

\strut

\begin{proof}
Since $G_{1}$ is a full-subgraph of $G_{2},$ the shadowed graph $G_{1}^{%
\symbol{94}}$ of $G_{1}$ is also a full-subgraph of the shadowed graph $%
G_{2}^{\symbol{94}}$ of $G_{2}.$ Therefore, the graph groupoid $\Bbb{G}_{1}$
of $G_{1}$ is a subset of the graph groupoid $\Bbb{G}_{2}$ of $G_{2}.$
Clearly, $\Bbb{G}_{1},$ itself, is a graph groupoid, and hence it is a
groupoid. Therefore, $\Bbb{G}_{1}$ is a subgroupoid of $\Bbb{G}_{2}$
\end{proof}

\strut \strut \strut

\strut \strut

\strut

\subsection{Automata and Fractal Groups}

\strut \strut

\strut \strut

Automata Theory is the study of abstract machines, and we are using it in
the formulation given by von Neumann. It is related to the theory of formal
languages. In fact, automata may be thought of as the class of formal
languages they are able to recognize. In von Neumann's version, an automaton
is a finite state machine (FSM). i.e., a machine with input of symbols,
transitions through a series of states according to a transition function
(often expressed as a table). The transition function tells the automata
which state to go to next, given a current state and a current symbol. The
input is read sequentially, symbol by symbol, for example as a tape with a
word written on it, registered by the head of the automaton; the head moves
forward over the tape one symbol at a time. Once the input is depleted, the
automaton stops. Depending on the state in which the automaton stops, it is
said that the automaton either accepts or rejects the input. The set of all
the words accepted by the automaton is called the language of the automaton.
For the benefit for the readers, we offer the following references for the
relevant part of Automata Theory: [1], [33], [34], [35], [48] and [49].

\strut

Let the quadruple $\mathcal{A}$ $=$ $<D,$ $Q,$ $\varphi ,$ $\psi >$ be
given, where $D$ and $Q$ are sets and

\strut

\begin{center}
$\varphi $ $:$ $D$ $\times $ $Q$ $\rightarrow $ $Q$ \ \ \ and \ \ $\psi $ $:$
$D$ $\times $ $Q$ $\rightarrow $ $D$
\end{center}

\strut

are maps. We say that $D$ and $Q$ are the (finite) alphabet and the state
set of $\mathcal{A},$ respectively and we say that $\varphi $ and $\psi $
are the output function and the state transition function, respectively. In
this case, the quadruple $\mathcal{A}$ is called an automaton. If the map $%
\psi (\bullet ,$ $q)$ is bijective on $D,$ for any fixed $q$ $\in $ $Q,$
then we say that the automaton $A$ is invertible. Similarly, if the map $%
\varphi (x,$ $\bullet )$ is bijective on $Q,$ for any fixed $x$ $\in $ $D,$
then we say that the automaton $\mathcal{A}$ is reversible. If the automaton 
$\mathcal{A}$ is both invertible and reversible, then $\mathcal{A}$ is said
to be bi-reversible.

\strut

To help visualize the use of automata, a few concrete examples may help.
With some oversimplification, they may be drawn from the analysis and
synthesis of input / output models in Engineering, often referred to as
black box diagram: excitation variables, response variables, and
intermediate variables (e.g., see [52] and [53]).

\strut

In our presentation above, the $D$ (the chosen finite alphabet) often takes
different forms on the side of input $D_{i}$ and output $D_{o}.$ In popular
automata that models stimuli of organisms, the three sets input $D_{i},$
output $D_{o}$, and the state set $Q,$ could be as in the following
prtotypical three examples:

\strut

\begin{example}
Models stimuli of organisms:

\strut

$\qquad D_{i}$ $=$ $\{$positive stimulus, negative stimulus$\},$

$\qquad D_{o}$ $=$ $\{$reaction, no reaction$\},$

and

$\qquad Q$ $=$ $\left\{ 
\begin{array}{c}
\text{reaction to last positive stimulus,} \\ 
\text{no reaction to last positive stimulus}
\end{array}
\right\} .$
\end{example}

\strut

\begin{example}
In a control model for say a steering mechanism in a vehicle:

\strut

$\qquad D_{i}$ $=$ $\{$right, left$\},$

$\qquad D_{o}$ $=$ $\{$switch on, switch off$\}$ or $\{$lamp on, lamp off$\}$

and

$\qquad Q$ $=$ $\left\{ 
\begin{array}{c}
\text{right-turning direction signal} \\ 
\text{left-turning direction signal}
\end{array}
\right\} .$
\end{example}

\strut

\begin{example}
In a model for quantization in Signal Processing:

\strut

\begin{center}
$D$ $=$ $D_{i}$ $=$ $D_{o}$ $=$ $\left\{ 
\begin{array}{c}
\text{assignments from a bit alphabet,} \\ 
\text{with the bits referring to the value} \\ 
\text{of pulses-in and pulses-out in a} \\ 
\text{signal processing algorithm}
\end{array}
\right\} ,$
\end{center}

\strut

for example, on a discrete multiresolution (e.g., [50]), and

\strut

\begin{center}
$Q$ $=$ $
\begin{array}{l}
\text{a subset of the Cartesian product of} \\ 
\text{copies of }D,\text{ fixing finite number of} \\ 
\text{times, i.e., }D\text{ }\times \text{ ... }\times \text{ }D
\end{array}
$
\end{center}
\end{example}

\strut \strut

Recently, various algebraists have studied automata and the corresponding
automata groups (Also, see [1], [20], [33] and [35]). We will consider a
certain special case, where $Q$ is a free semigroupoid of a shadowed graph.

\strut

Roughly speaking, a undirected tree is a connected simplicial graph without
loop finite paths. Recall that a (undirected) graph is simplicial, if the
graph has neither loop-edges nor multi-edges. A directed tree is a connected
simplicial graph without loop finite paths, with directed edges (See Section
5.1, more in detail). In particular, we say that a directed tree $\mathcal{T}%
_{n}$ is a $n$-regular tree, if $\mathcal{T}_{n}$ is rooted and one-flowed
infinite directed tree, having the same out-degrees $n,$ for all vertices
(Also, see Section 5.1, for details). For example, the $2$-regular tree $%
\mathcal{T}_{2}$ can be depicted by

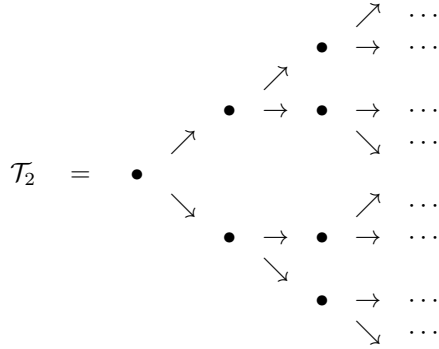
\strut

\begin{center}
$\mathcal{T}_{2}\quad =$\quad $
\begin{array}{lllllll}
&  &  &  &  & \nearrow & \cdots \\ 
&  &  &  & \bullet & \rightarrow & \cdots \\ 
&  &  & \nearrow &  &  &  \\ 
&  & \bullet & \rightarrow & \bullet & \rightarrow & \cdots \\ 
& \nearrow &  &  &  & \searrow & \cdots \\ 
\bullet &  &  &  &  &  &  \\ 
& \searrow &  &  &  & \nearrow & \cdots \\ 
&  & \bullet & \rightarrow & \bullet & \rightarrow & \cdots \\ 
&  &  & \searrow &  &  &  \\ 
&  &  &  & \bullet & \rightarrow & \cdots \\ 
&  &  &  &  & \searrow & \cdots
\end{array}
$
\end{center}

\strut

Let $\mathcal{A}$ $=$ $<D,$ $Q,$ $\varphi ,$ $\psi >$ be an automaton with $%
\left| D\right| $ $=$ $n.$ Then, we can construct automata actions of $%
\mathcal{A}$ on $\mathcal{T}_{n}.$ Let's fix $q$ $\in $ $Q.$ Then the action
of $\mathcal{A}_{q}$ is defined on the finite words $D_{*}$ of $D$ by

\strut

\begin{center}
$\mathcal{A}_{q}\left( x\right) $ $\overset{def}{=}$ $\varphi (x,$ $q),$ for
all $x$ $\in $ $D,$
\end{center}

\strut

and recursively,

$\strut $

\begin{center}
$\mathcal{A}_{q}\left( (x_{1},\text{ }x_{2},\text{ ..., }x_{m})\right) $ $=$ 
$\varphi \left( x_{1},\text{ }\mathcal{A}_{q}(x_{2},...,x_{m})\right) ,$
\end{center}

\strut

for all $(x_{1},$ ..., $x_{m})$ $\in $ $D_{*},$ where

\strut

\begin{center}
$D_{*}$ $\overset{def}{=}$ $\cup _{m=1}^{\infty }$ $\left( \left\{ (x_{1},%
\text{ ..., }x_{m})\in D^{m}\left| 
\begin{array}{c}
\text{ }x_{k}\in D,\text{ for all} \\ 
k=1,...,n
\end{array}
\right. \right\} \right) .$
\end{center}

\strut

Then the automata actions $\mathcal{A}_{q}$'s are acting on the $n$-regular
tree $\mathcal{T}_{n}$. In other words, all images of automata actions are
regarded as an elements in the free semigroupoid $\Bbb{F}^{+}(\mathcal{T}%
_{n})$ of the $n$-regular tree. i.e.,

\strut

\begin{center}
$V(\mathcal{T}_{n})$ $\supseteq $ $D_{*}$
\end{center}

\strut

and its edge set

\strut

\begin{center}
$\strut E(\mathcal{T}_{n})$ $\supseteq $ $\{\mathcal{A}_{q}(x)$ $:$ $x$ $\in 
$ $D,$ $q$ $\in $ $Q\}.$
\end{center}

\strut

This makes us to illustrate how the automata actions work.

\strut

Let $\mathcal{C}$ $=$ $\{\mathcal{A}_{q}$ $:$ $q$ $\in $ $Q\}$ be the
collection of automata actions of the given automaton $\mathcal{A}$ $=$ $<D,$
$Q,$ $\varphi ,$ $\psi >$. Then we can create a group $G(\mathcal{A})$
generated by the collection $\mathcal{C}.$ This group $G(\mathcal{A})$ is
called the automata group generated by $\mathcal{A}.$ The generator set $%
\mathcal{C}$ of $G(\mathcal{A})$ acts \emph{fully} on the $\left| D\right| $%
-regular tree $\mathcal{T}_{\left| D\right| },$ we say that this group $G(%
\mathcal{A})$ is a fractal group. There are many ways to define fractal
groups, but we define them in the sense of automata groups. (See [1] and
[35]. In fact, in [35], Batholdi, Grigorchuk and Nekrashevych did not define
the term ``fractal'', but they provide the fractal properties.)

\strut

Now, we will define a fractal group more precisely (Also see [1]). Let $%
\mathcal{A}$ be an automaton and let $\Gamma $ $=$ $G(\mathcal{A})$ be the
automata group generated by the automata actions acting on the $n$-regular
tree $\mathcal{T}_{n},$ where $n$ is the cardinality of the alphabet of $%
\mathcal{A}.$ By $St_{\Gamma }(k),$ denote the subgroup of $\Gamma $ $=$ $G(%
\mathcal{A})$, consisting of those elements of $\Gamma ,$ acting trivially
on the $k$-th level of $\mathcal{T}_{n},$ for all $k$ $\in $ $\Bbb{N}$ $\cup 
$ $\{0\}.$

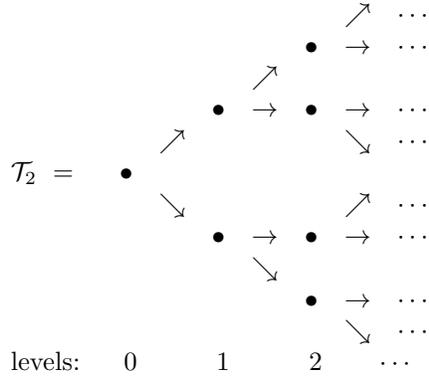
\strut

\begin{center}
$
\begin{array}{ll}
\mathcal{T}_{2}\text{ }= & 
\begin{array}{lllllll}
&  &  &  &  & \nearrow & \cdots \\ 
&  &  &  & \bullet & \rightarrow & \cdots \\ 
&  &  & \nearrow &  &  &  \\ 
&  & \bullet & \rightarrow & \bullet & \rightarrow & \cdots \\ 
& \nearrow &  &  &  & \searrow & \cdots \\ 
\bullet &  &  &  &  &  &  \\ 
& \searrow &  &  &  & \nearrow & \cdots \\ 
&  & \bullet & \rightarrow & \bullet & \rightarrow & \cdots \\ 
&  &  & \searrow &  &  &  \\ 
&  &  &  & \bullet & \rightarrow & \cdots \\ 
&  &  &  &  & \searrow & \cdots
\end{array}
\\ 
\text{levels:} & \,\,\,\,0\qquad \quad 1\qquad \quad 2\qquad \cdots
\end{array}
$
\end{center}

\strut

Analogously, for a vertex $u$ in $\mathcal{T}_{n},$ define $St_{\Gamma }(u)$
by the subgroup of $\Gamma ,$ consisting of those elements of $\Gamma ,$
acting trivially on $u.$ Then

$\strut $

\begin{center}
$St_{\Gamma }(k)$ $=$ $\underset{u\,:\,\text{vertices of the }k\text{-th
level of }\mathcal{T}_{n}}{\cap }$ $\left( St_{\Gamma }(u)\right) .$
\end{center}

\strut

For any vertex $u$ of $\mathcal{T}_{n},$ we can define the algebraic
projection $p_{u}$ $:$ $St_{\Gamma }(u)$ $\rightarrow $ $\Gamma .$

\strut

\begin{definition}
Let $\Gamma $ $=$ $G(\mathcal{A})$ be the automata group given as above. We
say that this group $\Gamma $ is a fractal group if, for any vertex $u$ of $%
\mathcal{T}_{n},$ the image of the projection $p_{u}\left( St_{\Gamma
}(u)\right) $ is group-isomorphic to $\Gamma ,$ after the identification of
the tree $\mathcal{T}_{n}$ with its subtree $\mathcal{T}_{u}$ with the root $%
u.$
\end{definition}

\strut

For instance, if $u$ is a vertex of the $2$-regular tree $\mathcal{T}_{2}$,
then we can construct a subtree $\mathcal{T}_{u},$ as follows:

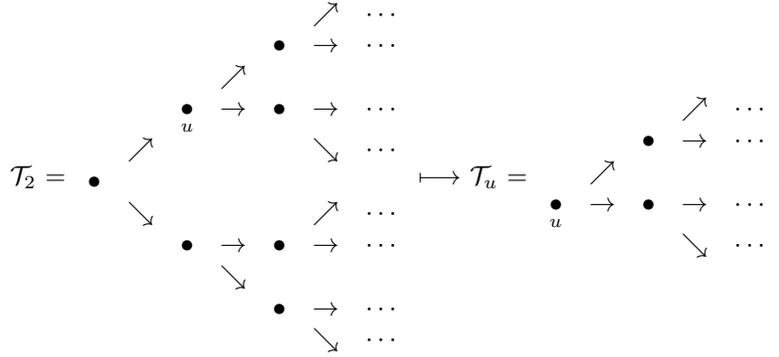
\strut

\begin{center}
$\mathcal{T}_{2}$ $=$ $
\begin{array}{lllllll}
&  &  &  &  & \nearrow & \cdots \\ 
&  &  &  & \bullet & \rightarrow & \cdots \\ 
&  &  & \nearrow &  &  &  \\ 
&  & \underset{u}{\bullet } & \rightarrow & \bullet & \rightarrow & \cdots
\\ 
& \nearrow &  &  &  & \searrow & \cdots \\ 
\bullet &  &  &  &  &  &  \\ 
& \searrow &  &  &  & \nearrow & \cdots \\ 
&  & \bullet & \rightarrow & \bullet & \rightarrow & \cdots \\ 
&  &  & \searrow &  &  &  \\ 
&  &  &  & \bullet & \rightarrow & \cdots \\ 
&  &  &  &  & \searrow & \cdots
\end{array}
$ $\longmapsto $ $\mathcal{T}_{u}$ $=$ $
\begin{array}{lllll}
&  &  & \nearrow & \cdots \\ 
&  & \bullet & \rightarrow & \cdots \\ 
& \nearrow &  &  &  \\ 
\underset{u}{\bullet } & \rightarrow & \bullet & \rightarrow & \cdots \\ 
&  &  & \searrow & \cdots
\end{array}
.$
\end{center}

\strut

As we can check, the graphs $\mathcal{T}_{2}$ and $\mathcal{T}_{u}$ are
graph-isomorphic. So, the above definition shows that if the automata
actions $\mathcal{A}_{q}$'s of $\mathcal{A}$ are acting \emph{fully} on $%
\mathcal{T}_{n},$ then the automata group $G(\mathcal{A})$ is a fractal
group. There are lots of famous fractal groups, but we introduce the
following example, for our purpose.

\strut

\begin{example}
Let $\mathcal{A}$ $=$ $<X_{2n},$ $F_{n},$ $\varphi ,$ $\psi >$ be an
automaton, where $F_{n}$ is the free group with its generator set $X_{2n}$ $%
= $ $\{g_{1}^{\pm 1},$ ..., $g_{n}^{\pm 1}\}.$ Then the automata group $G(%
\mathcal{A})$ is group-isomorphic to $F_{n}.$ It is easy to check that all
elements in $F_{n}$ acts fully on the $2n$-regular tree $\mathcal{T}_{2n},$
and hence $G(\mathcal{A})$ is a fractal group. For example, if $n$ $=$ $2,$
then we can get the following $0$-th and $1$-st levels of $\mathcal{T}_{4}$:

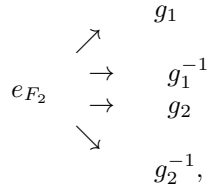
\strut

\begin{center}
$
\begin{array}{lll}
&  & g_{1} \\ 
& \nearrow &  \\ 
e_{F_{2}} & 
\begin{array}{l}
\rightarrow \\ 
\rightarrow
\end{array}
& 
\begin{array}{l}
g_{1}^{-1} \\ 
g_{2}
\end{array}
\\ 
& \searrow &  \\ 
&  & g_{2}^{-1},
\end{array}
$
\end{center}

\strut

where $e_{F_{2}}$ is the group-identity of $F_{2}$ $=$ $<g_{1},$ $g_{2}>.$
\end{example}

\strut \strut

\strut

\strut

\subsection{Free Probability}

\strut

\strut

Early versions of free probability models were motivated in part by finitely
correlated states which were built mathematically to realize quantum spin
chains; see e.g., [47], [54], [55], [56], and [57]; the purpose being to
analyze phase transition: Again large composite systems are built from
finite alphabets, with recursive algorithms run on finite-state machines,
with the finite data representing for example spin variables on a lattice,
in the simplest case, the lattice may be the set of integers with its usual
order, and the forward correlations depending on a finite number of states
going back in the state iterations. To make this precise, we consider a pair
of von Neumann algebras $B$ $\subset $ $A$, and with a prescribed
conditional expectation $E$ from $A$ onto $B$, mapping the unit of $A$ to
that of $B$. The free construction of states on infinite models built on the
pair $B$ $\subset $ $A$ is then done by repeating a prescribed automaton
which we now proceed to describe.

\strut

Let $B$ $\subset $ $A$ be von Neumann algebras with $1_{B}$ $=$ $1_{A}$ and
assume that there is a conditional expectation $E_{B}$ $:$ $A$ $\rightarrow $
$B$ satisfying that (i) $E_{B}$ is a ($\Bbb{C}$-)linear map, (ii) $E_{B}(b)$ 
$=$ $b,$ for all $b$ $\in $ $B$, (iii) $E_{B}(b_{1}$ $a$ $b_{2})$ $=$ $b_{1}$
$E_{B}(a)$ $b_{2},$ for all $b_{1},$ $b_{2}$ $\in $ $B$ and $a$ $\in $ $A,$
(iv) $E_{B}$ is continuous under the given topologies of $A$ and $B,$ and
(v) $E_{B}(a^{*})$ $=$ $E_{B}(a)^{*}$ in $B,$ for all $a$ $\in $ $A.$ The
algebraic pair $(A,$ $E_{B})$ is said to be a $B$-valued $W^{*}$-probability
space. Every operator in $(A,$ $E_{B})$ is called a $B$-valued (free) random
variable. Any $B$-valued random variables have their $B$-valued free
distributional data: $B$-valued $*$-moments and $B$-valued $*$-cumulants of
them. Suppose $a_{1},$ ..., $a_{s}$ are $B$-valued random variables in $(A,$ 
$E_{B}),$ where $s$ $\in $ $\Bbb{N}.$ The $(i_{1},$ ..., $i_{n})$-th joint $%
B $-valued $*$-moments of $a_{1},$ ..., $a_{s}$ are defined by

\strut

\begin{center}
$E_{B}\left(
(b_{1}a_{i_{1}}^{r_{i_{1}}})(b_{2}\,a_{i_{2}}^{r_{i_{2}}})\,...(b_{n}%
\,a_{i_{n}}^{r_{i_{n}}})\right) ,$
\end{center}

\strut

and the $(j_{1},$ ..., $j_{k})$-th joint $B$-valued $*$-cumulants of $a_{1},$
..., $a_{s}$ are defined by

\strut

\begin{center}
$k_{k}^{B}\left( (b_{1}a_{j_{1}}^{r_{i_{1}}}),\text{ ..., }%
(b_{k}a_{j_{k}}^{r_{j_{k}}})\right) =\underset{\pi \in NC(k)}{\sum_{n}^{(N)}}%
E_{B:\,\pi }\left( b_{1}a_{j_{1}}^{r_{j_{1}}}\text{ },...,\text{ }%
b_{k}a_{j_{k}}^{r_{j_{k}}}\right) \mu (\pi ,$ $1_{k}),$
\end{center}

\strut

for all $(i_{1},$ ..., $i_{n})$ $\in $ $\{1,$ ..., $s\}^{n}$ and for all $%
(j_{1},$ ..., $j_{k})$ $\in $ $\{1,$ ..., $s\}^{k},$ for $n,$ $k$ $\in $ $%
\Bbb{N},$ where $b_{j}$ $\in $ $B$ are arbitrary and $r_{i_{1}},$ ..., $%
r_{i_{n}},$ $r_{j_{1}},$ ..., $r_{j_{k}}$ $\in $ $\{1,$ $*\}$ and $NC(k)$ is
the lattice of all noncrossing partitions with its minimal element $0_{k}$ $%
= $ $\{(1),$ $(2),$ ..., $(k)\}$ and its maximal element $1_{k}$ $=$ $\{(1,$ 
$2,$ ..., $k)\},$ for all $k$ $\in $ $\Bbb{N},$ and $\mu $ is the Moebius
functional in the incidence algebra $\mathcal{I}.$ Here, $E_{B:\pi }(...)$
is the partition-depending $B$-valued moment. For example, if $\pi $ $=$ $%
\{(1,$ $4),$ $(2,$ $3),$ $(5)\}$ in $NC(5),$ then

\strut

\begin{center}
$E_{B:\pi }\left( a_{1},a_{2},a_{3},a_{4},a_{5}\right) =E_{B}\left(
a_{1}E_{B}(a_{2}a_{3})a_{4}\right) E_{B}(a_{5}).$
\end{center}

\strut \strut \strut \strut

Recall that the lattice $NC(n)$ of all noncrossing partitions over $\{1,$
..., $n\}$ has its partial ordering ``$\leq $'',

\strut

\begin{center}
$\pi \leq \theta $ $\overset{def}{\Longleftrightarrow }$ for each block $V$
in $\pi ,$ $\exists $ a block $B$ in $\theta $ s.t.,. $V$ $\subseteq $ $B,$
\end{center}

\strut

for $\pi ,\theta $ $\in $ $NC(n),$ where ``$\subseteq $'' means the usual
set inclusion, for all $n$ $\in $ $\Bbb{N}.$ Also recall that the incidence
algebra $\mathcal{I}$ is the collection of all functionals

\strut

\begin{center}
$\xi $ $:$ $\cup _{n=1}^{\infty }\left( NC(n)\times NC(n)\right) $ $%
\rightarrow $ $\Bbb{C}$
\end{center}

\strut

satisfying that $\xi (\pi ,$ $\theta )$ $=$ $0,$ whenever $\pi $ $>$ $\theta 
$, with the usual function addition $(+)$ and the convolution $(*)$ defined
by

\strut

\begin{center}
$\left( \xi _{1}*\xi _{2}\right) $ $(\pi ,$ $\theta )$ $\overset{def}{=}$ $%
\underset{\pi \leq \sigma \leq \theta }{\sum_{n}^{(N)}}$ $\xi _{1}(\pi
,\sigma )\xi _{2}(\sigma ,\theta ),$
\end{center}

\strut

for all $\xi _{1},$ $\xi _{2}$ $\in $ $\mathcal{I}.$ If we define the zeta
functional $\zeta $ $\in $ $\mathcal{I}$ by

\strut

\begin{center}
$\zeta (\pi ,$ $\theta )$ $=$ $1,$ \ for all \ $\pi $ $\leq $ $\theta $ in $%
NC(n),$ for all $n$ $\in $ $\Bbb{N},$
\end{center}

\strut

then its convolution-inverse in $\mathcal{I}$ is the Moebius functional $\mu 
$. Thus the Moebius functional $\mu $ satisfies that

\strut

\begin{center}
$\mu (0_{n},$ $1_{n})$ $=$ $(-1)^{n-1}$ $c_{n-1}$ \ and \ $\underset{\pi \in
NC(n)}{\sum_{n}^{(N)}}$ $\mu (\pi ,$ $1_{n})$ $=$ $0,$
\end{center}

\strut

where $c_{k}$ $\overset{def}{=}$ $\frac{1}{k+1}$ $\left( 
\begin{array}{l}
2k \\ 
\,\,k
\end{array}
\right) $ is the $k$-th Catalan number, for all $k$ $\in $ $\Bbb{N}$ (See
[21]).

\strut \strut

The $B$-valued freeness on $(A,$ $E_{B})$ is characterized by the $B$-valued 
$*$-cumulants. Let $A_{1}$ and $A_{2}$ be $W^{*}$-subalgebras of $A$ having
their common $W^{*}$-subalgebra $B.$ We say that $A_{1}$ and $A_{2}$ are
free over $B$ in $(A,$ $E_{B})$ if all mixed $B$-valued $*$-cumulants of $%
A_{1}$ and $A_{2}$ vanish. The subsets $X_{1}$ and $X_{2}$ of $A$ are said
to be free over $B$ in $(A,$ $E_{B})$ if the $W^{*}$-subalgebras $vN(X_{1},$ 
$B)$ and $vN(X_{2},$ $B)$ of $A$ are free over $B$ in $(A,$ $E_{B}),$ where $%
vN(S_{1},$ $S_{2})$ means the von Neumann algebra generated by sets $S_{1}$
and $S_{2}.$ Similarly, we say that the $B$-valued random variables $x$ and $%
y$ are free over $B$ in $(A,$ $E_{B})$ if the subsets $\{x\}$ and $\{y\}$
are free over $B$ in $(A,$ $E_{B})$ (Also See [21]).

\strut

Let $A_{1}$ and $A_{2}$ be $W^{*}$-subalgebra of $A$ containing their common 
$W^{*}$-subalgebra $B,$ and assume that they are free over $B$ in $(A,$ $%
E_{B}).$ Then we can construct a $W^{*}$-algebra $vN(A_{1},$ $A_{2})$ of $A$
generated by $A_{1}$ and $A_{2}.$ We denote it by $A_{1}$ $*_{B}$ $A_{2}.$
Suppose there exists a family $\{A_{i}$ $:$ $i$ $\in $ $\Lambda \}$ of $%
W^{*} $-subalgebras of $A$ containing their common $W^{*}$-subalgebra $B.$
If they generate $A$ and if they are free over $B$ from each other in $(A,$ $%
E_{B})$, then the von Neumann algebra $A$ is the $B$-valued free product
algebra $\underset{i\in \Lambda }{\,*_{B}}$ $A_{i}.$

\strut

Suppose a von Neumann algebra $A$ is a $B$-free product algebra $\underset{%
i\in \Lambda }{\,*_{B}}$ $A_{i}.$ Then $A$ is Banach-space isomorphic to the
Banach space

\strut

\begin{center}
$B$ $\oplus $ $\left( \oplus _{n=1}^{\infty }\left( \underset{i_{1}\neq
i_{2},\,i_{2}\neq i_{3},\,...,\,i_{n-1}\neq i_{n}}{\oplus }\left(
A_{i_{1}}^{o}\text{ }\otimes _{B}\text{ ... }\otimes _{B}\text{ }%
A_{i_{n}}^{o}\right) \right) \right) $
\end{center}

with

\begin{center}
$A_{i_{j}}^{o}$ $\overset{def}{=}$ $A_{i_{j}}$ $\ominus $ $B,$ for all $j$ $%
= $ $1,$ ..., $n,$
\end{center}

\strut

where $\otimes _{B}$ is the $B$-valued tensor product.\strut

\strut \strut \strut \strut \strut

\strut \strut \strut \strut

\strut \strut \strut \strut

\subsection{Graph Groupoids and Right Graph Von Neumann algebras}

\strut

\strut

\strut

Let $G$ be a countable directed graph with its vertex set $V(G)$ and its
edge set $E(G).$ Denote the set of all finite paths of $G$ by $FP(G).$
Clearly, the edge set $E(G)$ is contained in $FP(G).$ Let $w$ be a finite
path in $FP(G).$ Then it is represented as a word in $E(G),$ depending on
the directions of edges. If $e_{1},$ ..., $e_{n}$ are connected directed
edges in the order $e_{1}$ $\rightarrow $ $e_{2}$ $\rightarrow $ ... $%
\rightarrow $ $e_{n},$ for $n$ $\in $ $\Bbb{N},$ then we can express $w$ by $%
e_{1}$ ... $e_{n}$ in $FP(G).$ If there exists a finite path $w$ $=$ $e_{1}$
... $e_{n}$ in $FP(G),$ where $n$ $\in $ $\Bbb{N}$ $\setminus $ $\{1\},$ we
say that the directed edges $e_{1},$ ..., $e_{n}$ are admissible. The length 
$\left| w\right| $ of $w$ is defined to be $n,$ which is the cardinality of
the admissible edges generating $w.$ Also, we say that finite paths $w_{1}$ $%
=$ $e_{11}$ ... $e_{1k_{1}}$ and $w_{2}$ $=$ $e_{21}$ ... $e_{2k_{2}}$ are
admissible, if $w_{1}$ $w_{2}$ $=$ $e_{11}$ ... $e_{1k_{1}}$ $e_{21}$ ... $%
e_{2k_{2}}$ is again in $FP(G),$ where $e_{11},$ ..., $e_{1k_{1}},$ $e_{21},$
..., $e_{2k_{2}}$ $\in $ $E(G).$ Otherwise, we say that $w_{1}$ and $w_{2}$
are not admissible. Suppose that $w$ is a finite path in $FP(G)$ with its
initial vertex $v_{1}$ and its terminal vertex $v_{2}.$ Then we write $w$ $=$
$v_{1}$ $w$ or $w$ $=$ $w$ $v_{2}$ or $w$ $=$ $v_{1}$ $w$ $v_{2},$ for
emphasizing the initial vertex of $w,$ respectively the terminal vertex of $%
w,$ respectively both the initial vertex and the terminal vertex of $w$.
Suppose $w$ $=$ $v_{1}$ $w$ $v_{2}$ in $FP(G)$ with $v_{1},$ $v_{2}$ $\in $ $%
V(G).$ Then we also say that [$v_{1}$ and $w$ are admissible] and [$w$ and $%
v_{2}$ are admissible]. Notice that even though the elements $w_{1}$ and $%
w_{2}$ in $V(G)$ $\cup $ $FP(G)$ are admissible, $w_{2}$ and $w_{1}$ are not
admissible, in general. For instance, if $e_{1}$ $=$ $v_{1}$ $e_{1}$ $v_{2}$
is an edge with $v_{1}$, $v_{2}$ $\in $ $V(G)$ and $e_{2}$ $=$ $v_{2}$ $%
e_{2} $ $v_{3}$ is an edge with $v_{3}$ $\in $ $V(G)$ such that $v_{3}$ $%
\neq $ $v_{1},$ then there is a finite path $e_{1}$ $e_{2}$ in $FP(G),$ but
there is no finite path $e_{2}$ $e_{1}$.

\strut

The free semigroupoid $\Bbb{F}^{+}(G)$ of $G$ is defined by a set

\strut

\begin{center}
$\Bbb{F}^{+}(G)=\{\emptyset \}\cup V(G)\cup FP(G),$
\end{center}

\strut

with its binary operation $(\cdot )$ on $\Bbb{F}^{+}(G)$, defined by

\strut

\begin{center}
$(w_{1},w_{2})\mapsto w_{1}\cdot w_{2}=\left\{ 
\begin{array}{ll}
w_{1} & \text{if }w_{1}=w_{2}\text{ in }V(G) \\ 
w_{1} & \text{if }w_{1}\in FP(G),\text{ }w_{2}\in V(G)\text{ and }%
w_{1}=w_{1}w_{2} \\ 
w_{2} & \text{if }w_{1}\in V(G),\text{ }w_{2}\in FP(G)\text{ and }%
w_{2}=w_{1}w_{2} \\ 
w_{1}w_{2} & \text{if }w_{1},\text{ }w_{2}\text{ in }FP(G)\text{ and }%
w_{1}w_{2}\in FP(G) \\ 
\emptyset & \text{otherwise,}
\end{array}
\right. $
\end{center}

\strut

where $\emptyset $ is the empty word in $V(G)$ $\cup $ $E(G).$ (Sometimes,
the free semigroupoid $\Bbb{F}^{+}(G)$ of a certain graph $G$ does not
contain the empty word $\emptyset .$ For instance, the free semigroupoid of
the one-vertex--multi-loop-edge graph does not have the empty word. But, in
general, the empty word $\emptyset $ is contained in the free semigroupoid,
whenever $\left| V(G)\right| $ $\geq $ $2$. So, if there is no confusion,
then we usually assume that the empty word is contained in free
semigroupoids.) This binary operation $(\cdot )$ on $\Bbb{F}^{+}(G)$ is
called the admissibility. i.e., the algebraic structure $(\Bbb{F}^{+}(G),$ $%
\cdot )$ is the free semigroupoid of $G.$ For convenience, we denote $(\Bbb{F%
}^{+}(G),$ $\cdot )$ simply by $\Bbb{F}^{+}(G).$

\strut

For the given countable directed graph $G,$ we can define a new countable
directed graph $G^{-1}$ which is the opposite directed graph of $G,$ with

\strut

\begin{center}
$V(G^{-1})=V(G)$ \ \ and $\ \ E(G^{-1})=\{e^{-1}:e\in E(G)\}$,
\end{center}

\strut

where $e^{-1}$ $\in $ $E(G^{-1})$ is the opposite directed edge of $e$ $\in $
$E(G)$, called the shadow of $e$ $\in $ $E(G).$ i.e., if $e$ $=$ $v_{1}$ $e$ 
$v_{2}$ in $E(G)$ with $v_{1},$ $v_{2}$ $\in $ $V(G),$ then $e^{-1}$ $=$ $%
v_{2}$ $e^{-1}$ $v_{1}$ in $E(G^{-1})$ with $v_{2},$ $v_{1}$ $\in $ $%
V(G^{-1})$ $=$ $V(G).$ This new directed graph $G^{-1}$ is said to be the
shadow of $G$. It is trivial that $(G^{-1})^{-1}$ $=$ $G.$ This relation
shows that the admissibility on the shadow $G^{-1}$ is oppositely preserved
by that on $G.$

\strut

A new countable directed graph $G^{\symbol{94}}$ is called the shadowed
graph of $G$ if it is a directed graph with

\strut

\begin{center}
$V(G^{\symbol{94}})=V(G)=V(G^{-1})$
\end{center}

and

\begin{center}
$E(G^{\symbol{94}})=E(G)\cup E(G^{-1}).$
\end{center}

\strut \strut \strut

\begin{definition}
Let $G$ be a countable directed graph and $G^{\symbol{94}},$ the shadowed
graph of $G,$ and let $\Bbb{F}^{+}(G^{\symbol{94}})$ be the free
semigroupoid of $G^{\symbol{94}}.$ Define the reduction (RR) on $\Bbb{F}%
^{+}(G^{\symbol{94}})$ by

\strut

(RR)$\ \ \ \ \ \ \ \ \ \ \ \ \ \ \ \ \ \ ww^{-1}=v$ \ \ and \ \ $%
w^{-1}w=v^{\prime },$

\strut

whenever $w$ $=$ $vwv^{\prime }$ in $FP(G^{\symbol{94}}),$ with $v,$ $%
v^{\prime }$ $\in $ $V(G^{\symbol{94}})$. The set $\Bbb{F}^{+}(G^{\symbol{94}%
})$ with this reduction (RR) is denoted by $\Bbb{F}_{r}^{+}(G^{\symbol{94}%
}). $ And this set $\Bbb{F}_{r}^{+}(G^{\symbol{94}})$ with the inherited
admissibility $(\cdot )$ from $\Bbb{F}^{+}(G^{\symbol{94}})$ is called the
graph groupoid of $G.$ Denote $(\Bbb{F}_{r}^{+}(G^{\symbol{94}}),$ $\cdot )$
of $G$ by $\Bbb{G}.$
\end{definition}

\strut \strut \strut \strut

The graph groupoid $\Bbb{G}$ of $G$ is a (categorial) groupoid with its base 
$V(G^{\symbol{94}}).$ Define the reduced finite path set $FP_{r}(G^{\symbol{%
94}})$ of $\Bbb{G}$ by

$\strut $

\begin{center}
$FP_{r}(G^{\symbol{94}})$ $\overset{def}{=}$ $\Bbb{F}_{r}^{+}(G^{\symbol{94}%
})$ $\setminus $ $\left( V(G^{\symbol{94}})\text{ }\cup \text{ }\{\emptyset
\}\right) .$
\end{center}

\strut

The subset $FP_{r}(G^{\symbol{94}})$ of $\Bbb{G}$ is called the reduced
finite path set of $G.$ Let $w_{1}$ and $w_{2}$ be reduced finite paths in $%
FP_{r}(G^{\symbol{94}})$ $\subset $ $\Bbb{G}.$ We will use the same notation 
$w_{1}$ $w_{2}$ for the product of $w_{1}$ and $w_{2}$ in $\Bbb{G}.$ But we
have to keep in mind that the product $w_{1}$ $w_{2}$ in the graph groupoid $%
\Bbb{G}$ is different from the product $w_{1}$ $w_{2}$ in the free
semigroupoid $\Bbb{F}^{+}(G^{\symbol{94}}).$ Suppose $e_{1}$ and $e_{2}$ are
edges in $E(G^{\symbol{94}})$ and assume that they are admissible, and hence 
$e_{1}$ $e_{2}$ is a finite path in $FP(G^{\symbol{94}}).$ Then the product $%
e_{1}$ $e_{2}$ $e_{2}^{-1}$ of $e_{1}$ $e_{2}$ and $e_{2}^{-1}$ is the
length-3 finite path in $FP(G^{\symbol{94}})$ $\subset $ $\Bbb{F}^{+}(G^{%
\symbol{94}}),$ but the product $e_{1}$ $e_{2}$ $e_{2}^{-1}$ of them is $%
e_{1}$ $(e_{2}$ $e_{2}^{-1})$ $=$ $e_{1}$ in $FP_{r}(G^{\symbol{94}})$ $%
\subset $ $\Bbb{G}.$

\strut \strut \strut

Now, we will define and consider certain operators on graph Hilbert spaces.
As we defined and observed, we can regard all elements in $\Bbb{G}$ as
reduced words in $E(G^{\symbol{94}}),$ under the admissibility with the
reduction (RR).

\strut

\begin{definition}
Let $G$ be a countable directed graph and let $\Bbb{G}$ be the corresponding
graph groupoid. Define the Hilbert space $H_{G}$ of $G$ by

\strut

\begin{center}
$H_{G}$ $\overset{def}{=}$ $\ \underset{w\in FP_{r}(G^{\symbol{94}})}{\oplus 
}\left( \Bbb{C}\xi _{w}\right) $
\end{center}

\strut

with its Hilbert basis $\mathcal{B}_{H_{G}}$ $=$ $\{\xi _{w}$ $:$ $w$ $\in $ 
$FP_{r}(G^{\symbol{94}})\}$ in $H_{G},$ where $FP_{r}(G^{\symbol{94}})$ is
the reduced finite path set.
\end{definition}

\strut \strut \strut

We have the following multiplication rule on $H_{G}$:

\strut

\begin{center}
$\xi _{w_{1}}\xi _{w_{2}}=\left\{ 
\begin{array}{ll}
\xi _{w_{1}w_{2}} & \text{if }w_{1}\text{ }w_{2}\text{ }\neq \text{ }%
\emptyset \\ 
\xi _{\emptyset }=0_{H_{G}} & \text{otherwise,}
\end{array}
\right. $
\end{center}

\strut

for all $\xi _{w_{1}},$ $\xi _{w_{2}}$ $\in $ $\mathcal{B}_{H_{G}}.$ Suppose 
$w_{1}$ $=$ $w$ and $w_{2}$ $=$ $w^{-1}$ in $FP_{r}(G^{\symbol{94}}).$ Then,
by the above multiplication rule, we can have $\xi _{w}$ $\xi _{w^{-1}}$ $=$ 
$\xi _{ww^{-1}}$ and $w$ $w^{-1}$ is a vertex in $V(G^{\symbol{94}}).$ So,
we can determine $\xi _{v}$ $\in $ $H_{G},$ for all $v$ $\in $ $V(G^{\symbol{%
94}}).$ And hence, we can extend the above multiplication rule to the case
where $w_{1},$ $w_{2}$ $\in $ $\Bbb{G}.$ This multiplication rule let us
define multiplication operators on $H_{G}.$

\strut

\begin{definition}
An operator $R_{w}$ on $H_{G}$ is defined by the right multiplication
operator with its symbol $\xi _{w}$ on $H_{G},$ for $w$ $\in $ $\Bbb{G}.$
i.e.,

\strut

\begin{center}
$R_{w}\xi _{w^{\prime }}\overset{def}{=}\left\{ 
\begin{array}{ll}
\xi _{w^{\prime }}\xi _{w}=\xi _{w^{\prime }w} & \text{if }w^{\prime }\text{ 
}w\neq \emptyset \\ 
\xi _{\emptyset }=0_{H_{G}} & \text{otherwise,}
\end{array}
\right. $
\end{center}

\strut

for all $w,$ $w^{\prime }$ $\in $ $\Bbb{G}.$ The adjoint $R_{w}^{*}$ of $%
R_{w}$ is defined by $R_{w}^{*}$ $=$ $R_{w^{-1}},$ for all $w$ $\in $ $\Bbb{G%
}.$
\end{definition}

\strut

By definition, the product of two right multiplication operators $R_{w_{1}}$ 
$R_{w_{2}}$ is the multiplication operator $R_{w_{2}w_{1}}.$ i.e.,

$\strut $

\begin{center}
$R_{w_{1}}$ $R_{w_{2}}$ $=$ $R_{w_{2}w_{1}},$ \ for all $w_{1},$ $w_{2}$ $%
\in $ $\Bbb{G}.$
\end{center}

\strut

So, it is easy to check that the multiplication operators $R_{v},$ for all $%
v $ $\in $ $V(G^{\symbol{94}})$, are projections on $H_{G},$ since

$\strut $

\begin{center}
$R_{v}^{2}$ $=$ $R_{v}$ $R_{v}$ $=$ $R_{v^{2}}$ $=$ $R_{v}$ and $R_{v}^{*}$ $%
=$ $R_{v^{-1}}$ $=$ $R_{v}.$
\end{center}

\strut

Thus we can get that the right multiplication operators $R_{w},$ for all $w$ 
$\in $ $FP_{r}(G^{\symbol{94}}),$ are partial isometries. Indeed,

\strut

\begin{center}
$R_{w}R_{w}^{*}R_{w}=R_{ww^{-1}w}=R_{ww^{-1}}R_{w}=R_{w}$ \ 
\end{center}

and \ 

\begin{center}
$%
R_{w}^{*}R_{w}R_{w}^{*}=R_{w^{-1}ww^{-1}}=R_{w^{-1}w}R_{w^{-1}}=R_{w^{-1}}=R_{w}^{*}, 
$
\end{center}

\strut

for all $w$ $\in $ $FP_{r}(G^{\symbol{94}}).$\strut

\strut

\begin{remark}
In [10] and [11], we defined the same graph Hilbert space $H_{G},$ for any
countable directed graph $G,$ and we defined the (left) multiplication
operators $L_{w}$ $\in $ $B(H_{G})$ by

\strut

\begin{center}
$L_{w}$ $\xi _{w^{\prime }}$ $\overset{def}{=}$ $\left\{ 
\begin{array}{ll}
\xi _{ww^{\prime }} & \text{if }w\text{ }w^{\prime }\neq \emptyset \\ 
0 & \text{otherwise,}
\end{array}
\right. $
\end{center}

\strut

for all $w,$ $w^{\prime }$ $\in $ $\Bbb{G}.$ Then they have the similar
properties like the right multiplication operators $R_{w}$ on $H_{G}.$ i.e.,
if $w$ is a vertex, then $L_{w}$ is a projection: and if $w$ is a reduced
finite path, then $L_{w}$ is a partial isometries on $H_{G}.$
\end{remark}

\strut

Now, we define certain groupoid actions $\beta $ of the graph groupoid $\Bbb{%
G}$ on $H_{G},$ called right graph-representations.

\strut \strut \strut \strut

\begin{definition}
Let $\Bbb{G}$ be a graph groupoid of a countable directed graph $G$ and let $%
M$ be a von Neumann algebra acting on a Hilbert space $K$ (i.e., $M$ $%
\subseteq $ $B(K)$)$.$ Define a right graph-representation (in short, a
right $G$-representation) $\beta $ $:$ $\Bbb{G}$ $\rightarrow $ $B(K$ $%
\otimes $ $H_{G})$, by a nonunital intertwined partial representation
satisfying that

\strut

\begin{center}
$\beta _{w}(m)$ $R_{w}R_{w}^{*}=R_{w}^{*}mR_{w}=R_{w^{-1}}mR_{w},$
\end{center}

\strut

for all $m$ $\in $ $M$ and $w$ $\in $ $\Bbb{G},$ and

\strut

\begin{center}
$\beta _{v}(m)$ $=$ $m,$ for all $m$ $\in $ $M$ and $v$ $\in $ $V(G^{\symbol{%
94}}).$
\end{center}

\strut

In this setting, we can regard $R_{w}$'s as $1_{K}$ $\otimes $ $R_{w}$'s on $%
K$ $\otimes $ $H_{G}$, for all $w$ $\in $ $\Bbb{G}.$
\end{definition}

\strut \strut \strut \strut \strut \strut \strut \strut \strut \strut \strut
\strut \strut

For a fixed von Neumann algebra $M$ and a graph groupoid $\Bbb{G},$ we can
construct a groupoid crossed product algebra $M$ $\times _{\beta }$ $\Bbb{G}%
, $ for a fixed right graph-representation $\beta .$ This new von Neumann
algebra is called a right graph von Neumann algebra induced by $G$ over $M.$

\strut

\begin{definition}
Let $M$ be a von Neumann algebra and $\Bbb{G},$ the graph groupoid of a
countable directed graph $G$ and let $\beta $ be a right $G$-representation.
Define the crossed product $\Bbb{M}_{G}$ $=$ $M$ $\times _{\beta }$ $\Bbb{G}$
of $M$ and $\Bbb{G}$ via $\beta ,$ by the von Neumann algebra generated by $%
M $ and $\{R_{w}$ $:$ $w$ $\in $ $\Bbb{G}\}$ in $B(K$ $\otimes $ $H_{G}),$
satisfying the above $\beta $-rules on $M$. This von Neumann algebra $\Bbb{M}%
_{G}$ is called a right graph von Neumann algebra induced by $G$ over $M.$
\end{definition}

\strut \strut

\begin{remark}
In [10] and [11], we defined a graph von Neumann algebra $\Bbb{M}%
_{G}^{(left)}$ by the groupoid crossed product $M$ $\times _{\alpha }$ $\Bbb{%
G},$ where $\alpha $ $:$ $\Bbb{G}$ $\rightarrow $ $B(K$ $\otimes $ $H_{G})$
is a $G$-representation satisfying the $\alpha $-rules on $M$:

\strut

\begin{center}
$\alpha _{w}(m)$ $L_{w}$ $L_{w}^{*}$ $=$ $L_{w}^{*}$ $m$ $L_{w},$ for all $w$
$\in $ $FP_{r}(G^{\symbol{94}}),$
\end{center}

and

\begin{center}
$\alpha _{v}(m)$ $=$ $m,$ for all $v$ $\in $ $V(G^{\symbol{94}}),$
\end{center}

\strut

for all $m$ $\in $ $M,$ where $L_{w}$'s are the (left) multiplication
operators, for $w$ $\in $ $\Bbb{G}$ (See the previous remark). Notice that,
in fact, the von Neumann algebra $\Bbb{M}_{G}^{(left)}$ is the opposite
algebra $\Bbb{M}_{G}^{op}$ of a right graph von Neumann algebra $\Bbb{M}_{G}$
$=$ $M$ $\times _{\beta }$ $\Bbb{G},$ whenever $\beta $ $=$ $\alpha $ on $M$
(See [13]).
\end{remark}

\strut \strut

Every operator $x$ in a right graph von Neumann algebra $\Bbb{M}_{G}$ $=$ $M$
$\times _{\beta }$ $\Bbb{G}$ has its expression,

\strut

\begin{center}
$x=\underset{w\in \Bbb{G}}{\sum_{n}^{(N)}}m_{w}$ $R_{w},$ \ for $m_{w}$ $\in 
$ $M.$
\end{center}

\strut \strut

\begin{lemma}
Let $m_{1}R_{w_{1}},$ ..., $m_{n}R_{w_{n}}$ be operators in a right graph
von Neumann algebra $\Bbb{M}_{G},$ for $n$ $\in $ $\Bbb{N}.$ Then

\strut \strut

$\Pi _{k=1}^{n}\left( m_{k}R_{w_{k}}\right) $

\strut

\begin{center}
$=\left\{ 
\begin{array}{ll}
\left(
m_{1}m_{2}^{w_{1}^{-1}}m_{3}^{(w_{1}w_{2})^{-1}}...m_{n}^{(w_{1}w_{2}...w_{n-1})^{-1}}\right) R_{w_{n}...w_{1}}
& \text{if }w_{n}...w_{1}\neq \emptyset \\ 
&  \\ 
0_{\Bbb{M}_{G}} & \text{otherwise,}
\end{array}
\right. $
\end{center}

\strut

where $m^{w}$ $\overset{def}{=}$ $\beta _{w}(m),$ for all $m$ $\in $ $M$ and 
$w$ $\in $ $\Bbb{G}.$ $\square $
\end{lemma}

\strut \strut

\strut \strut

\subsection{$M$-Valued Right Graph $W^{*}$-Probability Spaces}

\strut

\strut

Let $G$ be a countable directed graph and let $M$ be an arbitrary von
Neumann algebra acting on a Hilbert space $K$. In this section, we will
define a $M$-diagonal right graph $W^{*}$-probability space $(\Bbb{M}_{G},$ $%
E),$ over its $M$-diagonal subalgebra $\Bbb{D}_{G}$, where $\Bbb{M}_{G}$ $=$ 
$M$ $\times _{\beta }$ $\Bbb{G}$ is a right graph von Neumann algebra.

\strut

Let $v$ $\in $ $V(G^{\symbol{94}}).$ Then we can define a conditional
expectation $E_{v}$ $:$ $\Bbb{M}_{G}$ $\rightarrow $ $M$ $\cdot $ $R_{v},$
where $M$ $\cdot $ $R_{v}$ $=$ $\{m$ $R_{v}$ $:$ $m$ $\in $ $M\},$ by

\strut

\begin{center}
$E_{v}\left( \underset{w\in \Bbb{G}}{\sum_{n}^{(N)}}m_{w}R_{w}\right)
=m_{v}R_{v},$
\end{center}

\strut

for all $\underset{w\in \Bbb{G}}{\sum_{n}^{(N)}}$ $m_{w}R_{w}$ $\in $ $\Bbb{M%
}_{G}.$ Notice that each $M$ $\cdot $ $R_{v}$ is a $W^{*}$-subalgebra of $%
\Bbb{M}_{G},$ which is $*$-isomorphic to $M,$ since $R_{v}$ is a projection,
for all $v$ $\in $ $V(G^{\symbol{94}}).$ A pair $(\Bbb{M}_{G},$ $E_{v})$ is
a $M$-valued $W^{*}$-probability space, for $v$ $\in $ $V(G^{\symbol{94}}).$
We call it a vertex-depending (or the $v$-depending) right graph $W^{*}$%
-probability space over $M.$ The conditional expectation $E_{v}$ is said to
be a vertex-depending (or the $v$-depending) conditional expectation, for $v$
$\in $ $V(G^{\symbol{94}}).$ By the very definition, the right graph von
Neumann algebra $\Bbb{M}_{G}$ has $\left| V(G^{\symbol{94}})\right| $-many
vertex-depending $W^{*}$-probability spaces over $M.$

\strut

\begin{definition}
By $\Bbb{D}_{G},$ denote a $W^{*}$-subalgebra $\underset{v\in V(G^{\symbol{94%
}})}{\oplus }$ $(M$ $\cdot $ $R_{v})$ of $\Bbb{M}_{G}$ $=$ $M$ $\times
_{\beta }$ $\Bbb{G}.$ This subalgebra $\Bbb{D}_{G}$ is called the $M$%
-diagonal subalgebra of $\Bbb{M}_{G}.$ Define a conditional expectation $E$ $%
:$ $\Bbb{M}_{G}$ $\rightarrow $ $\Bbb{D}_{G}$, by $\underset{v\in V(G^{%
\symbol{94}})}{\oplus }$ $E_{v},$ where $E_{v}$'s are the $v$-depending
conditional expectations, for all $v$ $\in $ $V(G^{\symbol{94}}).$ i.e.,

\strut 

\begin{center}
$E\left( \underset{w\in \Bbb{G}}{\sum_{n}^{(N)}}m_{w}R_{w}\right) =\underset{%
v\in V(G^{\symbol{94}})}{\sum_{n}^{(N)}}m_{v}R_{v},$
\end{center}

\strut 

for all $\underset{w\in \Bbb{G}}{\sum_{n}^{(N)}}m_{w}R_{w}$ $\in $ $\Bbb{M}%
_{G}.$ The pair $(\Bbb{M}_{G},$ $E)$ is called the $M$-diagonal right graph $%
W^{*}$-probability space over $\Bbb{D}_{G}.$
\end{definition}

\strut \strut \strut

\begin{definition}
Let $G$ be a countable directed graph and let $\Bbb{G}$ be the graph
groupoid of $G.$ Define a map $\delta $ $:$ $\Bbb{G}$ $\rightarrow $ $\Bbb{G}
$ by mapping $w$ $\in $ $\Bbb{G}$ to the graphical image $\delta (w)$ of $w,$
for all $w$ $\in $ $\Bbb{G}.$ The map $\delta $ on $\Bbb{G}$ is called the
diagram map and the graphical image $\delta (w)$ of $w$ is called the
diagram of $w$, for all $w$ $\in $ $\Bbb{G}.$ We say that the elements $w_{1}
$ and $w_{2}$ in $\Bbb{G}$ are diagram-distinct if (i) $w_{1}$ $\neq $ $%
w_{2}^{-1}$ and (ii) $\delta (w_{1})$ $\neq $ $\delta (w_{2}).$ Suppose that 
$X_{1}$ and $X_{2}$ are subsets in $\Bbb{G}.$ They are said to be
diagram-distinct if, for any pair $(w_{1},$ $w_{2})$ in $X_{1}$ $\times $ $%
X_{2},$ $w_{1}$ and $w_{2}$ are diagram-distinct.
\end{definition}

\strut

Let $l$ be a loop edge in $E(G^{\symbol{94}}),$ and let $w_{1}$ $=$ $%
l^{k_{1}}$ and $w_{2}$ $=$ $l^{k_{2}}$ in $\Bbb{G},$ for $k_{1}$ $\neq $ $%
k_{2}$ $\in $ $\Bbb{N}.$ Clearly, $w_{1}$ $\neq $ $w_{2}^{-1}.$ But the
elements $w_{1}$ and $w_{2}$ are not diagram-distinct, because their
diagrams $\delta (w_{1})$ and $\delta (w_{2})$ are identical to the diagram $%
\delta (l)$ $=$ $l$ of $l.$ Suppose either $w_{1}$ or $w_{2}$ is a non-loop
finite path such that (i) $w_{1}$ $\neq $ $w_{2}^{-1},$ and (ii) $\delta
(w_{k})$'s are not loop finite paths, for $k$ $=$ $1,$ $2.$ Then they are
diagram-distinct whenever $w_{1}$ and $w_{2}$ are distinct in $\Bbb{G}.$

\strut \strut \strut

\begin{theorem}
The subsets $M$ $\cdot $ $R_{w_{1}}$ and $M$ $\cdot $ $R_{w_{2}}$ of a right
graph von Neumann algebra $\Bbb{M}_{G}$ $=$ $M$ $\times _{\beta }$ $\Bbb{G}$
are free over $D_{G}$ in $(\Bbb{M}_{G},$ $E)$ if and only if $w_{1}$ and $%
w_{2}$ are diagram-distinct in $\Bbb{G}$, where $M$ $\cdot $ $R_{w}$ $%
\overset{def}{=}$ $\{m$ $L_{w}$ $:$ $m$ $\in $ $M\},$ for all $w$ $\in $ $%
\Bbb{G}.$ $\square $
\end{theorem}

\strut \strut \strut \strut \strut

We will skip the proof. The readers can find the proof showing that $M$ $%
\cdot $ $L_{w_{1}}$ and $M$ $\cdot $ $L_{w_{2}}$ are free over $\Bbb{D}_{G}$
in a graph von Neumann algebra $M$ $\times _{\alpha }$ $\Bbb{G}$, in [10]
and [11]. By the slight modification, we can prove the above theorem. As
corollary of the previous theorem, we can see that if $e_{1}$ and $e_{2}$
are edges in $E(G)$ ($\subset $ $E(G^{\symbol{94}})$), then $M$ $\cdot $ $%
R_{e_{1}}$ and $M$ $\cdot $ $R_{e_{2}}$ are free over $\Bbb{D}_{G}$ in $(%
\Bbb{M}_{G},$ $E).$ Indeed, if two edges $e_{1}$ and $e_{2}$ are distinct in 
$E(G),$ then $e_{1}$ and $e_{2}^{\pm 2}$ (resp., $e_{2}$ and $e_{1}^{\pm 1}$%
) are diagram-distinct in $E(G^{\symbol{94}})$ $\subset $ $\Bbb{G}$.

\strut

Let's consider the subset $\{e,$ $e^{-1}\}$ of $\Bbb{G},$ for $e$ $\in $ $%
E(G)$ $\subset $ $E(G^{\symbol{94}}).$ Then we can define a subset $\Bbb{G}%
_{e}$ of $\Bbb{G}$, by the collection of all reduced words in $\{e,$ $%
e^{-1}\}.$ It is easy to check that: if $e$ is a loop edge in $E(G)$, then $%
\Bbb{G}_{e}$ is a group, which is group-isomorphic to the infinite abelian
cyclic group $\Bbb{Z},$ and if $e$ is a non-loop edge in $E(G)$, then $\Bbb{G%
}_{e}$ is a groupoid $\{\emptyset ,$ $v_{1},$ $v_{2},$ $e,$ $e^{-1}\},$
where $e$ $=$ $v_{1}$ $e$ $v_{2},$ with $v_{1},$ $v_{2}$ $\in $ $V(G).$
Notice that the graph groupoid $\Bbb{G}$ is the reduced free product $%
\underset{e\in E(G)}{*}$ $\Bbb{G}_{e}$ of $\Bbb{G}_{e}$'s (in the sense of
[10]). Moreover, each $\Bbb{G}_{e}$ can be regarded as a new graph groupoid
induced by a one-edge graph $G_{e},$ with $V(G_{e})$ $=$ $\{v_{1},$ $v_{2}\}$
and $E(G_{e})$ $=$ $\{e$ $=$ $v_{1}$ $e$ $v_{2}\}$, where $v_{1}$ and $v_{2}$
are not necessarily distinct in $V(G).$ Thus we can construct the right
graph von Neumann algebras,

\strut

\begin{center}
$\Bbb{M}_{e}$ $=$ $M$ $\times _{\beta _{e}}$ $\Bbb{G}_{e}$ \strut with $%
\beta _{e}$ $\overset{def}{=}$ $\beta $ $\mid _{\Bbb{G}_{e}}$ on $M,$
\end{center}

\strut

for all $e$ $\in $ $E(G).$ Then they are embedded $W^{*}$-subalgebras of $%
\Bbb{M}_{G}$ $=$ $M$ $\times _{\beta }$ $\Bbb{G}.$ Moreover, by the previous
theorem, we can have that:

\strut

\begin{corollary}
The $W^{*}$-subalgebras $\Bbb{M}_{e}$'s of $\Bbb{M}_{G}$ are free over $\Bbb{%
D}_{G}$, from each other, for all $e$ $\in $ $E(G),$ in $(\Bbb{M}_{G},$ $E).$
$\square $
\end{corollary}

\strut \strut \strut \strut \strut \strut

\begin{definition}
The $W^{*}$-subalgebra $\Bbb{M}_{e}$ of $\Bbb{M}_{G}$ are called $\Bbb{D}_{G}
$-free blocks of $\Bbb{M}_{G},$ for all $e$ $\in $ $E(G).$
\end{definition}

\strut \strut \strut

Notice that $\Bbb{M}_{e}$ is identically same as $\Bbb{M}_{e^{-1}}$, for all 
$e$ $\in $ $E(G).$\ So, we can get the following theorem characterizing the $%
\Bbb{D}_{G}$-free structure of $(\Bbb{M}_{G},$ $E).$\strut 

\strut \strut 

The following three theorems are about the crossed product algebra $\Bbb{M}%
_{G}$ from Theorem 2.4. The first two theorems both involve the generalized
diagonal $\Bbb{D}_{G},$ and they offer decompositions of $\Bbb{M}_{G}$ into
certain ``blocks''. Each of the two decompositions generalizes what is known
for finite matrices. The first one (Theorem 2.6) is an algebraic
decomposition of $\Bbb{M}_{G},$ while the second (Theorem 2.7) offers a
canonical decomposition of $\Bbb{M}_{G}$ as a Banach space. Theorem 2.8 in
turn proves that the blocks in the decompositions fall in two specific
classes.

\strut

\begin{theorem}
Let $(\Bbb{M}_{G},$ $E)$ be the $M$-diagonal right graph $W^{*}$-probability
space over its $M$-diagonal subalgebra $\Bbb{D}_{G}.$ Then

\strut 

\begin{center}
$\Bbb{M}_{G}$ $\overset{*\text{-isomorphic}}{=}$ $\underset{e\in E(G)}{%
\,\,*_{\Bbb{D}_{G}}}$ $\Bbb{M}_{e},$
\end{center}

\strut 

where $\Bbb{M}_{e}$'s are the $\Bbb{D}_{G}$-free blocks of $\Bbb{M}_{G},$
for all $e$ $\in $ $E(G).$ $\square $
\end{theorem}

\strut \strut

\strut Recall that if $\{A_{i}\}_{i\in I}$ are von Neumann algebras over
their common $W^{*}$-subalgebra $B,$ then the $B$-free product $\underset{%
i\in I}{\,\,*_{B}}$ $A_{i}$ can be expressed by

\strut

\begin{center}
$B\oplus \left( \oplus _{n=1}^{\infty }\left( \underset{i_{1}\neq
i_{2},i_{2}\neq i_{3},...,i_{n-1}\neq i_{n}\in I}{\oplus }%
(A_{i_{1}}^{o}\otimes _{B}...\otimes _{B}A_{i_{n}}^{o})\right) \right) ,$
\end{center}

\strut

as a Banach space, where $A_{i_{k}}^{o}$ $\overset{def}{=}$ $A_{i_{k}}$ $%
\ominus $ $B,$ for all $k$ $=$ $1,$ ..., $n.$ Since our right graph von
Neumann algebra $\Bbb{M}_{G}$ is $*$-isomorphic to the $\Bbb{D}_{G}$-free
product $\underset{e\in E(G)}{\,\,*_{\Bbb{D}_{G}}}$ $\Bbb{M}_{e},$ where $%
\Bbb{M}_{e}$'s are the $\Bbb{D}_{G}$-free blocks of $\Bbb{M}_{G},$ we can
get the Banach space expression. However, the $\Bbb{D}_{G}$-tensor products
in our $\Bbb{D}_{G}$-free product is determined by the admissibility on $%
\Bbb{G}.$

\strut

\begin{definition}
Let $G$ be a countable directed graph. Define the subset $E(G^{\symbol{94}%
})_{r}^{*}$ of $FP_{r}(G^{\symbol{94}})$ by

\strut \strut 

\begin{center}
$E(G)_{r}^{*}$ $\overset{def}{=}$ $E(G)$ $\cup $ $\left( \cup _{k=2}^{\infty
}\left\{ e_{1}e_{2}...e_{k-1}e_{k}\left| 
\begin{array}{c}
e_{1}\text{ ... }e_{k}\in FP_{r}(G^{\symbol{94}}), \\ 
e_{1}\neq e_{2}^{\pm 1},\text{ }e_{2}\neq e_{3}^{\pm 1} \\ 
\text{ \ \ }...,\text{ }e_{k-1}\neq e_{k}^{\pm 1}
\end{array}
\right. \right\} \right) .$
\end{center}
\end{definition}

\strut \strut

\strut The following theorem modifies a result from two papers [10] and
[11]. In [10], a generalized diagonal $\Bbb{D}_{G}$ is introduced, and
associated $\Bbb{D}_{G}$-free blocks in a (left) graph von Neumann algebra
are characterized. 

\strut \strut

\begin{theorem}
Let $G$ be a countable directed graph and let $\Bbb{M}_{G}$ $=$ $M$ $\times
_{\beta }$ $\Bbb{G}$ be the corresponding right graph von Neumann algebra,
for a fixed von Neumann algebra $M.$ Then, as a Banach space,

\strut 

\begin{center}
$\Bbb{M}_{G}$ $=$ $\Bbb{D}_{G}\oplus \left( \underset{w\in E(G)_{r}^{*}}{%
\oplus }\,\Bbb{M}_{w}^{o}\right) .$
\end{center}

\strut 

$\square $
\end{theorem}

\strut \strut \strut

\textbf{Notation} \ To emphasize that $\underset{e\in E(G)}{\,*_{\Bbb{D}_{G}}%
}$ $\Bbb{M}_{e}$ has the above Banach-space expression, we will denote it by 
$\underset{e\in E(G)}{\,*_{\Bbb{D}_{G}}^{r}}$ $\Bbb{M}_{e}$. And we will
call it the $\Bbb{D}_{G}$-valued \emph{reduced} free product of the $\Bbb{D}%
_{G}$-free blocks $\Bbb{M}_{e}$'s.

\strut

In [11], the $\Bbb{D}_{G}$-free blocks in a (left) graph von Neumann algebra
are characterized. Similarly, we have the following results.

\strut

\begin{theorem}
Let $\Bbb{M}_{G}$ $=$ $M$ $\times _{\beta }$ $\Bbb{G}$ be a right graph von
Neumann algebra. Then its $\Bbb{D}_{G}$-free block $\Bbb{M}_{e}$ is $*$%
-isomorphic to either $M$ $\times _{\lambda _{e}}$ $\Bbb{Z}$ (whenever $e$
is a loop-edge) or $M$ $\otimes _{\beta _{e}}$ $M_{2}(\Bbb{C})$ (whenever $e$
is a non-loop edge), where $\lambda _{e}$ is the group action of $\Bbb{Z},$
such that $\lambda _{n}(m)$ $=$ $\beta _{e^{n}}(m),$ for all $m$ $\in $ $M$
and $n$ $\in \Bbb{Z}.$ $\square $
\end{theorem}

\strut \strut

\textbf{Assumption} In the rest of this paper, we will restrict our
interests to the case where $M$ $=$ $\Bbb{C}.$

\strut

Suppose $M$ $=$ $\Bbb{C}.$ Then a right graph von Neumann algebra $\Bbb{M}%
_{G}$ $=$ $\Bbb{C}$ $\times _{\beta }$ $\Bbb{G}$ is $*$-isomorphic to $%
\overline{\Bbb{C}[\Bbb{G}]}^{w},$ for any graph-representation $\beta ,$ in $%
B(H_{G}).$ Also, in such case, the $\Bbb{D}_{G}$-free blocks are $*$%
-isomorphic to either $L^{\infty }(\Bbb{T})$ or $M_{2}(\Bbb{C}).$

\strut \strut

\begin{example}
Let $G_{N}$ be a one-vertex-$N$-loop-edge graph with its vertex set $V(G_{N})
$ $=$ $\{v\}$ and $E(G_{N})$ $=$ $\{e_{j}$ $=$ $v$ $e_{j}$ $v$ $:$ $j$ $=$ $%
1,$ ..., $N\}.$ Then the right graph von Neumann algebra $\Bbb{C}_{G_{N}}$ $=
$ $\Bbb{C}$ $\times _{\beta }$ $\Bbb{G}_{N},$ where $\Bbb{G}_{N}$ is the
graph groupoid of $G_{N},$ is $*$-isomorphic to the free group factor $%
L(F_{N})$ (Also See [10]). Notice that the graph groupoid $\Bbb{G}_{N}$ of $%
G_{N}$ is a group and moreover it is group-isomorphic to the free group $%
F_{N}$ $=$ $<$ $g_{1},$ ..., $g_{N}$ $>$ with $N$-generators $\{g_{1},$ ..., 
$g_{N}\}$. Indeed, the edges $e_{1},$ ..., $e_{N}$ are generators of $\Bbb{G}%
_{N}$ as a group. So, there is a natural generator-preserving
group-isomorphism between $\Bbb{G}_{N}$ and $F_{N}.$ Also, notice that the
vertex $v$ in $\Bbb{G}_{N}$ is the group identity. And all partial
isometries $R_{e_{j}^{n}}$ correspond to the unitary representation $%
u_{g_{j}^{n}}$ of $g_{j}^{n},$ for all $n$ $\in $ $\Bbb{Z}$ $\setminus $ $%
\{0\}$ and for $j$ $=$ $1,$ ..., $N.$ Therefore, the right graph von Neumann
algebra $\Bbb{C}_{G_{N}}$ $=$ $\Bbb{C}$ $\times _{\beta }$ $\Bbb{G}_{N}$ is $%
*$-isomorphic to the classical crossed product algebra $\Bbb{C}$ $\times
_{\lambda }$ $F_{N},$ whenever $\lambda $ $=$ $\beta .$ Notice that $\Bbb{C}$
$\times _{\beta }$ $\Bbb{G}_{N}$ is $*$-isomorphic to $\overline{\Bbb{C}[%
\Bbb{G}_{N}]}^{w},$ for any $\beta ,$ by the linearity of $\beta $ on $\Bbb{C%
}.$ Also, it is well-known that a crossed product algebra $\Bbb{C}$ $\times
_{\gamma }$ $\Gamma $ is $*$-isomorphic to $\overline{\Bbb{C}[\Gamma ]}^{w},$
for any group action $\gamma ,$ for all groups $\Gamma .$ Therefore, the
right graph von Neumann algebra $\Bbb{C}_{G_{N}}$ is $*$-isomorphic to the
group von Neumann algebra $\overline{\Bbb{C}[F_{N}]}^{w}.$ More generally,
the right graph von Neumann algebra $\Bbb{M}_{G_{N}}$ $=$ $M$ $\times
_{\beta }$ $\Bbb{G}_{N}$ is $*$-isomorphic to the classical crossed product
von Neumann algebra $M$ $\times _{\lambda }$ $F_{N},$ whenever $\lambda $ $=$
$\beta .$
\end{example}

\strut \strut \strut \strut

Notice that if $G$ is a countable directed graph with its graph groupoid $%
\Bbb{G}$ and if $\Bbb{M}_{G}$ $=$ $\Bbb{C}$ $\times _{\beta }$ $\Bbb{G}$ is
a right graph von Neumann algebra induced by $G$ over $\Bbb{C},$ then we can
see that $\Bbb{M}_{G}$ is $*$-isomorphic to the von Neumann algebra $%
\overline{\Bbb{C}[\Bbb{G}]}^{w},$ generated by $\Bbb{G},$ by the linearity
of $\beta $ on $\Bbb{C}$. Thus, for any graph-representation $\beta ,$ the
corresponding right graph von Neumann algebra $\Bbb{C}$ $\times _{\beta }$ $%
\Bbb{G}$ is $*$-isomorphic to $\overline{\Bbb{C}[\Bbb{G}]}^{w}.$

\strut \strut

\begin{example}
Let $C_{N}$ be the one-flow circulant graph with $V(C_{N})$ $=$ $\{v_{1},$
..., $v_{N}\}$ and $E(C_{N})$ $=$ $\{e_{j}$ $=$ $v_{j}$ $e_{j}$ $v_{j+1}$ $:$
$j$ $=$ $1,$ ..., $N,$ $v_{N+1}$ $\overset{def}{=}$ $v_{1}\}.$ Then the
right graph von Neumann algebra $\Bbb{C}_{C_{N}}$ $=$ $\Bbb{C}$ $\times
_{\beta }$ $\Bbb{G}_{C_{N}}$ contains $W^{*}$-subalgebras which are $*$%
-isomorphic to the group von Neumann algebra $\overline{\Bbb{C}[\Bbb{Z}]}%
^{w}.$ Since a finite path $w$ $=$ $e_{1}$ ... $e_{N}$ induces the subset $%
\Bbb{S}_{w}$ of $\Bbb{G}_{C_{N}}$ consisting of all reduced words in $S_{w}$ 
$=$ $\{w,$ $w^{-1}\}$ in $\Bbb{G}_{C_{N}},$ we have the $W^{*}$-subalgebra $%
N $ $=$ $\Bbb{C}$ $\times _{\beta }$ $\Bbb{S}_{w}$ in $\Bbb{C}_{C_{N}}.$ We
can regard this $W^{*}$-subalgebra $N$ as a right graph von Neumann algebra
induced by the graph $G_{w}$ with its vertex set $V(G_{w})$ $=$ $\{v_{1}\}$
and its edge set $E(G_{w})$ $=$ $\{w$ $=$ $v_{1}$ $w$ $v_{1}\}.$ Then it is
an one-vertex-one-loop-edge graph. So, by the previous example, the von
Neumann algebra $N$ is $*$-isomorphic to $L(F_{1}),$ where $F_{1}$ is the
free group with $1$-generator which is group-isomorphic to $\Bbb{Z}.$ So, $N$
is $*$-isomorphic to $L(\Bbb{Z})$ which is $*$-isomorphic to $L^{\infty }(%
\Bbb{T}),$ where $\Bbb{T}$ is the unit circle in $\Bbb{C}.$ Consider $w^{k},$
for all $k$ $\in $ $\Bbb{Z}$ $\setminus $ $\{0\}.$ Then $W^{*}$-subalgebras $%
N$ $=$ $\Bbb{C}$ $\times _{\beta }$ $\Bbb{G}_{w^{k}}$ of $\Bbb{C}_{C_{N}}$
are all $*$-isomorphic to $L(\Bbb{Z}).$

\strut

We can easily check that $E(C_{N})_{r}^{*}$ $=$ $FP(C_{N})$ $\cup $ $%
FP(C_{N}^{-1}),$ for all $N$ $\in $ $\Bbb{N}.$ Therefore,

\strut

$\ \ \ \ \ \ \ \ \Bbb{C}_{C_{N}}$ $=$ $\Bbb{C}$ $\times _{\beta }$ $\Bbb{G}%
_{C_{N}}$ $=$ $\underset{k=1}{\overset{N}{\,\,\,*_{\Bbb{D}_{C_{N}}}}}$ $\Bbb{%
C}_{e_{k}}$

\strut

where $\Bbb{C}_{e_{k}}$ $\overset{def}{=}$ $vN\left( \Bbb{C}\times _{\beta }%
\text{ }\Bbb{G}_{e},\text{ }\Bbb{D}_{C_{N}}\right) ,$ for all $k$ $=$ $1,$
..., $N$

\strut

\ \ $\ \ \ \ \ \ \ \ \ \ \ \ =\Bbb{D}_{C_{N}}$ $\oplus $ $\left( \underset{%
w^{*}\in E(C_{N})_{r}^{*}}{\oplus }\text{ }\Bbb{C}_{w^{*}}^{o}\right) $

\strut

$\ \ \ \ \ \ \ \ \ \ \ \ \ \ =$ $\Bbb{D}_{C_{N}}\oplus \left( \underset{w\in
FP(C_{N})\cup FP(C_{N}^{-1})}{\oplus }\,\Bbb{C}_{w}^{o}\right) ,$

\strut

where $\Bbb{C}_{w}^{o}$ $\overset{def}{=}$ $\Bbb{C}_{e^{(1)}}^{o}$ $\otimes
_{\Bbb{D}_{C_{N}}}$ ... $\otimes _{\Bbb{D}_{C_{N}}}$ $\Bbb{C}_{e^{(n)}}^{o},$
whenever $w$ $=$ $e^{(1)}$ ... $e^{(n)}$ in $E(C_{N})_{r}^{*},$ for $n$ $\in 
$ $\Bbb{N}.$
\end{example}

\strut \strut

\strut

\strut

\section{Labeled Graph Groupoids and Graph Automata}

\strut

\strut \strut \strut

We will give suitable weights to the elements of graph groupoids induced by
the given locally finite connected directed graphs. Then these weights give
more accurate information of the admissibility on graph groupoids. Such
admissibility conditions are explained by the labeling map and the shifting
(map) of the automata constructed by graph groupoids. Also, this process
would show how the admissibility on a graph groupoid works on the
corresponding graph Hilbert space.

\strut

Recall that a countable directed graph $G$ is \emph{locally finite}, if
every vertex of $G$ has only finitely many incident edges (equivalently, the
degree of each vertex is finite). Also, we say that $G$ is \emph{connected},
if, for any pair $(v_{1},$ $v_{2})$ of distinct vertices ($v_{1}$ $\neq $ $%
v_{2}$), there always exists at least one reduced finite path $w$ $\in $ $%
FP_{r}(G^{\symbol{94}})$ such that $w$ $=$ $v_{1}$ $w$ $v_{2}$ and $w^{-1}$ $%
=$ $v_{2}$ $w^{-1}$ $v_{1}.$

\strut

\textbf{Assumption} From now on, all given directed graphs are locally
finite and connected. $\square $

\strut

Let $G$ be a locally finite connected directed graph and let $v_{0}$ be an
any fixed vertex of $G.$ Then we can define the out-degree $\deg
_{out}(v_{0})$, the in-degree $\deg _{in}(v_{0}),$ and the degree $\deg
(v_{0})$ of $v_{0}$ as follows:

\strut

\begin{center}
$\deg _{out}(v_{0})$ $\overset{def}{=}$ $\left| \{e\in E(G):e=v_{0}\text{ }%
e\}\right| ,$
\end{center}

\strut

\begin{center}
$\deg _{in}(v_{0})$ $\overset{def}{=}$ $\left| \{e\in E(G):e=\text{ }e\text{ 
}v_{0}\}\right| ,$
\end{center}

and

\begin{center}
$\deg (v_{0})$ $\overset{def}{=}$ $\deg _{out}(v_{0})$ $+$ $\deg
_{in}(v_{0}).$
\end{center}

By the locally-finiteness of $G,$ the above three cardinalities are less
than $\infty .$ We are interested in out-degrees.

\strut

\strut

\strut

\subsection{Labeled Graph Groupoids}

\strut

\strut

As before, let $G$ be a locally finite connected directed graph. By the
locally finiteness of $G,$ we can have that

\strut

\begin{center}
$\max \{\deg _{out}(v)$ $:$ $v$ $\in $ $V(G)\}$ $<$ $\infty .$
\end{center}

\strut \strut

Denote the maximal value of the out-degrees of all vertices by $N$. Then we
can define a set

\strut

\begin{center}
$Y$ $=$ $\{1,$ ..., $N\}$ and $Y_{0}$ $=$ $\{0\}$ $\cup $ $Y.$
\end{center}

\strut

Also, fix sets $X$ and $X_{0}$ of certain vectors in $\Bbb{C}^{\oplus \,N},$
indexed by $Y$ and $Y_{0},$ respectively:

\strut

\begin{center}
$X$ $=$ $\{x_{1},$ ..., $x_{N}\}$ and $X_{0}$ $=$ $\{x_{0},$ $x_{1},$ ..., $%
x_{N}\},$
\end{center}

\strut where

\begin{center}
$x_{0}$ $=$ $\left( \underset{N\text{-times}}{\underbrace{0,........,0}}%
\right) $ and $x_{j}$ $=$ $\left( 0,\text{ ..., }0,\text{ }\underset{j\text{%
-th}}{1},\text{ }0,\text{ ..., }0\right) ,$
\end{center}

\strut

for all $j$ $=$ $1,$ ..., $N.$ Remark that the number $N$ is the maximal
out-degree of vertices of $G.$

\strut

\begin{definition}
Let $X$ and $X_{0}$ be given as above. We will call both $X$ and $X_{0}$ $=$ 
$\{x_{0}\}$ $\cup $ $X,$ the (vector) labeling sets of $G$ (in $\Bbb{C}%
^{\oplus \,N}$).
\end{definition}

\strut

Let's fix $v_{0}$ $\in $ $V(G),$ and let $\deg _{out}(v_{0})$ $=$ $n$ $\leq $
$N.$ Define the subset $E_{out}^{v_{0}}$ of the edge set $E(G)$ by

\strut

\begin{center}
$E_{out}^{v_{0}}$ $\overset{def}{=}$ $\{e$ $\in $ $E(G)$ $:$ $e$ $=$ $v_{0}$ 
$e\}$ $\subseteq $ $E(G).$
\end{center}

\strut

i.e., $\deg _{out}(v)$ $=$ $\left| E_{out}^{v}\right| ,$ for all $v$ $\in $ $%
V(G).$ Then we decide the weights $\{(v_{0},$ $x_{1}),$ ..., $(v_{0},$ $%
x_{n})\},$ contained in the set $\{v_{0}\}$ $\times $ $X,$ on $%
E_{out}^{v_{0}}$. i.e.,

\strut

\begin{center}
$E_{out}^{v_{0}}$ $=$ $\left\{ e_{(v_{0},\text{ }x_{1})},...,e_{(v_{0},\text{
}x_{n})}:\left| 
\begin{array}{c}
\omega _{v_{0}}(e_{(v_{0},\text{ }x_{j})})=(v_{0},x_{j}), \\ 
\forall j=1,...,n
\end{array}
\right. \right\} ,$
\end{center}

\strut where

\begin{center}
$\omega _{v_{0}}$ $:$ $E_{out}^{v_{0}}$ $\rightarrow $ $\{v_{0}\}$ $\times $ 
$\{x_{1},$ ..., $x_{n}\}$
\end{center}

\strut

is the \emph{weighting process} putting the weight $(v_{0},$ $x_{j})$ to an
edge $e_{(v_{0},\text{ }x_{j})}.$

\strut

For any $v$ $\in $ $V(G),$ we can do the same process and we have the
corresponding weighting process $\omega _{v}$ and the corresponding set

$\strut $

\begin{center}
$E_{out}^{v}$ $=$ $\left\{ e_{(v,x_{1})},\text{ ..., }e_{(v,x_{k})}\left| 
\begin{array}{c}
x_{1},...,x_{k}\in X \\ 
\text{and }\omega _{v}(e_{(v,x_{j})})=(v,\text{ }x_{j}), \\ 
\forall \text{ }j=1,...,k,\text{ where} \\ 
k=\deg _{out}(v)\leq N
\end{array}
\right. \right\} .$
\end{center}

\strut

Notice that if $\deg _{out}(v)$ $=$ $0,$ then $E_{out}^{v}$ is empty, by the
fact that $\{v\}$ $\times $ $\varnothing $ $=$ $\varnothing ,$ where $%
\varnothing $ is the empty set. Notice also that, by the connectedness of $G$%
,

\strut

\begin{center}
$E(G)$ $=$ $\underset{v\in V(G)}{\sqcup }$ $E_{out}^{v},$
\end{center}

\strut

where ``$\sqcup $'' means the disjoint union. So, the each weighting
processes $(\omega _{v})_{v\in V(G)}$ extends to the \emph{weighting process
of }$E(G)$ by $X$:

\strut

\begin{center}
$\omega $ $:$ $E(G)$ $\rightarrow $ $V(G)$ $\times $ $X$
\end{center}

\strut by

\begin{center}
$\omega (e)$ $=$ $\omega (ve)$ $\overset{def}{=}$ $\omega _{v}(e),$
\end{center}

\strut \strut

for all $e$ $=$ $v$ $e$ $\in $ $E(G),$ with $v$ $\in $ $V(G).$

\strut

Now, let $G^{-1}$ be the shadow of $G.$ We can determine the similar
weighting process on $G^{-1}.$ Assume that the graph $G$ is weighted by the
weighting process $\omega $ for the labeling set $X.$ Define now the sets $%
-X $ and $-X_{0}$ by

\strut

\begin{center}
$-X$ $\overset{def}{=}$ $\{-x_{1},$ ..., $-x_{N}\}$ and $-X_{0}$ $=$ $%
\{x_{0}\}$ $\cup $ $(-X),$
\end{center}

where

\begin{center}
$x_{0}$ $=$ $\left( \underset{N\text{-times}}{\underbrace{0,\text{ ..., }0}}%
\right) $ and $-x_{j}$ $=$ $\left( 0,\text{ ..., }0,\text{ }\underset{j\text{%
-th}}{-1},\text{ }0,\text{ ..., }0\right) ,$
\end{center}

\strut

for all $j$ $=$ $1,$ ..., $N.$ i.e., the set $-X$ and $-X_{0}$ are also
contained in $\Bbb{C}^{\oplus \,N}.$ Now, consider the set $V(G^{-1})$ $%
\times $ $(-X)$ $=$ $V(G)$ $\times $ $(-X)$. Then, we can get the weighting
processes $(\omega _{v}^{-1})_{v\in V(G^{-1})}$ defined by

\strut

\begin{center}
$\omega _{v}^{-1}$ $:$ $E_{out}^{v^{-1}}$ $\rightarrow $ $\{v\}$ $\times $ $%
(-X)$
\end{center}

such that

\begin{center}
$\omega _{v}^{-1}(e^{-1})$ $=$ $\omega _{v}^{-1}(e^{-1}v)$ $\overset{def}{=}$
$(v,$ $-x_{i_{j}}),$
\end{center}

\strut

whenever $w_{v}(e)$ $=$ $(v$, $x_{i_{j}})$ $\in $ $\{v\}$ $\times $ $X.$
Therefore, we can define the \emph{weighting process }$\omega ^{-1}$\emph{\
of }$E(G^{-1})$ by

\strut

\begin{center}
$\omega ^{-1}$ $:$ $E(G^{-1})$ $\rightarrow $ $V(G^{-1})$ $\times $ $(-X)$
\end{center}

by

\begin{center}
$\omega ^{-1}(e^{-1})$ $=$ $\omega ^{-1}(e^{-1}v)$ $\overset{def}{=}$ $(v,$ $%
-x_{i_{j}}),$
\end{center}

\strut whenever

\begin{center}
$\omega (e)$ $=$ $(v,$ $x_{i_{j}})$ $\in $ $\{v\}$ $\times $ $X.$
\end{center}

\strut

Define the canonical projection

\strut

\begin{center}
$pr_{G}$ $:$ $V(G)$ $\times $ $X$ $\rightarrow $ $X$
\end{center}

by

\begin{center}
$pr_{G}\left( (v,\text{ }x)\right) $ $\overset{def}{=}$ $x,$ for all $(v,$ $%
x)$ $\in $ $V(G)\times X.$
\end{center}

\strut

Similarly, define the canonical projection

\strut

\begin{center}
$pr_{G^{-1}}$ $:$ $V(G)$ $\times $ $(-X)$ $\rightarrow $ $-X$
\end{center}

by

\begin{center}
$pr_{G^{-1}}\left( (v,\text{ }-x)\right) $ $\overset{def}{=}$ $-x,$ for all $%
(v,$ $-x)$ $\in $ $V(G^{-1})$ $\times $ $(-X).$
\end{center}

\strut \strut

\begin{definition}
Let $G$ be a locally finite connected directed graph equipped with its
weighting process $\omega $ of the edge set $E(G)$ by $V(G)$ $\times $ $X,$
where $X$ $=$ $\{1,$ ..., $N\},$ and $N$ $=$ $\max \{\deg _{out}(v)$ $:$ $v$ 
$\in $ $V(G)\}.$ Sometimes, we denote the graph $G$ equipped with the
weighting process $\omega $ by the pair $(G,$ $\omega ).$ Then the graph $%
(G, $ $\omega )$ is called the canonical weighted graph. Similarly, the
weighted graph $(G^{-1},$ $\omega ^{-1})$ equipped with the weighting
process $\omega ^{-1}$ of $E(G^{-1})$ by $(-X)$ $\times $ $V(G),$ where $%
G^{-1}$ is the shadow of $G,$ is called the canonical weighted shadow of $%
(G, $ $\omega ).$

\strut

Let $(G,$ $\omega )$ (resp. $(G^{-1},$ $\omega ^{-1})$) be the canonical
weighted graph (resp. the canonical weighted shadow of $(G,$ $\omega )$).
The image $pr_{G}$ $\circ $ $\omega $ of $E(G)$ (resp. $pr_{G^{-1}}$ $\circ $
$\omega ^{-1}$) contained in $X$ (resp. in $-X$) is called the canonical
labeling process of $(G,$ $\omega )$ (resp. of $(G^{-1},$ $\omega ^{-1})$).
\end{definition}

\strut

Let $(G,$ $\omega )$ be the canonical weighted graph with its canonical
weighted shadow $(G^{-1},$ $\omega ^{-1}).$ If $\Bbb{G}$ is the graph
groupoid of $G,$ then we can put the canonical weighting process, also
denoted by $\omega ,$ in terms of the canonical weighting process $\omega $
of the free semigroupoid $\Bbb{F}^{+}(G^{\symbol{94}})$ of the shadowed
graph $G^{\symbol{94}}$ of $G.$ Let $X^{*}$ be the collection of all finite
words in $X,$ containing the empty word $\emptyset _{X}$ in $X.$ Similarly,
we can determine the set $(-X)^{*}$ consisting of all finite words in $-X,$
containing the empty word $\emptyset _{-X}.$ Identify $\emptyset _{X}$ and $%
\emptyset _{-X}$ to $\emptyset _{*}$ in $(\pm X)^{*},$ the collection of all
words in $\pm X$ $\overset{def}{=}$ $X$ $\cup $ $(-X).$ i.e.,

\strut

\begin{center}
$\left( \pm X\right) ^{*}$ $\overset{def}{=}$ $\{\emptyset _{*}\}$ $\cup $ $%
\left( \underset{n=1}{\overset{\infty }{\cup }}\left\{ x_{i_{1}}\text{ ... }%
x_{i_{n}}\left| 
\begin{array}{c}
\text{ }x_{i_{1}},...,x_{i_{n}}\in \pm X,\text{ where} \\ 
\pm X=\{\pm x_{1},\text{ ..., }\pm x_{N}\}, \\ 
\text{with }N=\underset{v\in V(G)}{\max }\deg _{out}(v)
\end{array}
\right. \right\} \right) .$
\end{center}

\strut \strut \strut

Also define

\begin{center}
\strut $\pm X_{0}^{*}$ $\overset{def}{=}$ $\{x_{0}\}$ $\cup $ $\left( \pm
X\right) ^{*}.$
\end{center}

Define the map

\strut

\begin{center}
$\omega $ $:$ $FP(G^{\symbol{94}})$ $\cup $ $\{\emptyset \}$ $\rightarrow $ $%
V(G)^{2}\times $ $(\pm X_{0}^{*})$
\end{center}

by

\begin{center}
$\omega (\emptyset )$ $\overset{def}{=}$ $\emptyset _{*}$ $\overset{\text{%
regarded as}}{=}$ $\left( \varnothing ,\text{ }\emptyset _{*}\right) $
\end{center}

and

\begin{center}
$\omega \left( w\right) $ $=$ $\omega \left( v\text{ }e_{1}\text{ ... }%
e_{k}v^{\prime }\right) $ $\overset{def}{=}$ $\left( (v,\text{ }v^{\prime
}),\quad \underset{j=1}{\overset{k}{\Pi }}\widetilde{pr}\left( \omega _{\pm
}(e_{j})\right) \right) ,$
\end{center}

\strut

where $\varnothing $ is the empty set in $V(G)^{2}$ and $FP(G^{\symbol{94}})$
is the (non-reduced) finite path set of the shadowed graph $G^{\symbol{94}},$
and where

\strut

\begin{center}
$\widetilde{pr}$ $\in $ $\{pr_{G},$ $pr_{G^{-1}}\}$ and $\omega _{\pm }$ $%
\in $ $\{\omega ,$ $\omega ^{-1}\},$
\end{center}

whenever

\begin{center}
$w$ $=$ $v$ $w$ $v^{\prime }$ $=$ $e_{1}$ ... $e_{k},$
\end{center}

with

\begin{center}
$v,$ $v^{\prime }$ $\in $ $V(G^{\symbol{94}}),$ and $e_{1},$ ..., $e_{k}$ $%
\in $ $E(G^{\symbol{94}}),$
\end{center}

\strut

for $k$ $\in $ $\Bbb{N}.$ In the above, we have

\strut

\begin{center}
$\omega _{\pm }$ $:$ $E(G^{\symbol{94}})$ $=$ $E(G)$ $\sqcup $ $E(G^{-1})$ $%
\rightarrow $ $\pm X$
\end{center}

satisfies

\begin{center}
$\omega _{\pm }(e)$ $=$ $\left\{ 
\begin{array}{lll}
\omega (e) &  & \text{if }e\in E(G) \\ 
\omega ^{-1}(e) &  & \text{if }e\in E(G^{-1}),
\end{array}
\right. $
\end{center}

and

\begin{center}
$\widetilde{pr}$ $:$ $\left( V(G)\times X\right) \cup \left( V(G^{-1})\times
-X\right) \rightarrow X$ $\cup $ $(-X)$
\end{center}

satisfies

\begin{center}
$\widetilde{pr}\left( g\right) $ $=$ $\left\{ 
\begin{array}{lll}
pr_{G}(g) &  & \text{if }g\in V(G)\times X \\ 
pr_{G^{-1}}(g) &  & \text{if }g\in V(G)\times (-X).
\end{array}
\right. $
\end{center}

\strut

The above map $\omega $ from $FP(G^{\symbol{94}})$ $\cup $ $\{\emptyset \}$
can be extended to $\Bbb{F}^{+}(G^{\symbol{94}})$ into $V(G)^{2}$ $\times $ $%
\left( \pm X_{0}^{*}\right) $ by putting the weights on vertices:

\strut

\begin{center}
$\omega (v)$ $\overset{def}{=}$ $\left( (v,\text{ }v),\text{ }x_{0}\right) ,$
for all $v$ $\in $ $V(G).$
\end{center}

\strut \strut \strut

Such weighting process of the free semigroupoid $\Bbb{F}^{+}(G^{\symbol{94}%
}) $ of $G^{\symbol{94}}$ is called the \emph{canonical weighting process of 
}$\Bbb{F}^{+}(G^{\symbol{94}})$ or that of $\Bbb{G}$\emph{\ induced by the
canonical weighted graph} $(G,$ $\omega ).$

\strut

\textbf{Notation} From now, for convenience, we denote the above weighting
process of $\Bbb{F}^{+}(G^{\symbol{94}})$ simply by $\omega $. Remark that,
set-theoretically, the graph groupoid $\Bbb{G}$ is contained in the free
semigroupoid $\Bbb{F}^{+}(G^{\symbol{94}})$ of the shadowed graph $G^{%
\symbol{94}}$ of $G.$ So, we can determine the weighting process $\omega $
of $\Bbb{G}$ simply by restricting $\omega $ of $\Bbb{F}^{+}(G^{\symbol{94}%
}).$ i.e., the weighting process of $\Bbb{G}$ is $\omega $ $=$ $\omega $ $%
\mid _{\Bbb{G}}.$ We will use the notation $\omega $ for the weighting
processes of $\Bbb{F}^{+}(G^{\symbol{94}})$ and $\Bbb{G},$ alternatively. $%
\square $

\strut

Notice that the above weighting process $\omega $ of $\Bbb{F}^{+}(G^{\symbol{%
94}})$ represents the full admissibility conditions of $\Bbb{G}.$ i.e., the
weight $\omega (w)$ of $w$ $\in $ $\Bbb{G}$ contains its initial and
terminal vertices and the connecting pattern (admissibility of edges) of $w,$
in terms of labelings in $\pm X_{0}.$

\strut \strut

\begin{definition}
Let $(G,$ $\omega )$ be the canonical weighted graph with its canonical
weighted shadow $(G^{-1},$ $\omega ^{-1}).$ Let $\Bbb{G}$ be the graph
groupoid of $G$ and assume that $\omega $ is the canonical weighting process
of $\Bbb{G}$ induced by $(G,$ $\omega ).$ Then we call the pair $(\Bbb{G},$ $%
\omega )$ the labeled graph groupoid of $(G,$ $\omega ).$ We call the
processes $\widetilde{pr}$ $\circ $ $\omega $ on $(\Bbb{G},$ $\omega ),$ the
labeling process, where $\widetilde{pr}$ $\in $ $\{pr_{G},$ $pr_{G^{-1}}\},$
as in the previous paragraph.
\end{definition}

\strut

\textbf{Assumption} From now, if we mention (locally finite connected)
directed graphs, then they are automatically assumed to be the canonical
weighted graphs. Similarly, if we mention their graph groupoids, then they
are also automatically assumed to be the labeled graph groupoids. $\square $

\strut

\begin{example}
Let $G$ be the circulant graph with three vertices with

\strut

\begin{center}
$V(G)$ $=$ $\{v_{1},$ $v_{2},$ $v_{3}\}$
\end{center}

and

\begin{center}
$E(G)$ $=$ $\{e_{j}$ $=$ $v_{j}$ $e_{j}$ $v_{j+1}$ $:$ $j$ $=$ $1,$ $2,$ $3,$
with $v_{4}$ $\overset{def}{=}$ $v_{1}\}.$
\end{center}

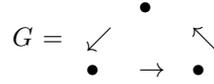
\strut \strut i.e.,

\begin{center}
$G$ $=$ $
\begin{array}{lll}
& \bullet &  \\ 
\swarrow &  & \nwarrow \\ 
\bullet & \rightarrow & \bullet
\end{array}
$
\end{center}

\strut

Then, since $\deg _{out}(v_{j})$ $=$ $1,$ for all $j$ $=$ $1,$ $2,$ $3,$ we
can get the labeling sets

\strut

\begin{center}
$X$ $=$ $\{x_{1}\}$ $=$ $\{1\}$ and $X_{0}$ $=$ $\{x_{0},$ $x_{1}\}$ $=$ $%
\{0,$ $1\},$
\end{center}

\strut

contained in $\Bbb{C}^{\oplus \,1}$ $=$ $\Bbb{C}.$ Relatively, we have that

\strut

\begin{center}
$-X$ $=$ $\{-x_{1}\}$ $=$ $\{-1\}$ and $-X_{0}$ $=$ $\{-x_{0},$ $-x_{1}\}$ $%
= $ $\{0,$ $-1\},$
\end{center}

\strut

for the shadow $G^{-1}$ of $G.$ Then we can get the weights

\strut

\begin{center}
$\omega (e_{j})$ $=$ $\left( (v_{j},\text{ }v_{j+1}),\text{ }1\right) ,$ for
all $j$ $=$ $1,$ $2,$ $3$
\end{center}

and

\begin{center}
$\omega ^{-1}(e_{j}^{-1})$ $=$ $\left( (v_{j+1},\text{ }v_{j}),\text{ }%
-1\right) ,$ for all $j$ $=$ $1,$ $2,$ $3,$
\end{center}

\strut

where $v_{4}$ $\overset{def}{=}$ $v_{1},$ in $V(G).$ We can construct the
labeled graph groupoid $(\Bbb{G},$ $\omega ).$ Let $w$ $=$ $e_{2}$ $e_{3}$ $%
e_{1}$ $\in $ $FP_{r}(G^{\symbol{94}}).$ Then we can get that

\strut

\begin{center}
$\omega (w)$ $=$ $((v_{2},$ $v_{2}),$ $111),$ and $\omega (w^{-1})$ $=$ $%
((v_{2},$ $v_{2}),$ $111).$
\end{center}

\strut

Take now the reduced finite path $y$ $=$ $e_{1}^{-1}e_{3}^{-1}$ $\in $ $%
FP_{r}(G^{\symbol{94}}).$ Then we can get

\strut

\begin{center}
$\omega (y)$ $=$ $((v_{2},$ $v_{3}),$ $(-1)(-1)).$

\strut

Suppose we have $w^{2}$ $=$ $e_{2}$ $e_{3}$ $e_{1}$ $e_{2}$ $e_{3}$ $e_{1}$
in $FP_{r}(G^{\symbol{94}}).$ Then the weight of $w^{2}$ is

\strut

$\omega (w^{2})$ $=$ $((v_{2},$ $v_{2}),$ $111111).$
\end{center}

\strut

Recall that $FP_{r}(G^{\symbol{94}})$ $=$ $FP(G)$ $\cup $ $FP(G^{-1})$ (See
[10] and [11]: Notice that, in general, $FP_{r}(G^{\symbol{94}})$ $%
\supsetneqq $ $FP(G)$ $\cup $ $FP(G^{-1})$).\strut
\end{example}

\strut

\textbf{Notation} For convenience, we re-define the set $X^{*}$ and $%
X_{0}^{*}$ of the labeling sets $X$ and $X_{0}$ as follows:

\strut

\begin{center}
$X^{*}$ $\overset{def}{=}$ $\{\emptyset _{*}\}$ $\cup $ $\left( \underset{k=1%
}{\overset{\infty }{\cup }}\left\{ (x_{i_{1}},\text{ }...,\text{ }%
x_{i_{k}})\in X^{k}\left| x_{i_{1}},\text{ }...,\text{ }x_{i_{n}}\in
X\right. \right\} \right) $
\end{center}

and

\begin{center}
$X_{0}^{*}$ $\overset{def}{=}$ $\{x_{0}\}$ $\cup $ $X^{*}.$
\end{center}

\strut

i.e., instead of writing the elements of $X^{*}$ as finite words in $X,$ we
will write them as a finite tuples having their entries in $X.$ $\square $%
\strut

\strut

\strut

\subsection{The Operation $\theta $ and $\omega _{+}$}

\strut

\strut

In this section, we will define a certain operation on the labeling set $\pm
X_{0}^{*}$ of the labeled graph groupoid $(\Bbb{G},$ $\omega )$ induced by
the canonical weighted graph $(G,$ $\omega ).$ Notice that the set $\pm
X_{0}^{*}$ is a subset of the $\Bbb{C}$-vector space $\Bbb{C}^{\oplus
\,\infty }.$ Let

\strut

\begin{center}
$N$ $\overset{def}{=}$ $\max \{\deg _{out}(v)$ $:$ $v$ $\in $ $V(G)\}$ $\in $
$\Bbb{N}.$
\end{center}

\strut

By regarding all elements $\pm X_{0}^{*}$ $\subset $ $\Bbb{C}^{\oplus
\,\infty },$ as finite tuples of their entries in $\Bbb{C}^{\oplus \,N},$ we
can define an operation $\theta $ from $\pm X_{0}^{*}$ into $\Bbb{C}^{\oplus
\,N}$ by

\strut

\begin{center}
$\theta \left( (t_{i_{1}},\text{ ..., }t_{i_{n}})\right) $ $\overset{def}{=}$
$\sum_{j=1}^{n}$ $t_{i_{j}}$ $\in $ $\Bbb{C}^{\oplus \,N},$
\end{center}

\strut \strut

for all $(t_{i_{1}},$ ..., $t_{i_{n}})$ $\in $ $\pm X_{0}^{*},$ for $n$ $\in 
$ $\Bbb{N},$ where

\strut

\begin{center}
$t_{i_{j}}$ $\overset{def}{=}$ $\left\{ 
\begin{array}{lll}
x_{i_{j}} &  & \text{if }t_{i_{j}}\in X_{0} \\ 
-x_{i_{j}} &  & \text{if }t_{i_{j}}\in -X.
\end{array}
\right. $
\end{center}

\strut

By help of $\theta ,$ we can define the operation $\omega _{+}$ by

$\strut $

\begin{center}
$\omega _{+}$ $:$ $V(G)^{2}$ $\times $ $\left( \pm X_{0}^{*}\right) $ $%
\rightarrow $ $V(G)^{2}$ $\times $ $\Bbb{C}^{\oplus \,N}$
\end{center}

by

\begin{center}
$
\begin{array}{ll}
\omega _{+}\left( (v,\text{ }v^{\prime }),\text{ }(t_{i_{1}},\text{ ..., }%
t_{i_{n}})\right) & \overset{def}{=}\left( (v,\text{ }v^{\prime }),\text{ }%
\theta \left( (t_{i_{1}},\text{ ..., }t_{i_{n}})\right) \right) \\ 
&  \\ 
& =\left( (v,\text{ }v^{\prime }),\text{ }\sum_{j=1}^{n}t_{i_{j}}\right) ,
\end{array}
$
\end{center}

\strut

for all $((v,$ $v^{\prime }),$ $(t_{i_{1}},$ ..., $t_{i_{n}}))$ $\in $ $%
V(G)^{2}$ $\times $ $(\pm X_{0}^{*})$. Later, we will use the above
operation $\omega _{+}$ to verify the amalgamated moments of the labeling
operators on right graph von Neumann algebras.

\strut

The operation $\theta $ provides the way to detect the reduction (RR) on $%
\Bbb{G}.$ Suppose we have two reduced finite path $w_{1}$ and $w_{2}$ having

\strut

\begin{center}
$\omega (w_{1})$ $=$ $\left( (v,\text{ }v^{\prime }),\text{ }X_{1}\right) $
and $\omega (w_{2})$ $=$ $\left( (v^{\prime },\text{ }v),\text{ }%
X_{2}\right) ,$
\end{center}

\strut

where $X_{1},$ $X_{2}$ $\in $ $\pm X_{0}^{*}$ $\setminus $ $\{\emptyset
_{*}\}.$ Since $w_{1}$ $=$ $w_{1}$ $v^{\prime }$ and $w_{2}$ $=$ $v^{\prime
} $ $w_{2},$ we can realize that $w_{1}$ and $w_{2}$ are admissible. i.e., $%
w_{1}$ $w_{2}$ $\neq $ $\emptyset $ in $\Bbb{G}.$ However, we do not know $%
w_{1}$ $w_{2}$ is either a vertex $v$ or a loop finite path $v$ $(w_{1}$ $%
w_{2})$ $v$ in $\Bbb{G}.$ Here, the operation $\theta $ (and hence the
operation $\omega _{+}$) provides the way to detect $w_{1}$ $w_{2}$ is
either the vertex $v$ or the loop finite path $w_{1}$ $w_{2}$ in $\Bbb{G}.$
In conclusion, we can get that:

\strut

\begin{proposition}
Let $w_{1}$ and $w_{2}$ be the reduced finite paths in $FP_{r}(G^{\symbol{94}%
})$ given in the previous paragraph. Then

\strut

\begin{center}
$\theta \left( \left( \widetilde{pr}\text{ }\circ \text{ }\omega \right)
(w_{1}w_{2})\right) $ $=$ $x_{0}$ in $\Bbb{C}^{\oplus \,N}$,
\end{center}

\strut or equivalently,

\begin{center}
$\omega _{+}\left( \omega (w_{1}w_{2})\right) $ $=$ $\left( (v,\text{ }v),%
\text{ }x_{0}\right) ,$
\end{center}

\strut

if and only if the element $w_{1}$ $w_{2}$ $\in $ $\Bbb{G}$ $\setminus $ $%
\{\emptyset \}$ is identical to the vertex $v.$
\end{proposition}

\strut

\begin{proof}
($\Rightarrow $) Suppose $\omega _{+}\left( \omega (w_{1}w_{2})\right) $ $=$ 
$\left( (v,\text{ }v),\text{ }x_{0}\right) .$ By hypothesis, we have that

\strut

\begin{center}
$\omega (w_{1}w_{2})$ $=$ $\left( (v,\text{ }v),\text{ }(t_{i_{1}},\text{
..., }t_{i_{n}},\text{ }t_{j_{1}},\text{ ..., }t_{j_{k}})\right) ,$
\end{center}

whenever

\begin{center}
$\omega (w_{1})$ $=$ $\left( (v,\text{ }v^{\prime }),\text{ }(t_{i_{1}},%
\text{ ..., }t_{i_{n}})\right) $
\end{center}

and

\begin{center}
$\omega (w_{2})$ $=$ $\left( (v^{\prime },\text{ }v),\text{ }(t_{j_{1}},%
\text{ ..., }t_{j_{k}})\right) .$
\end{center}

\strut

By assumption, $\theta \left( (t_{i_{1}},\text{ ..., }t_{i_{n}},\text{ }%
t_{j_{1}},\text{ ..., }t_{j_{k}})\right) $ $=$ $x_{0}.$ And since $w_{1}$
and $w_{2}$ are ``reduced'' finite paths in $\Bbb{G},$ we can verify that
(i) $n$ $=$ $k,$ and (ii) for any $t_{i_{p}}$, there should be the unique
entry $t_{j_{q}}$ $=$ $-t_{i_{p}},$ for $p,$ $q$ $=$ $1,$ ..., $n$ $=$ $k.$
This shows that, if $w_{1}$ $=$ $e_{i_{1}}$ ... $e_{i_{n}}$ and $w_{2}$ $=$ $%
e_{j_{1}}$ ... $e_{j_{n}}$ in $FP_{r}(G^{\symbol{94}}),$ with $e_{i_{1}},$
..., $e_{i_{n}},$ $e_{j_{1}},$ ..., $e_{j_{n}}$ $\in $ $E(G^{\symbol{94}}),$
then $e_{j_{q}}$ $=$ $e_{i_{p}}^{-1},$ uniquely, for $p,$ $q$ $=$ $1,$ ..., $%
n.$ Therefore, $w_{1}$ $w_{2}$ $=$ $v.$

\strut

($\Leftarrow $)\strut Suppose $w_{1}$ $w_{2}$ $=$ $v$ $\in $ $V(G).$ Then

\strut

\begin{center}
$\omega _{+}\left( \omega (w_{1}w_{2})\right) $ $=$ $\omega _{+}\left(
\omega (v)\right) $ $=$ $\omega _{+}\left( (v,\text{ }v),\text{ }%
x_{0}\right) $ $=$ $\left( (v,\text{ }v),\text{ }x_{0}\right) .$
\end{center}

\strut
\end{proof}

\strut

The above proposition shows how the labeling (or weighting) works on the
graph groupoid $\Bbb{G}.$ The labeling process on $\Bbb{G}$ gives a way to
check how admissibility works under (RR).

\strut

\strut

\subsection{Graph Automata}

\strut

\strut

In this section, we will construct an automaton induced by the labeling
graph groupoid $\Bbb{G}$ $=$ $(\Bbb{G},$ $\omega )$ of the canonical
weighted graph $G$ $=$ $(G,$ $\omega ).$ Let $X_{0}$ $=$ $\{x_{0},$ $x_{1},$
..., $x_{N}\}$ be the labeling set of the countable locally finite connected
directed graph $G$ contained in $\Bbb{C}^{\oplus \,N},$ where

\strut

\begin{center}
$N$ $\overset{def}{=}$ $\max \{\deg _{out}(v)$ $:$ $v$ $\in $ $V(G)\}.$
\end{center}

\strut \strut

Then we can create the automaton $\mathcal{A}_{G}$ $=$ $<\mathcal{X}_{0},$ $%
\Bbb{F}^{+}(G^{\symbol{94}}),$ $\varphi ,$ $\psi >$, where $G^{\symbol{94}}$ 
$=$ $G$ $\cup $ $G^{-1}$ is the shadowed graph of the given canonical
weighted graph $G$ and

\strut

\begin{center}
$\mathcal{X}_{0}$ $\overset{def}{=}$ $\{\emptyset _{G}\}$ $\cup $ $\left(
V(G)^{2}\times (\pm X_{0})\right) $
\end{center}

where

\begin{center}
$\pm X_{0}$ $=$ $X_{0}$ $\cup $ $(-X)$ $\subset $ $\Bbb{C}^{\oplus \,N},$
\end{center}

\strut and

\begin{center}
$\varphi \left( \left( (v_{1},\text{ }v_{2}),\text{ }t\right) ,\text{ }%
e\right) $ $\overset{def}{=}$ $\left\{ 
\begin{array}{ll}
\begin{array}{l}
\omega (e)
\end{array}
& 
\begin{array}{l}
\begin{array}{l}
\begin{array}{l}
\text{if }e=v_{2}\text{ }e\text{ in }E(G^{\symbol{94}}),\text{ equivalently,}
\\ 
\exists \text{ }v^{\prime }\in V(G)\text{ and }t^{\prime }\in \pm X_{0}\text{
s.t.,} \\ 
\quad \omega (e)=\left( (v_{2},\text{ }v^{\prime }),\text{ }t^{\prime
}\right) \neq \emptyset _{*}
\end{array}
\end{array}
\end{array}
\\ 
\,\,\,\,\emptyset _{G} & \quad \text{otherwise,}
\end{array}
\right. $
\end{center}

and

\begin{center}
$\psi \left( \left( (v_{1},\text{ }v_{2}),\text{ }t\right) ,\text{ }e\right) 
$ $\overset{def}{=}$ $\left\{ 
\begin{array}{ll}
\begin{array}{l}
e
\end{array}
& 
\begin{array}{l}
\text{if }e=v_{2}e\text{ in }E(G^{\symbol{94}})
\end{array}
\\ 
\emptyset & \quad \text{otherwise,}
\end{array}
\right. $
\end{center}

\strut

where $\emptyset _{G}$ is the empty element of $\mathcal{X}_{0}$ and $%
\emptyset $ is the empty element in $\Bbb{G}.$

\strut

\begin{definition}
We will say that the automaton $\mathcal{A}_{G}$ $=$ $<\mathcal{X}_{0},$ $%
\Bbb{F}^{+}(G^{\symbol{94}}),$ $\varphi ,$ $\psi >$ is the graph automaton
induced by the canonical weighted graph $G$ $=$ $(G,$ $\omega ).$ Sometimes,
we call $\varphi $ and $\psi ,$ the labeling map and the shift (or shifting
map), respectively.
\end{definition}

\strut

We can realize that the graph automaton $\mathcal{A}_{G}$ induced by $G$ is
identically same with the automaton $\mathcal{A}_{G^{-1}}$ induced by the
shadow $G^{-1}$ of $G.$ Now, fix an edge $e_{0}$ $\in $ $E(G^{\symbol{94}})$
and a reduced finite path $w$ $=$ $v_{1}e_{1}$ ... $e_{n}v_{n}^{\prime }$ $%
\in $ $FP_{r}(G^{\symbol{94}}),$ with $e_{j}$ $=$ $v_{j}$ $e_{j}$ $%
v_{j}^{\prime }$ $\in $ $E(G^{\symbol{94}}),$ and $v_{j},$ $v_{j}^{\prime }$ 
$\in $ $V(G^{\symbol{94}}),$ for all $j$ $=$ $1,$ ..., $n,$ for some $n$ $%
\in $ $\Bbb{N}.$ Then we can define $\varphi $ and $\psi $ inductively on
the free semigroupoid $\Bbb{F}^{+}(G^{\symbol{94}})$ of $G^{\symbol{94}}$ as
follows:

\strut

$\qquad \varphi \left( \left( (v_{1},\text{ }v_{n}^{\prime }),\text{ }%
(\omega (e_{1}),\text{ ..., }\omega (e_{n}))\right) ,\text{ }e_{0}\right) $

\strut

$\qquad \qquad \qquad =$ $\varphi \left( \omega (w),\text{ }e_{0}\right) $

\strut

$\qquad \qquad \qquad \overset{def}{=}$ $\left\{ 
\begin{array}{ll}
\omega (e_{0}) & \text{if }w\text{ }e_{0}\neq \emptyset \\ 
\emptyset _{G} & \text{otherwise,}
\end{array}
\right. $

\strut

and similarly,

\strut

$\qquad \psi \left( \left( (v_{1},\text{ }v_{n}^{\prime }),\text{ }(\omega
(e_{1}),\text{ ..., }\omega (e_{n}))\right) ,\text{ }e_{0}\right) $

\strut

$\qquad \qquad \qquad =$ $\psi \left( \omega (w),\text{ }e_{0}\right) $

\strut

$\qquad \qquad \qquad \overset{def}{=}$ $\left\{ 
\begin{array}{ll}
e_{0} & \text{if }w\text{ }e_{0}\neq \emptyset \\ 
\emptyset & \text{otherwise,}
\end{array}
\right. $

\strut $\strut $

for $n$ $\in $ $\Bbb{N}.$ Also, inductively, we have that, if $w^{\prime }$ $%
=$ $e_{1}^{\prime }$ ... $e_{m}^{\prime }$ $\in $ $FP_{r}(G^{\symbol{94}}),$
for $m$ $\in $ $\Bbb{N},$ then

\strut

$\qquad \varphi \left( \left( (v_{1},\text{ }v_{n}^{\prime }),\text{ }%
(\omega (e_{1}),\text{ ..., }\omega (e_{n}))\right) ,\text{ }w^{\prime
}\right) $

\strut

$\qquad \qquad \qquad \overset{def}{=}$ $\left\{ 
\begin{array}{ll}
\omega (e_{m}^{\prime }) & \text{if }w\text{ }w^{\prime }\neq \emptyset \\ 
\emptyset _{G} & \text{otherwise,}
\end{array}
\right. $

\strut

and similarly,

\strut

$\qquad \psi \left( \left( (v_{1},\text{ }v_{n}^{\prime }),\text{ }(\omega
(e_{1}),\text{ ..., }\omega (e_{n}))\right) ,\text{ }w^{\prime }\right) $

\strut

$\qquad \qquad \qquad \overset{def}{=}$ $\left\{ 
\begin{array}{ll}
w^{\prime } & \text{if }w\text{ }w^{\prime }\neq \emptyset \\ 
\emptyset & \text{otherwise.}
\end{array}
\right. $

\strut \strut \strut

Then, we can construct the automata actions $\{\mathcal{A}_{e}$ $:$ $e$ $\in 
$ $E(G^{\symbol{94}})\}$ of $\mathcal{A}_{G},$ acting on the set $V(G)^{2}$ $%
\times $ $(\pm X_{0}^{*})$ (consisting of all finite words in $\mathcal{X}%
_{0}$). Indeed, we can define the action $\mathcal{A}_{e}$ by

\strut

\begin{center}
$\mathcal{A}_{e}\left( \left( (v_{1},\text{ }v_{2}),\text{ }(x_{i_{1}},\text{
..., }x_{i_{n}})\right) \right) $ $\overset{def}{=}$ $\varphi \left( \left(
(v_{1},\text{ }v_{2}),\text{ }(x_{i_{1}},\text{ ..., }x_{i_{n}})\right) ,%
\text{ }e\right) $
\end{center}

\strut

for all $\left( (v_{1},v_{2}),\text{ }(x_{i_{1}},\text{ ..., }%
x_{i_{n}})\right) $ $\in $ $V(G)^{2}$ $\times $ $\left( \pm X_{0}^{*}\right)
,$ for all $e$ $\in $ $E(G^{\symbol{94}}).$ It is easy to be checked that
the actions $\mathcal{A}_{e}$'s satisfies

\strut

\begin{center}
$\mathcal{A}_{e_{1}}$ $\mathcal{A}_{e_{2}}$ $=$ $\mathcal{A}_{e_{2}e_{1}},$
for all $e_{1},$ $e_{2}$ $\in $ $E(G^{\symbol{94}}).$
\end{center}

\strut

Thus the actions $\{\mathcal{A}_{v}$ $:$ $v$ $\in $ $V(G^{\symbol{94}})\}$
are induced by $\mathcal{A}_{e}$'s since

\strut

\begin{center}
$\mathcal{A}_{v}$ $\overset{def}{=}$ $\mathcal{A}_{ee^{-1}}$ $=$ $\mathcal{A}%
_{e^{-1}}$ $\mathcal{A}_{e},$ whenever $e$ $=$ $v$ $e$ $\in $ $E(G^{\symbol{%
94}}).$
\end{center}

\strut

This guarantees that the automata actions $\{\mathcal{A}_{e}$ $:$ $e$ $\in $ 
$E(G^{\symbol{94}})\}$ on $V(G)^{2}$ $\times $ $\pm X_{0}^{*}$ generate the
groupoid which is groupoid-isomorphic to the graph groupoid $\Bbb{G}$ of $G.$
i.e., we can have that:

\strut

\begin{theorem}
The set $\{\mathcal{A}_{e}$ $:$ $e$ $\in $ $E(G^{\symbol{94}})\}$ of
automata actions generates the actions $\{\mathcal{A}_{w}$ $:$ $w$ $\in $ $%
\Bbb{F}^{+}(G^{\symbol{94}})\}$ of $\mathcal{A}_{G}.$ And the groupoid
generated by this set $\{\mathcal{A}_{e}$ $:$ $e$ $\in $ $E(G^{\symbol{94}%
})\}$ is groupoid-isomorphic to the graph groupoid $\Bbb{G}$ of $G.$ $%
\square $
\end{theorem}

\strut \strut \strut

\strut

\strut

\strut

\section{Fractaloids}

\strut \strut

\strut

As we assumed, all graphs in this Section are locally finite connected
countable directed graphs and they are canonically weighted by the labeling
set. \strut So, the corresponding graph groupoids are labeled.

\strut

\strut

\subsection{Fractaloids}

\strut

\strut \strut

As before, the standing assumption for our graphs $G$ are as follows:
locally finite, connected, and countable. The edges of $G$ are directed, and
they are assigned weights with the use of a labeling set. In this Section,
we will consider a certain special labeled graph groupoids. Let $\Bbb{G}$ be
the labeled graph groupoid induced by a locally finite directed graph $G,$
with the labeling set $\pm X_{0}$ $=$ $\{x_{0},$ $\pm x_{1},$ ..., $\pm
x_{N}\},$ and let the canonical weighting process $\omega $ of the free
semigroupoid $\Bbb{F}^{+}(G^{\symbol{94}})$ of the shadowed graph $G^{%
\symbol{94}}$ of $G,$ be given, where

$\strut $

\begin{center}
$N$ $\overset{def}{=}$ $\max \{\deg _{out}(v)$ $:$ $v$ $\in $ $V(G)\}.$
\end{center}

\strut

We observed that the canonical weighted graph $G$ creates the corresponding
graph automaton

$\strut $

\begin{center}
$\mathcal{A}_{G}$ $=$ $<\mathcal{X}_{0},$ $\Bbb{F}^{+}(G^{\symbol{94}}),$ $%
\varphi ,$ $\psi >,$
\end{center}

where

\begin{center}
$\mathcal{X}_{0}$ $=$ $\{\emptyset _{*}\}$ $\cup $ $\left( V(G)^{2}\times
(\pm X_{0})\right) ,$
\end{center}

\strut

and $\varphi $ is the labeling map and $\psi $ is the shifting map. We
already observed that $\Bbb{G}$ acts on a $\Bbb{F}^{+}(G^{\symbol{94}})$-set 
$\mathcal{X}_{0}^{*},$ where $\mathcal{X}_{0}^{*}$ is the collection of all
finite (non-reduced) words in $\mathcal{X}_{0}.$

\strut

Recall that two countable directed graphs $G_{1}$ and $G_{2}$ are \emph{%
graph-isomorphic}, if there exists a bijection

$\strut $

\begin{center}
$g$ $:$ $V(G_{1})$ $\cup $ $E(G_{1})$ $\rightarrow $ $V(G_{2})$ $\cup $ $%
E(G_{2}),$
\end{center}

\strut

such that (i) $g\left( V(G_{1})\right) $ $=$ $V(G_{2})$ and $g\left(
E(G_{1})\right) $ $=$ $E(G_{2}),$ and (ii) $g\left( e\right) $ $=$ $g(v_{1}$ 
$e$ $v_{2})$ $=$ $g(v_{1})$ $g(e)$ $g(v_{2})$ in $E(G_{2}),$ whenever $e$ $=$
$v_{1}$ $e$ $v_{2}$ $\in $ $E(G_{1}),$ with $v_{1},$ $v_{2}$ $\in $ $%
V(G_{1}).$

\strut

In [10] and [11], we showed that if two graphs $G_{1}$ and $G_{2}$ are
graph-isomorphic, then the corresponding graph grouopoids $\Bbb{G}_{1}$ and $%
\Bbb{G}_{1}$ are groupoid-isomorphic. More generally, if two graphs $G_{1}$
and $G_{2}$ have the graph-isomorphic shadowed graphs, then the graph
groupoids $\Bbb{G}_{1}$ and $\Bbb{G}_{1}$ are groupoid-isomorphic. Also, we
showed that if two graph groupoids $\Bbb{G}_{1}$ and $\Bbb{G}_{2}$ are
groupoid-isomorphic, then the (left) graph von Neumann algebras $\overline{%
\Bbb{C}[\Bbb{G}_{1}]}^{w}$ and $\overline{\Bbb{C}[\Bbb{G}_{2}]}^{w}$ are $*$%
-isomorphic, as $W^{*}$-subalgebras in the operator algebra $B(H_{G_{1}})$ $%
\overset{\text{identically same}}{=}$ $B(H_{G_{2}}).$ Therefore, if two
graphs $G_{1}$ and $G_{2}$ have the graph-isomorphic shadowed graphs, then
the right graph von Neumann algebras $\Bbb{M}_{G_{1}}$ and $\Bbb{M}_{G_{2}}$
are $*$-isomorphic, too.

\strut

Recall also that we say a directed graph $G$ is a (directed)\emph{\ tree} if
this graph $G$ is connected and it has no loop finite paths in $\Bbb{F}%
^{+}(G)$. Also, we say that a directed tree is \emph{rooted}, if we can
find-and-fix a vertex $v_{0}$ of $G,$ with $\deg _{in}(v_{0})$ $=$ $0.$ This
fixed vertex $v_{0}$ is called the\emph{\ root of} $G.$ For instance, a graph

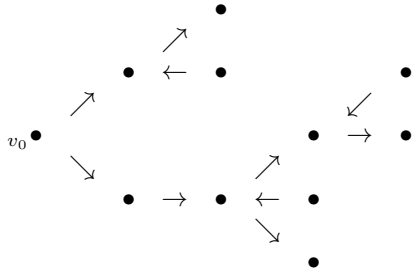
\strut

\begin{center}
$
\begin{array}{lllllllll}
&  &  &  & \bullet &  &  &  &  \\ 
&  &  & \nearrow &  &  &  &  &  \\ 
&  & \bullet & \leftarrow & \bullet &  &  &  & \bullet \\ 
& \nearrow &  &  &  &  &  & \swarrow &  \\ 
_{v_{0}}\bullet &  &  &  &  &  & \bullet & \rightarrow & \bullet \\ 
& \searrow &  &  &  & \nearrow &  &  &  \\ 
&  & \bullet & \rightarrow & \bullet & \leftarrow & \bullet &  &  \\ 
&  &  &  &  & \searrow &  &  &  \\ 
&  &  &  &  &  & \bullet &  & 
\end{array}
$
\end{center}

\strut

is a rooted tree with its (fixed) root $v_{0}$. A rooted tree $G$ is \emph{%
one-flow}, if the directions of edges oriented only one way from the root $%
v_{0}.$ For example, a graph

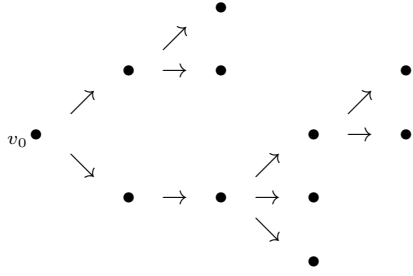
\strut

\begin{center}
$
\begin{array}{lllllllll}
&  &  &  & \bullet &  &  &  &  \\ 
&  &  & \nearrow &  &  &  &  &  \\ 
&  & \bullet & \rightarrow & \bullet &  &  &  & \bullet \\ 
& \nearrow &  &  &  &  &  & \nearrow &  \\ 
_{v_{0}}\bullet &  &  &  &  &  & \bullet & \rightarrow & \bullet \\ 
& \searrow &  &  &  & \nearrow &  &  &  \\ 
&  & \bullet & \rightarrow & \bullet & \rightarrow & \bullet &  &  \\ 
&  &  &  &  & \searrow &  &  &  \\ 
&  &  &  &  &  & \bullet &  & 
\end{array}
$
\end{center}

\strut

is a one-flow rooted tree with its root $v_{0}$. A one-flow rooted tree $G$
is said to be \emph{growing} if $G$ is an infinitely countable directed
graph. Finally, we will say that a one-flow growing rooted tree $G$ is \emph{%
regular} if the out-degrees of all vertices are identical. For example, a
graph

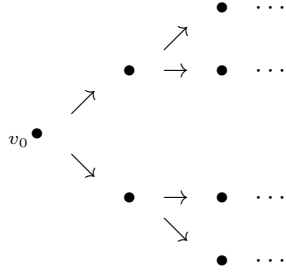
\strut

\begin{center}
$
\begin{array}{llllll}
&  &  &  &  &  \\ 
&  &  &  & \bullet & \cdots \\ 
&  &  & \nearrow &  &  \\ 
&  & \bullet & \rightarrow & \bullet & \cdots \\ 
& \nearrow &  &  &  &  \\ 
_{v_{0}}\bullet &  &  &  &  &  \\ 
& \searrow &  &  &  &  \\ 
&  & \bullet & \rightarrow & \bullet & \cdots \\ 
&  &  & \searrow &  &  \\ 
&  &  &  & \bullet & \cdots
\end{array}
$
\end{center}

\strut

is a regular one-flow growing rooted tree. In particular, if $\deg _{out}(v)$
$=$ $N,$ for all vertices $v$, then this regular one-flow growing rooted
tree is called \emph{the} $N$-\emph{regular tree}. The very above example is
the $2$-regular tree.

\strut

Let $N$ be the maximal out-degree of the graph $G,$ and let $\mathcal{T}%
_{2N} $ be the $2N$-regular tree. Then the automata actions $\{\mathcal{A}%
_{w}$ $:$ $w$ $\in $ $\Bbb{F}^{+}(G^{\symbol{94}})\}$ of the graph automaton 
$\mathcal{A}_{G}$ acts on this $2N$-regular tree $\mathcal{T}_{2N}.$ For
instance, for the set $\{\mathcal{A}_{w}\}_{w\in \Bbb{F}^{+}(G^{\symbol{94}%
})}$ of automata actions of the graph automaton $\mathcal{A}_{G},$ we can
create a one-flow growing rooted tree $\mathcal{T}_{G}$ having its
arbitrarily fixed root $\omega (v_{0})$ $=$ $\left( (v_{0},v_{0}),\text{ }%
x_{0}\right) $ $\in $ $V(G)^{2}$ $\times $ $(\pm X_{0}^{*})$, where

\strut

\begin{center}
$V(\mathcal{T}_{G})$ $\overset{def}{=}$ $\left\{ \varphi (X,\text{ }e)\left| 
\begin{array}{c}
X\in V(G)^{2}\times (\pm X_{0}^{*}), \\ 
e\in E(G^{\symbol{94}})
\end{array}
\right. \right\} $
\end{center}

and

\begin{center}
$E(\mathcal{T}_{G})$ $\overset{def}{=}$ $\left\{ \psi (X,e)\left| 
\begin{array}{c}
X\in V(G)^{2}\times (\pm X_{0}^{*}), \\ 
e\in E(G^{\symbol{94}})
\end{array}
\right. \right\} ,$
\end{center}

\strut

where $G^{\symbol{94}}$ is the shadowed graph of $G.$

\strut

\textbf{Observation and Notation} Notice that, by the connectedness of the
shadowed graph $G^{\symbol{94}}$ of the connected graph $G,$ we can fix any
weight $\omega (v)$ of $v$ $\in $ $V(G^{\symbol{94}}),$ as the root of the
tree $\mathcal{T}_{G}$. Suppose $\mathcal{T}_{v_{1}}$ and $\mathcal{T}%
_{v_{2}}$ are the one-flow growing trees with their roots $((v_{1},$ $%
v_{1}), $ $x_{0})$ and $((v_{2},$ $v_{2}),$ $x_{0}),$ respectively. Then, $%
\mathcal{T}_{v_{i}}$ is embedded in $\mathcal{T}_{v_{j}},$ as a
full-subgraph (See the definition of full-subgraphs below), whenever $i$ $%
\neq $ $j$ $\in $ $\{1,$ $2\}.$ In general, the graphs $\mathcal{T}_{v_{1}}$
and $\mathcal{T}_{v_{2}}$ have no graph-isomorphic relation, but they are
embedded from each other, by the connectedness of $G^{\symbol{94}}.$ Since $%
v_{1}$ and $v_{2}$ are arbitrary, we can consider only one choice $\mathcal{T%
}_{v},$ for $v$ $\in $ $V(G),$ as a candidate of the tree, denoted by $%
\mathcal{T}_{G},$ as the tree where the automata actions act. i.e., whenever
we choose one tree $\mathcal{T}_{v},$ for $v$ $\in $ $V(G),$ then the trees $%
\mathcal{T}_{v^{\prime }}$'s are embedded in $\mathcal{T}_{v},$ for all $%
v^{\prime }$ $\in $ $V(G).$ Without loss of generality, if we write $%
\mathcal{T}_{G}$ from now, then it means a tree $\mathcal{T}_{v},$ for a
fixed $v$ $\in $ $V(G).$ Remark that, the tree $\mathcal{T}_{G}$ has its
root $\omega (v)$ if and only if $\mathcal{T}_{G}$ $=$ $\mathcal{T}_{v}.$ $%
\square $

\strut

We can easily check that

\strut

\begin{center}
$FP(\mathcal{T}_{G})$ $=$ $\left\{ \psi (X,\text{ }w)\left| 
\begin{array}{c}
X\in V(G)^{2}\times (\pm X_{0}^{*}), \\ 
w\in FP(G^{\symbol{94}})
\end{array}
\right. \right\} ,$
\end{center}

and hence

\begin{center}
$
\begin{array}{ll}
\Bbb{F}^{+}(\mathcal{T}_{G})= & \{\emptyset _{*}\}\cup \left\{ \varphi (X,%
\text{ }e)\left| 
\begin{array}{c}
X\in V(G)^{2}\times (\pm X_{0}^{*}), \\ 
e\in E(G^{\symbol{94}})
\end{array}
\right. \right\} \\ 
&  \\ 
& \quad \quad \;\cup \left\{ \psi (X,\text{ }w)\left| 
\begin{array}{c}
X\in V(G)^{2}\times (\pm X_{0}^{*}), \\ 
w\in FP(G^{\symbol{94}})
\end{array}
\right. \right\} .
\end{array}
$
\end{center}

\strut \strut \strut

By the connectedness of the graph $G,$ and by the definition of the automata
actions, every finite paths in $FP(G^{\symbol{94}})$ is embedded in the
finite paths of the tree $\mathcal{T}_{G},$ via the automata actions. Then,
we can construct the full-subgraphs $\{\mathcal{T}_{w}$ $:$ $w$ $\in $ $%
FP(G^{\symbol{94}})\}$ of $\mathcal{T}_{G},$ where $\mathcal{T}_{w}$'s are
the one-flow growing rooted tree with their roots $\omega \left( w\right) ,$
for all $w$ $\in $ $FP(G^{\symbol{94}}).$

\strut

Recall that, we say that a countable directed graph $G_{1}$ is a \emph{%
full-subgraph} of a countable directed graph $G_{2},$ if

\strut

\begin{center}
$E(G_{1})$ $\subseteq $ $E(G_{2})$
\end{center}

and

\begin{center}
$V(G_{1})$ $=$ $\left\{ v\in V(G_{2})\left| 
\begin{array}{c}
e=ve\text{ or }e=ev, \\ 
\forall e\in E(G_{1})
\end{array}
\right. \right\} .$
\end{center}

\strut

Notice the difference between full-subgraphs and subgraphs. We say that $%
G_{1}$ is a \emph{subgraph} of $G_{2},$ if

\strut

\begin{center}
$V(G_{1})$ $\subseteq $ $V(G_{2})$
\end{center}

and

\begin{center}
$E(G_{1})$ $=$ $\left\{ e\in E(G_{2})\left| 
\begin{array}{c}
e=v\text{ }e\text{ }v^{\prime }, \\ 
\forall v,\text{ }v^{\prime }\in V(G_{2})
\end{array}
\right. \right\} .$
\end{center}

\strut

Every subgraph is a full-subgraph, but the converse does not hold true.\strut

\strut

\begin{definition}
Let $\mathcal{T}_{G}$ be the above one-flow growing rooted tree, where the
automata actions of the graph automaton $\mathcal{A}_{G}$ act. We say that
the tree $\mathcal{T}_{G}$ is the $\mathcal{A}_{G}$-tree, which is a
full-subgraph of the $2N$-regular tree $\mathcal{T}_{2N}.$ The
full-subgraphs $\mathcal{T}_{w}$'s for $w$ $\in $ $FP(G^{\symbol{94}})$ of $%
\mathcal{T}_{G}$ are called the $w$-parts of $\mathcal{T}_{G},$ which are
also full-subgraphs of $\mathcal{T}_{2N}.$
\end{definition}

\strut \strut

\strut The important thing is now that the $w$-parts $\mathcal{T}_{w}$'s are
embedded in the $2N$-regular tree $\mathcal{T}_{2N},$ and $\mathcal{T}%
_{w^{\prime }}$'s are embedded in $\mathcal{T}_{w}$'s, whenever

\strut

\begin{center}
$\varphi \left( \omega (v),\text{ }w^{\prime }\right) $ $=$ $\left( (v,\text{
}v^{\prime }),\text{ }X^{\prime }\right) ,$
\end{center}

where

\begin{center}
$\varphi \left( \omega (v),\text{ }w\right) $ $=$ $\left( (v,\text{ }%
v^{\prime \prime }),\text{ }X\right) $ and $X^{\prime }$ $=$ $(X,$ $%
X^{\prime \prime }),$
\end{center}

\strut

for some $X^{\prime \prime }$ $\in $ $(\pm X_{0}^{*}),$ where $\omega (v)$
is the root of $\mathcal{T}_{G}.$

\strut

\begin{remark}
The construction of the $\mathcal{A}_{G}$-tree $\mathcal{T}_{G}$ is nothing
but the rearrangement of the finite paths in the free semigroupoid $\Bbb{F}%
^{+}(G^{\symbol{94}})$ of the shadowed graph $G^{\symbol{94}}$ inside the $%
2N $-regular tree $\mathcal{T}_{2N},$ up to the admissibility on $\Bbb{F}%
^{+}(G^{\symbol{94}})$. Notice that, in fact, the $\mathcal{A}_{G}$-tree $%
\mathcal{T}_{G}$ contains the information about the vertices in $\Bbb{F}%
^{+}(G^{\symbol{94}}),$ too, since the vertices of $\mathcal{T}_{G}$ are
contained in $V(G)^{2}$ $\times $ $(\pm X_{0}^{*}).$ Remark that, by the
connectedness of $G^{\symbol{94}},$ the $w$-parts $\mathcal{T}_{w}$'s of $%
\mathcal{T}_{G}$ for $w$ $\in $ $FP(G^{\symbol{94}})$ are well-constructed
as a one-flow growing tree with their roots $\varphi (\omega (v),$ $w)$,
where $\omega (v)$ is the root of $\mathcal{T}_{G}.$ Moreover, each tree $%
\mathcal{T}_{w}$ is embedded in the other trees $\mathcal{T}_{G}.$
\end{remark}

\strut \strut \strut \strut \strut

Notice that, even though $\left| \Bbb{G}\right| $ $<$ $\infty ,$ $\left| 
\Bbb{F}^{+}(G^{\symbol{94}})\right| $ $=$ $\infty ,$ whenever $\left|
E(G)\right| $ $\geq $ $1.$ By identifying their roots, the $\mathcal{A}_{G}$%
-tree $\mathcal{T}_{G}$ is a full-subgraph of the $2N$-regular tree $%
\mathcal{T}_{2N}.$ Therefore, we can conclude that:

\strut

\begin{lemma}
The automata actions $\{\mathcal{A}_{w}$ $:$ $w$ $\in $ $\Bbb{F}^{+}(G^{%
\symbol{94}})\}$ of the graph automaton $\mathcal{A}_{G},$ induced by the
given graph $G,$ act on the $2N$-regular tree $\mathcal{T}_{2N}.$
\end{lemma}

\strut

\begin{proof}
By definition, the automata actions $\{\mathcal{A}_{w}\}_{w\in \Bbb{F}%
^{+}(G^{\symbol{94}})}$ act on the $\mathcal{A}_{G}$-tree $\mathcal{T}_{G}.$
And the tree $\mathcal{T}_{G}$ is a full-subgraph of $\mathcal{T}_{2N}.$
\end{proof}

\strut \strut

Let $\mathcal{A}_{G}$ $=$ $<\mathcal{X}_{0},$ $\Bbb{F}^{+}(G^{\symbol{94}}),$
$\varphi ,$ $\psi >$ be the graph automaton induced by the canonical
weighted graph $G$ $=$ $(G,$ $\omega ).$ In Section 3.3, we showed that the
groupoid $\Bbb{G}(\mathcal{A}_{G})$ generated by the automata actions $\{%
\mathcal{A}_{w}$ $:$ $w$ $\in $ $\Bbb{F}^{+}(G^{\symbol{94}})\}$ is
groupoid-isomorphic to the graph groupoid $\Bbb{G}$ of $G.$ This is the
graph-groupoid version of the fact that: if $\mathcal{A}$ $=$ $<X^{\pm },$ $%
X^{\pm *},$ $\varphi ,$ $\psi >$ is an automaton, where $X$ is the generator
set, with $X^{\pm }$ $=$ $\{x^{\pm 1}$ $:$ $x$ $\in $ $X\},$ of a group $%
\Gamma $ and $X^{*}$ is the free monoid of the group $\Gamma $, which is
monoid-homomorphic to $\Gamma ,$ then the automata actions $\{\mathcal{A}%
_{y} $ $:$ $y$ $\in $ $X^{\pm *}\}$ generate the automata group $\Gamma (%
\mathcal{A}),$ which is group-isomorphic to $\Gamma $.

\strut

In fact, this fact is not so interesting in Automata Theory and Groupoid
Theory, however, for our works, it plays a key role, since we can consider
the labeled graph groupoid $\Bbb{G}$ and the automata groupoid $\Bbb{G}(%
\mathcal{A}_{G}),$ alternatively, as groupoid-isomorphic objects. Also, it
provides the technique to see the (reduced or non-reduced) finite paths of
the shadowed graph $G^{\symbol{94}}$ as those of $\mathcal{T}_{2N}$ (as
embedded elements).

\strut

\begin{definition}
Let $G$, $\Bbb{F}^{+}(G^{\symbol{94}}),$ and $\Bbb{G}$ be given as before
and let $\mathcal{A}_{G}$ $=$ $<\mathcal{X}_{0},$ $\Bbb{F}^{+}(G^{\symbol{94}%
}),$ $\varphi ,$ $\psi >$ be the graph automaton induced by $G,$ acting on
the $2N$-regular tree $\mathcal{T}_{2N}$. For any fixed $w$ $\in $ $FP(G^{%
\symbol{94}}),$ the tree $\mathcal{T}_{w}$ is the $w$-part of the $\mathcal{A%
}_{G}$-tree $\mathcal{T}_{G},$ with its root $\varphi (\omega (v),$ $w),$
where $\omega (v)$ is the root of the $\mathcal{A}_{G}$-tree $\mathcal{T}%
_{G}.$ Let $\Bbb{G(\mathcal{T}}_{w}\Bbb{)}$ be the groupoid generated by the
actions $\mathcal{A}_{y}$'s acting only on $\mathcal{T}_{w}$. If $\Bbb{G}(%
\mathcal{T}_{w})$'s are groupoid-isomorphic to $\Bbb{G}(\mathcal{T}_{G}),$
for all $w$ $\in $ $FP(G^{\symbol{94}}),$ then we say that the groupoid $%
\Bbb{G}(\mathcal{A}_{G})$ is a fractaloid. Equivalently, we say that the
(labeled) graph groupoid $\Bbb{G}$ is a fractaloid.
\end{definition}

\strut

Readers can understand the above definition of fractaloids as the
graph-groupoid version of the fractal groups (See [1]). The following
theorem provides the graph-theoretical characterization of fractaloids.

\strut

\begin{theorem}
Let $G$ be a canonical weighted graph with its labeled graph groupoid $\Bbb{G%
},$ and let $\mathcal{A}_{G}$ be the graph automaton induced by $G$ and $%
\mathcal{T}_{G},$ the $\mathcal{A}_{G}$-tree. Every $w$-part $\mathcal{T}%
_{w} $ of $\mathcal{T}_{G}$ is graph-isomorphic to $\mathcal{T}_{G},$ for
all $w$ $\in $ $FP(G^{\symbol{94}}),$ if and only if $\Bbb{G}$ is a
fractaloid.
\end{theorem}

\strut

\begin{proof}
($\Leftarrow $) By definition, if the graph groupoid $\Bbb{G}$ $\overset{%
\text{Groupoid}}{=}$ $\Bbb{G}(\mathcal{A}_{G})$ is a fractaloid, then, for
any $w$-parts $\mathcal{T}_{w}$, the automata actions acting only on $%
\mathcal{T}_{w}$ generate the groupoids $\Bbb{G}(\mathcal{T}_{w}),$ which is
groupoid-isomorphic to $\Bbb{G},$ where ``$\overset{\text{Groupoid}}{=}$''
means ``being groupoid-isomorphic''. Assume now that there exists a finite
path $w_{0}$ $\in $ $FP(G^{\symbol{94}})$ such that the $w_{0}$-part $%
\mathcal{T}_{w_{0}}$ is not graph-isomorphic to the $\mathcal{A}_{G}$-tree $%
\mathcal{T}_{G}.$ Then, clearly, the groupoids $\Bbb{G(\mathcal{T}}_{w_{0}}%
\Bbb{)}$ and $\Bbb{G}(\mathcal{T}_{G})$ are not groupoid-isomorphic. This
contradicts our assumption.

\strut

($\Rightarrow $) Assume now that every $w$-part $\mathcal{T}_{w}$ of the $%
\mathcal{A}_{G}$-tree $\mathcal{T}_{G}$ is graph-isomorphic to $\mathcal{T}%
_{G}$, for all $w$ $\in $ $FP(G^{\symbol{94}}).$ This shows that the
groupoids $\Bbb{G}(\mathcal{T}_{w})$'s and $\Bbb{G}(\mathcal{T}_{G})$ are
groupoid-isomorphic, for all $w$ $\in $ $FP(G^{\symbol{94}}),$ since the
automata actions $\{\mathcal{A}_{y}\}$ acting on $\mathcal{T}_{w}$ and the
automata action $\{\mathcal{A}_{d}\}$ acting on $\mathcal{T}_{G}$ are
equivalent. Since $w$ is arbitrary in $FP(G\symbol{94})$, the graph groupoid 
$\Bbb{G}$ is a fractaloid.
\end{proof}

\strut

Let $G$ be a canonical weighted graph with its labeled graph groupoid $\Bbb{G%
},$ and assume that the automata actions $\{\mathcal{A}_{w}$ $:$ $w$ $\in $ $%
FP(G^{\symbol{94}})\}$ of the graph automaton $\mathcal{A}_{G}$\emph{\ act
fully on the }$2N$\emph{-regular tree} $\mathcal{T}_{2N},$ in the sense that
the $\mathcal{A}_{G}$-tree $\mathcal{T}_{G},$ which is a full-subgraph of $%
\mathcal{T}_{2N},$ is identical to $\mathcal{T}_{2N}.$ i.e., the automata
actions of $\mathcal{A}_{G}$ act fully on $\mathcal{T}_{2N},$ if $\mathcal{T}%
_{G}$ $\overset{\text{Graph}}{=}$ $\mathcal{T}_{2N}$.

\strut

\begin{corollary}
Let $G$ be a canonical weighted graph with its labeled graph groupoid $\Bbb{G%
},$ and let $\mathcal{A}_{G}$ be the graph automaton induced by $G.$ If the
automata actions $\{\mathcal{A}_{w}$ $:$ $w$ $\in $ $FP(G^{\symbol{94}})\}$
of $\mathcal{A}_{G}$ act fully on the $2N$-regular tree $\mathcal{T}_{2N},$
then $\Bbb{G}$ is a fractaloid.
\end{corollary}

\strut

\begin{proof}
Suppose the automata actions act fully on the $2N$-regular tree $\mathcal{T}%
_{2N}.$ Then, by definition, the $\mathcal{A}_{G}$-tree $\mathcal{T}_{G}$ is
identical to $\mathcal{T}_{2N}.$ This shows that, every $w$-part $\mathcal{T}%
_{w}$ of $\mathcal{T}_{G}$ is graph-isomorphic to $\mathcal{T}_{2N},$ for
all $w$ $\in $ $FP(G^{\symbol{94}}).$ Therefore, by the previous theorem,
the graph groupoid $\Bbb{G}$ of $G$ is a fractaloid.
\end{proof}

\strut

Also, the above corollary provide an easy technique to construct
fractaloidal examples. For example, the one-flow circulant graphs and
one-vertex-multi-loop-edge graphs have their graph groupoids which are
fractaloids, as connected ``finite'' directed graphs (See Section 4.2).
Recall that a graph $G$ is \emph{finite}, if $\left| V(G)\right| $ $<$ $%
\infty $ and $\left| E(G)\right| $ $<$ $\infty $.

\strut

How about the converse of the previous corollary? Let $\mathcal{T}$ be a
one-flow growing infinite rooted tree with its root $v_{0}.$ Then the
terminal vertices of the incident edges of $v_{0}$ are called the $1$-st
level of $\mathcal{T}.$ Similarly, if $v_{1}$, ..., $v_{t}$ are the vertices
in the 1-st level, then the terminal vertices of the incident edges of $%
v_{1},$ ..., $v_{t}$ are called the 2-nd level of $\mathcal{T}.$
Inductively, we can decide the $n$-th level of $\mathcal{T},$ for $n$ $\in $ 
$\Bbb{N}.$ For instance,

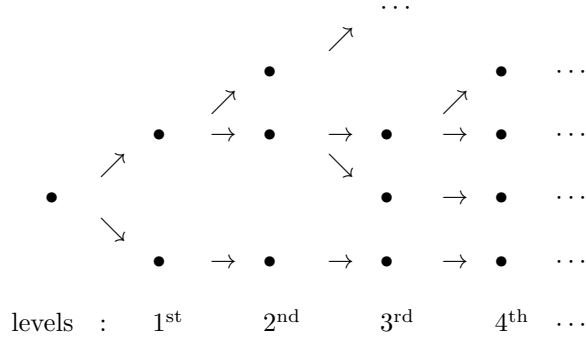
\strut

\begin{center}
$
\begin{array}{llllllllll}
&  &  &  &  &  & \cdots &  &  &  \\ 
&  &  &  &  & \nearrow &  &  &  &  \\ 
&  &  &  & \bullet &  &  &  & \bullet & \cdots \\ 
&  &  & \nearrow &  &  &  & \nearrow &  &  \\ 
&  & \bullet & \rightarrow & \bullet & \rightarrow & \bullet & \rightarrow & 
\bullet & \cdots \\ 
& \nearrow &  &  &  & \searrow &  &  &  &  \\ 
\quad \;\bullet &  &  &  &  &  & \bullet & \rightarrow & \bullet & \cdots \\ 
& \searrow &  &  &  &  &  &  &  &  \\ 
&  & \bullet & \rightarrow & \bullet & \rightarrow & \bullet & \rightarrow & 
\bullet & \cdots \\ 
&  &  &  &  &  &  &  &  &  \\ 
\text{levels} & : & 1^{\text{st}} &  & 2^{\text{nd}} &  & 3^{\text{rd}} &  & 
4^{\text{th}} & \cdots
\end{array}
$
\end{center}

\strut

\begin{theorem}
If $\Bbb{G}$ is a fractaloid, then the automata actions $\{\mathcal{A}_{w}$ $%
:$ $w$ $\in $ $FP(G^{\symbol{94}})\}$ act fully on the $\mathcal{A}_{G}$%
-tree $\mathcal{T}_{G},$ and $\mathcal{T}_{G}$ is graph-isomorphic to the $%
2N $-regular tree $\mathcal{T}_{2N},$ by identifying their roots, where $N$ $%
=$ $\max \{$ $\deg _{out}(v)$ $:$ $v$ $\in $ $V(G)\}.$
\end{theorem}

\strut

\begin{proof}
By the graph-theoretical characterization of fractaloids, the labeled graph
groupoid $\Bbb{G}$ is a fractaloid if and only if every $w$-part $\mathcal{T}%
_{w}$ of the $\mathcal{A}_{G}$-tree $\mathcal{T}_{G}$ is graph-isomorphic to 
$\mathcal{T}_{G},$ for all $w$ $\in $ $FP(G^{\symbol{94}}),$ by identifying
their roots. Equivalently, for any vertices $v_{1}$ and $v_{2},$ the
one-flow growing rooted trees $\mathcal{T}_{v_{1}}$ and $\mathcal{T}_{v_{2}}$
having their roots $\omega (v_{1})$ and $\omega (v_{2})$ are
graph-isomorphic. Indeed, both $\mathcal{T}_{v_{1}}$ and $\mathcal{T}%
_{v_{2}} $ are graph-isomorphic to $\mathcal{T}_{G},$ by identifying their
roots, respectively, and hence the trees $\mathcal{T}_{v_{1}}$ and $\mathcal{%
T}_{v_{2}}$ are graph-isomorphic.

\strut \strut \strut

Let $N$ $=$ $\max \{\deg _{out}(v)$ $:$ $v$ $\in $ $V(G)\}.$ Now, choose $%
v_{0}$ $\in $ $V(G)$ having its out-degree $\deg _{out}(v_{0})$ $=$ $N.$
Then we can take the $\mathcal{A}_{G}$-tree $\mathcal{T}_{G}$ as a one-flow
growing rooted tree $\mathcal{T}_{v_{0}},$ having its root $\omega (v_{0})$ $%
=$ $\left( (v_{0},\text{ }v_{0}),\text{ }x_{0}\right) .$ Remember that the
tree $\mathcal{T}_{v_{0}}$ $=$ $\mathcal{T}_{G}$ is a full-subgraph of the $%
2N$-regular tree $\mathcal{T}_{2N}.$ By the connectedness of the shadowed
graph $G^{\symbol{94}}$ of the canonical weighted graph $G,$ for any vertex $%
v$ $\in $ $V(G),$ we can find the finite path $w$ $\in $ $FP(G^{\symbol{94}%
}) $ such that $w$ $=$ $v_{0}$ $w$ $v$ and $w^{-1}$ $=$ $v$ $w$ $v_{0}.$ And
since $\Bbb{G}$ is a fractaloid, the full-subgraph $\mathcal{T}_{v}$ of $%
\mathcal{T}_{v_{0}}$ $=$ $\mathcal{T}_{G}$ are graph-isomorphic. The
existence of $w$ and $w^{-1},$ we can get that

\strut

\begin{center}
$\omega (w)$ $=$ $\left( (v_{0},\text{ }v),\text{ }(t_{i_{1}},\text{ ..., }%
t_{i_{n}})\right) $
\end{center}

and

\begin{center}
$\omega (w^{-1})$ $=$ $\left( (v,\text{ }v_{0}),\text{ }(t_{i_{n}},\text{
..., }t_{i_{1}})\right) ,$
\end{center}

\strut

for some $n$ $=$ $\left| w\right| $ $\in $ $\Bbb{N}.$ This shows that

\strut

\begin{center}
$\deg _{out}^{\mathcal{T}_{G}}\left( \omega (v_{0})\right) $ $=$ $\deg
_{out}^{\mathcal{T}_{G}}\left( (v_{0},\text{ }v_{0}),\text{ }x_{0}\right) $ $%
=$ $2N,$
\end{center}

\strut

where $\deg _{out}^{\mathcal{T}_{G}}(.)$ means the out-degree of vertices of 
$\mathcal{T}_{G}.$ Equivalently, the root $\omega (v_{0})$ of the fixed $%
\mathcal{A}_{G}$-tree $\mathcal{T}_{G}$ $=$ $\mathcal{T}_{v_{0}}$ has $(2N)$%
-incident edges, which are the first level of the tree. This means that if $%
\deg _{out}(v_{0})$ $=$ $N,$ then the root $\omega (v_{0})$ of $\mathcal{T}%
_{G}$ has $\deg _{out}^{\mathcal{T}_{G}}(v_{0})$ $=$ $2N.$

\strut

Let $v_{0}$ be given as before and let $\psi _{\pm 1},$ ..., $\psi _{\pm N}$
be the $2N$-incident edges of the root $\omega (v_{0})$ of $\mathcal{T}%
_{v_{0}}$ $=$ $\mathcal{T}_{G}.$ For convenience, let $(\widetilde{pr}$ $%
\circ $ $\omega )(\psi _{k})$ $=$ $t_{k},$ where $t_{k}$ $=$ $x_{k}$ if $k$ $%
\in $ $\{1,$ ..., $N\},$ and $t_{k}$ $=$ $-x_{k},$ if $k$ $\in $ $\{-1,$
..., $-N\}.$ Then, for the terminal vertices

\strut

\begin{center}
$\varphi _{k}$ $=$ $\varphi \left( \left( (v_{0},\text{ }v_{k}),\text{ }%
t_{k}\right) ,\text{ }\psi _{k}\right) ,$ for $k$ $=$ $\pm 1,$ ..., $\pm N,$
\end{center}

\strut

\strut create the full-subgraphs $\mathcal{T}_{v_{k}}$'s of $\mathcal{T}%
_{v_{0}},$ having their roots $\varphi _{k},$ for all $k$ $=$ $\pm 1,$ ..., $%
\pm N,$ and they are graph-isomorphic to $\mathcal{T}_{v_{0}},$ since $\Bbb{G%
}$ is a fractaloid. Note that the full-subgraphs $\mathcal{T}_{v_{k}}$'s are
determined by the rooted one-flow growing tree induced by the $1$-st level
of $\mathcal{T}_{v_{0}}.$ Thus, $\mathcal{T}_{v_{k}}$'s have their $1$-st
level consisting of $2N$-elements, for all $k$ $=$ $\pm 1,$ ..., $\pm N.$
Inductively, we can conclude that, for any $w$ $\in $ $FP(G^{\symbol{94}}),$
the $w$-part $\mathcal{T}_{w}$ should have its $1$-st level consisting of $%
2N $-elements. Equivalently, the tree $\mathcal{T}_{v_{0}}$ should be
graph-isomorphic to the $2N$-regular tree $\mathcal{T}_{2N}.$ Since $%
\mathcal{T}_{v_{0}}$ is our $\mathcal{A}_{G}$-tree $\mathcal{T}_{G},$ we can
get that the $\mathcal{A}_{G}$-tree $\mathcal{T}_{G}$ is graph-isomorphic to
the $2N$-regular tree $\mathcal{T}_{2N},$ by identifying their roots.
\end{proof}

\strut

The above theorem shows that the converse of the previous corollary also
holds true. Therefore, we can get the following other characterization of
fractaloids.

\strut

\begin{corollary}
Let $G$ be a canonical weighted graph with its labeled graph groupoid $\Bbb{G%
},$ and let $\mathcal{A}_{G}$ be the graph automaton induced by $G.$ Then
the $\mathcal{A}_{G}$-actions $\{\mathcal{A}_{w}$ $:$ $w$ $\in $ $FP(G^{%
\symbol{94}})\}$ act fully on the $2N$-regular tree $\mathcal{T}_{2N},$ if
and only if $\Bbb{G}$ is a fractaloid.
\end{corollary}

\strut

\begin{proof}
We know that the labeled graph groupoid $\Bbb{G}$ is a fractaloid if the $%
\mathcal{A}_{G}$-actions act fully on the $2N$-regular tree $\mathcal{T}%
_{2N},$ by the previous corollary. Now, by the previous theorem, we can
conclude that, if $\Bbb{G}$ is a fractaloid, then the $\mathcal{A}_{G}$-tree 
$\mathcal{T}_{G}$ is identified with the $2N$-regular tree $\mathcal{T}%
_{2N}. $ Therefore, if $\Bbb{G}$ is a fractaloid, then the $\mathcal{A}_{G}$%
-actions act fully on $\mathcal{T}_{2N},$ since they act fully on $\mathcal{T%
}_{G}.$
\end{proof}

\strut \strut

\strut

\subsection{Examples}

\strut

\strut

In this section, we consider some examples of fractaloids.

\strut

\begin{example}
Let $G$ be the one-flow circulant graph with three vertices. i.e., it is a
graph with

\strut

\begin{center}
$V(G)$ $=$ $\{v_{1},$ $v_{2},$ $v_{3}\}$
\end{center}

and

\begin{center}
$E(G)$ $=$ $\{e_{j}$ $=$ $v_{j}$ $e_{j}$ $v_{j+1}$ $:$ $j$ $=$ $1,$ $2,$ $3,$
$v_{4}$ $\overset{def}{=}$ $v_{1}\}.$
\end{center}

\strut

Then we have $\max \{\deg _{out}(v_{j})$ $:$ $j$ $=$ $1,$ $2,$ $3\}$ $=$ $1.$
So, as in Section 3.1, we can decide the labeling set $\{0,$ $1\}$, and the
labeled graph groupoid $\Bbb{G}$ of $G,$ labeled by $\{0,$ $\pm 1\}^{*}.$
Then the graph automaton $\mathcal{A}_{G}$ of $G$ is determined and the
automata actions of $\mathcal{A}_{G}$ act fully on the $2$-regular tree $%
\mathcal{T}_{2},$ by the existence of the pairs $(\mathcal{A}_{e_{1}},$ $%
\mathcal{A}_{e_{3}^{-1}}),$ $(\mathcal{A}_{e_{2}},$ $\mathcal{A}%
_{e_{1}^{-1}})$ and $(\mathcal{A}_{e_{3}},$ $\mathcal{A}_{e_{2}^{-1}})$ of
actions. Therefore, the graph groupoid $\Bbb{G}$ of $G$ is a fractaloid.
\end{example}

\strut

By the inductive modification of the previous example, we can get the
following proposition.

\strut

\begin{proposition}
Let $G$ be a one-flow circulant graph with $N$-vertices with $V(G)$ $=$ $%
\{v_{1},$ ..., $v_{N}\}$ and $E(G)$ $=$ $\{e_{j}$ $=$ $v_{j}$ $e_{j}$ $%
v_{j+1}$ $:$ $j$ $=$ $1,$ ..., $N,$ $v_{N+1}$ $\overset{def}{=}$ $v_{1}\}.$
Then the graph groupoid $\Bbb{G}$ of $G$ is a fractaloid. $\square $
\end{proposition}

\strut

The following example is about certain fractal groups.

\strut

\begin{example}
Let $G$ be the one-vertex-$N$-loop-edge graph with $V(G)$ $=$ $\{v\}$ and $%
E(G)$ $=$ $\{e_{j}$ $=$ $v$ $e_{j}$ $v$ $:$ $j$ $=$ $1,$ ..., $N\}.$ In [10]
and [11], we showed that the graph groupoid $\Bbb{G}$ of $G$ is a group
which is group-isomorphic to the free group $F_{N}$ with $N$-generators. The
free group $F_{N}$ is a fractal group, for $N$ $\in $ $\Bbb{N}$ (See [1]).
And hence, the graph groupoid $\Bbb{G}$ is a fractaloid which is a fractal
group.
\end{example}

\strut \strut

By the previous example, we can get the following proposition.

\strut

\begin{proposition}
Let $G$ be a one-vertex-multi-loop-edge graph. Then the graph groupoid $\Bbb{%
G}$ of $G$ is a fractaloid. $\square $
\end{proposition}

\strut \strut

Recall that, we say a directed graph $G$ is \emph{finite}, if $\left|
V(G)\right| $ $<$ $\infty ,$ and $\left| E(G)\right| $ $<$ $\infty .$ The
previous two examples show that it is possible that, even though the given
graph is finite, its graph groupoid can be a fractaloid.

\strut

The following example is trivial. It shows that the automata actions induced
by the $N$-regular tree $\mathcal{T}_{N}$ act fully on the $2N$-regular tree 
$\mathcal{T}_{2N},$ and hence they generate a fractaloid.

\strut

\begin{example}
Let $\mathcal{T}_{N}$ be the $N$-regular tree. Then it induces the graph
automaton $\mathcal{A}_{T_{N}}$, and the automata actions act fully on the $%
2N$-regular tree $\mathcal{T}_{2N}.$ Therefore, the graph groupoid $\Bbb{T}%
_{N}$ of $\mathcal{T}_{N}$ is a fractaloid.
\end{example}

\strut \strut

\begin{example}
Every infinitely countable linear graph $L$ induces a fractaloid $\Bbb{L}.$
i.e., the graph $L$ is graph-isomorphic to the following graph,

\strut

\begin{center}
$\bullet \longrightarrow \bullet \longrightarrow \bullet \longrightarrow
\bullet \longrightarrow \cdot \cdot \cdot .$
\end{center}

\strut

Then the graph $L$ induces the graph automaton $\mathcal{A}_{L}.$ Since the
(non-reduced) finite path set $FP(L^{\symbol{94}})$ of the shadowed graph $%
L^{\symbol{94}}$ is bijective to the finite path set $FP(\mathcal{T}_{2})$
of $2$-regular tree $\mathcal{T}_{2}.$ This guarantees that automata actions
of $\mathcal{A}_{L}$ act fully on $\mathcal{T}_{2}.$ Therefore, $\Bbb{L}$ is
a fractaloid. Indeed, the graph $L$ is regarded as the $1$-regular tree, and
hence, by the previous example, it induces a fractaloid.
\end{example}

\strut \strut \strut \strut

\strut

\strut \strut \strut

\section{Labeling Operators of Fractaloids\strut}

\strut

\strut

Let $G$ be a canonical weighted graph with its labeled graph groupoid $\Bbb{G%
},$ and let

\strut

\begin{center}
$N$ $=$ $\max \{\deg _{out}(v)$ $:$ $v$ $\in $ $V(G)\}.$
\end{center}

\strut

We will use the same notations used in the previous Sections. In Section
2.4, we defined right graph von Neumann algebras $\Bbb{M}_{G}$ $=$ $M$ $%
\times _{\beta }$ $\Bbb{G},$ where $M$ is an arbitrary fixed von Neumann
algebra, contained in the operator algebra $B(K),$ and $\beta $ is a right $%
G $-representation, which is an intertwined groupoid action of $\Bbb{G},$
acting on the Hilbert space $K$ $\otimes $ $H_{G}.$ As we mentioned in
Section 2.4, we are only interested in the case where $M$ $=$ $\Bbb{C}.$
Then the right graph von Neumann algebras $\Bbb{M}_{G}$ $=$ $\Bbb{C}$ $%
\times _{\beta }$ $\Bbb{G},$ for all $\beta ,$ are $*$-isomorphic to the
groupoid von Neumann algebra

$\strut $

\begin{center}
$\overline{\Bbb{C}[\Bbb{G}]}^{w}$ $=$ $\overline{\Bbb{C}[\beta (\Bbb{G})]}%
^{w},$ for all $\beta ,$
\end{center}

\strut

as a $W^{*}$-subalgebra of $B(H_{G}),$ where $H_{G}$ is the graph Hilbert
space induced by the graph $G.$

\strut

\textbf{Notation} From now, if we denote $\Bbb{M}_{G},$ then it is
automatically assumed to be the $W^{*}$-subalgebra $\overline{\Bbb{C}[\Bbb{G}%
]}^{w}$ of $B(H_{G}).$ And we call it ``the'' right graph von Neumann
algebra induced by $G.$ $\square $

\strut

\strut

\strut

\subsection{Labeling Operators}

\strut

\strut

Now, we will define certain operators on the graph Hilbert space $H_{G}.$
Then we can easily realize that these operators are elements of $\Bbb{M}%
_{G}. $

\strut

\begin{definition}
Let $G$ be a canonical weighted graph and let $H_{G}$ be the corresponding
graph Hilbert space. Define an operator $T_{i}$ on $H_{G}$ by

\strut

\begin{center}
$T_{i}(\xi _{w})$ $\overset{def}{=}$ $\xi _{w}$ $\xi _{e}$ $=$ $\xi _{we}$ $%
= $ $\left\{ 
\begin{array}{ll}
\begin{array}{l}
\xi _{we}
\end{array}
& 
\begin{array}{l}
\begin{array}{l}
\text{if }\exists e\in E(G^{\symbol{94}})\text{ s.t.} \\ 
\,\,\,\,\,\,\,\,\,w\text{ }e\neq \emptyset \text{ and } \\ 
\omega (e)=\left( (v,\text{ }v^{\prime }),\text{ }t_{i}\right)
\end{array}
\end{array}
\\ 
\quad 0 & \quad \text{otherwise,}
\end{array}
\right. $
\end{center}

\strut

for all $\xi _{w}$ $\in $ $\mathcal{B}_{H_{G}}$, with $w$ $=$ $w$ $v$ in $%
FP_{r}(G^{\symbol{94}}),$ and for all $i$ $=$ $\pm 1,$ ..., $\pm N,$ where $%
\mathcal{B}_{H_{G}}$ is the Hilbert basis of $H_{G}.$ Here, $t_{i}$ $=$ $%
x_{i}$ if $t_{i}$ $\in $ $X,$ and $t_{i}$ $=$ $-x_{i},$ if $t_{i}$ $\in $ $%
-X.$ The operators $T_{k}$'s are called the $k$-th labeling operators, for $%
k $ $=$ $\pm 1,$ ..., $\pm N.$ The operator $T_{G}$ on $H_{G}$ is said to be
the labeling operator, if

\strut

\begin{center}
$T_{G}$ $\overset{def}{=}$ $\sum_{j=-N}^{-1}$ $T_{j}$ $+$ $\sum_{i=1}^{N}$ $%
T_{i}.$
\end{center}
\end{definition}

\strut

The labeling operator is similar to the Hecke-type operators (e.g., [35]) or
the Ruelle operators (e.g., [31] and [32]) or the radial operators (e.g.,
[36]) The definition shows that each labeling $t_{k}$ $\in $ $\pm X$
generates the $k$-th labeling operator $T_{k}$ on $H_{G},$ for $k$ $=$ $\pm
1,$ ..., $\pm N.$ And the labeling set $\pm X$ $=$ $\{\pm x_{1},$ ..., $\pm
x_{N}\},$ itself, generates the labeling operator $T_{G}$ on $H_{G}.$

\strut

By the very definition, we can realize that the $k$-th labeling operators $%
T_{k}$'s, for $k$ $=$ $\pm 1,$ ..., $\pm N,$ and the labeling operator $%
T_{G} $ of $\Bbb{G}$ are contained in the right graph von Neumann algebra $%
\Bbb{M}_{G}.$

\strut

\begin{lemma}
Let $G$ be a canonical weighted graph with its labeled graph groupoid $\Bbb{G%
},$ and let $\Bbb{M}_{G}$ be the right graph von Neumann algebra induced by $%
G.$ Also, let $T_{k}$'s and $T_{G}$ be the $k$-th labeling operators and the
labeling operator of $\Bbb{G}$ on the graph Hilbert space $H_{G},$ where $k$ 
$=$ $\pm 1,$ ..., $\pm N.$ Then they are contained in $\Bbb{M}_{G}.$
\end{lemma}

\strut

\begin{proof}
Let $\Bbb{M}_{G}$ be the right graph von Neumann algebra induced by $G.$
Construct the elements $\tau _{k}$ and $\tau $ of $\Bbb{M}_{G}$:

\strut

\begin{center}
$\tau _{k}$ $=$ $\underset{e\in E(G^{\symbol{94}}),\text{ }(\widetilde{pr}%
\text{ }\circ \text{ }\omega )(e)=t_{k}}{\sum_{n}^{(N)}}$ $R_{e},$
\end{center}

where

\begin{center}
$t_{k}$ $\overset{def}{=}$ $\left\{ 
\begin{array}{lll}
x_{k} &  & \text{if }k=1,\text{ }...,\text{ }N \\ 
-x_{k} &  & \text{if }k=-1,\text{ }...,\text{ }-N,
\end{array}
\right. $
\end{center}

\strut

for all $k$ $=$ $\pm 1,$ ..., $\pm N.$ Then these elements $\tau _{k}$'s are
identified with the $k$-th labeling operators $T_{k},$ on the graph Hilbert
space $H_{G}$, for all $k$ $=$ $\pm 1,$ ..., $\pm N.$ Therefore, the $k$-th
labeling operators $T_{k}$'s are contained in $\Bbb{M}_{G},$ for all $k.$

\strut

Define now the element $\tau $ of $\Bbb{M}_{G}$ by

\strut

\begin{center}
$\tau $ $=$ $\sum_{i=-N}^{-1}$ $\tau _{i}$ $+$ $\sum_{j=1}^{N}$ $\tau _{j},$
\end{center}

\strut

where $\tau _{k}$'s are defined as above, for all $k$ $=$ $\pm 1,$ ..., $\pm
N.$ Then, by definition, $\tau $ is identified with $T_{G}$ on $H_{G}.$
\end{proof}

\strut

By the previous lemma, we can regard the labeling operators $T_{k}$'s and $%
T_{G}$ as elements of the right graph von Neumann algebra $\Bbb{M}_{G}.$

\strut

\begin{corollary}
Let $T_{G}$ be the labeling operator of $\Bbb{G}$ on $H_{G}.$ Then it is
identified with the operator $\underset{e\in E(G^{\symbol{94}})}{\sum }$ $%
R_{e}$ generated by all edges of the shadowed graph $G^{\symbol{94}}.$ $%
\square $
\end{corollary}

\strut

By the previous lemma, we can get the following adjoint property of the $k$%
-th labeling operators $T_{k}$'s.

\strut

\begin{lemma}
Let $k$ $\in $ $\{\pm 1,$ ..., $\pm N\}$ be fixed, and let $T_{k}$ be the $k$%
-th labeling operator on $H_{G}.$ Then the adjoint $T_{k}^{*}$ of $T_{k}$ is
identical to the $(-k)$-th labeling operator $T_{-k}.$
\end{lemma}

\strut

\begin{proof}
By the previous lemma, we can consider the $k$-th labeling operator $T_{k}$
as an element in $\Bbb{M}_{G}.$ i.e., let

\strut

\begin{center}
$T_{k}$ $=$ $\underset{e\in E(G^{\symbol{94}}),\,(\widetilde{pr}\text{ }%
\circ \text{ }\omega )(e)=t_{k}}{\sum }$ $R_{e}$ in $\Bbb{M}_{G}.$
\end{center}

\strut

Then we can have that

\strut

\begin{center}
$
\begin{array}{ll}
T_{k}^{*} & =\underset{e\in E(G^{\symbol{94}}),\text{ }(\widetilde{pr}\text{ 
}\circ \text{ }\omega )(e)=t_{k}}{\sum }R_{e}^{*}=\underset{e\in E(G^{%
\symbol{94}}),\text{ }(\widetilde{pr}\text{ }\circ \text{ }\omega )(e)=t_{k}%
}{\sum }R_{e^{-1}} \\ 
&  \\ 
& =\underset{e^{-1}\in E(G^{\symbol{94}}),\text{ }(\widetilde{pr}\text{ }%
\circ \text{ }\omega )(e^{-1})=-t_{k}}{\sum }R_{e^{-1}}=\underset{e\in E(G^{%
\symbol{94}}),\text{ }(\widetilde{pr}\text{ }\circ \text{ }\omega )(e)=-t_{k}%
}{\sum }R_{e}=T_{-k}.
\end{array}
$
\end{center}

\strut \strut

Indeed, by the canonical weighting process, for any $e$ $\in $ $E(G^{\symbol{%
94}}),$ $\widetilde{pr}\left( \omega (e)\right) $ $=$ $t_{k}$ if and only if 
$\widetilde{pr}\left( \omega (e^{-1})\right) $ $=$ $-t_{k},$ where $t_{k}$ $%
\in $ $\pm X.$
\end{proof}

\strut \strut

Recall that, by Section 2.5, [10] and [11], we know that two elements $w_{1}$
and $w_{2}$ of a graph groupoid $\Bbb{G}$ are diagram-distinct, in the sense
that (i) $w_{1}$ $\neq $ $w_{2}^{-1},$ and (ii) the diagrams (which are the
graphical images) of $w_{1}$ and $w_{2}$ (in $\Bbb{R}^{2}$) are distinct, if
and only if the right multiplication operators $R_{w_{1}}$ and $R_{w_{2}}$
are free over $\Bbb{D}_{G}$ in $(\Bbb{M}_{G},$ $E),$ where $\Bbb{D}_{G}$ is
the $\Bbb{C}$-diagonal subalgebra of $\Bbb{M}_{G}$ and $E$ is the canonical
conditional expectation of $\Bbb{M}_{G}$ onto $\Bbb{D}_{G}.$ Thus, more
precisely, two subsets $X_{1}$ and $X_{2}$ of $\Bbb{G}$ are
diagram-distinct, in the sense that, for any pair $(w_{1},$ $w_{2})$ $\in $ $%
X_{1}$ $\times $ $X_{2},$ the elements $w_{1}$ and $w_{2}$ are
diagram-distinct if and only if the operators $\underset{w_{1}\in X_{1}}{%
\sum }$ $R_{w_{1}}$ and $\underset{w_{2}\in X_{2}}{\sum }$ $R_{w_{2}}$ are
free over $\Bbb{D}_{G}$ in $(\Bbb{M}_{G},$ $E)$. So, we can get the
following proposition.

\strut

\begin{proposition}
The operators $T_{k}$ $+$ $T_{-k}$ $\in $ $\Bbb{M}_{G},$ for $k$ $=$ $1,$
..., $N,$ are free over $\Bbb{D}_{G}$ from each other in $(\Bbb{M}_{G},$ $%
E), $ where $T_{k}$'s are the $k$-th labeling operators, for $k$ $=$ $\pm 1,$
..., $\pm N.$
\end{proposition}

\strut

\begin{proof}
Let's denote the operators $T_{k}$ $+$ $T_{-k}$ by $S_{k},$ for $k$ $=$ $1,$
..., $N.$ And, for the fixed $k$, construct the set $\mathcal{S}_{k}$ by

\strut

\begin{center}
$\mathcal{S}_{k}$ $=$ $\{e$ $\in $ $E(G^{\symbol{94}})$ $:$ $(\widetilde{pr}$
$\circ $ $\omega )(e)$ $=$ $x_{k}$ or $-x_{k}\},$
\end{center}

\strut

for $k$ $=$ $1,$ ..., $N.$ Then the operator $S_{k}$ can be re-written by

\strut

\begin{center}
$S_{k}$ $=$ $\underset{e\in \mathcal{S}_{k}}{\sum }$ $R_{e},$ for all $k$ $=$
$1,$ ..., $N.$
\end{center}

\strut

By the canonical weighting process, the sets $\mathcal{S}_{k}$'s are
diagram-distinct from each other, for all $k$ $=$ $1,$ ..., $N.$ Therefore,
the operators $S_{k}$'s are free over $\Bbb{D}_{G}$ from each other in $(%
\Bbb{M}_{G},$ $E).$
\end{proof}

\strut

\begin{corollary}
(1) The labeling operator $T_{G}$ of $\Bbb{G}$ on $H_{G}$ is a $\Bbb{D}_{G}$%
-free sum of the $\Bbb{D}_{G}$-free elements $T_{k}$ $+$ $T_{-k},$ in $\Bbb{M%
}_{G},$ for all $k$ $=$ $1,$ ..., $N.$

\strut

(2) The labeling operator $T_{G}$ is self-adjoint on $H_{G}.$ $\square $
\end{corollary}

\strut

\strut

\strut

\subsection{Amalgamated Free Distributional Data of Labeling Operators}

\strut

\strut

\strut

Let $B$ $\subset $ $A$ be von Neumann algebras with $1_{B}$ $=$ $1_{A},$ and
assume that there exists a conditional expectation $E_{B}$ $:$ $A$ $%
\rightarrow $ $B.$ Then the pair $(A,$ $E_{B})$ is a $B$-valued $W^{*}$%
-probability space. Let $a$ $\in $ $(A,$ $E_{B})$ be a $B$-valued random
variable. Then the $B$-valued free distributional data of $a$ is
characterized by the $B$-valued joint $*$-moments

$\strut $

\begin{center}
$\cup _{n=1}^{\infty }\left\{
E_{B}(b_{1}a^{r_{1}}b_{2}a^{r_{2}}...b_{n}a^{r_{n}})\left| 
\begin{array}{c}
(r_{1},...,r_{n})\in \{1,*\}^{n} \\ 
b_{1},\text{ ..., }b_{n}\in B\text{ arbitrary}
\end{array}
\right. \right\} .$
\end{center}

\strut

i.e., the above $B$-values generate the $B$-valued free distribution $\rho
_{a}$ of $a$ (See [5], [12] and [21]). Suppose now that the operator $a$ is
self-adjoint. Then its $B$-valued free distribution $\rho _{a}$ is
completely determined by the $B$-valued moments of $a,$

\strut

\begin{center}
$\left\{ E_{B}(b_{1}ab_{2}a...b_{n}a):\left| 
\begin{array}{c}
n\in \Bbb{N},\text{ and} \\ 
b_{1},\text{ ..., }b_{n}\in B\text{ arbitrary}
\end{array}
\right. \right\} .$
\end{center}

\strut

Since the labeling operator $T_{G}$ of $\Bbb{G}$ on the graph Hilbert space $%
H_{G}$ is an element in the right graph von Neumann algebra $\Bbb{M}_{G},$
we can consider its $\Bbb{D}_{G}$-valued free distributional data, by
observing the $\Bbb{D}_{G}$-valued moments $\{E(T_{G}^{n})$ $:$ $n$ $\in $ $%
\Bbb{N}\}.$ Since

\strut

\begin{center}
$T_{G}$ $=$ $\left( \sum_{i=-N}^{-1}T_{i}\right) $ $+$ $\left(
\sum_{j=1}^{N}T_{j}\right) ,$
\end{center}

\strut

where $T_{k}$'s are the $k$-th labeling operators, for $k$ $=$ $\pm 1,$ ..., 
$\pm N,$ we can have that

\strut

$\quad E\left( T_{G}^{n}\right) $ $=$ $\underset{(i_{1},\text{ ..., }%
i_{n})\in \{\pm 1,\text{ ..., }\pm N\}^{n}}{\sum }$ $E\left(
T_{i_{1}}...T_{i_{n}}\right) $

\strut (5.1)

$\qquad =$ $\underset{(i_{1},\text{ ..., }i_{n})\in \{\pm 1,\text{ ..., }\pm
N\}^{n}}{\sum }\left( \underset{(e_{1},...,e_{n})\in E(G^{\symbol{94}})^{n},%
\text{ }(\widetilde{pr}\text{ }\circ \text{ }\omega )(e_{j})=t_{i_{j}}}{\sum 
}\left( E(R_{e_{1}}...R_{e_{n}})\right) \right) ,$

\strut

for all $n$ $\in $ $\Bbb{N}.$

\strut

\strut

\strut

\subsection{Labeling Operators of Fractaloids}

\strut

\strut

In this section, we will consider the $\Bbb{D}_{G}$-valued moments of the
labeling operator $T_{G}$ of $\Bbb{G}$ in $\Bbb{M}_{G},$ where $\Bbb{G}$ is
a fractaloid. We showed that, in general, if $T_{G}$ is the labeling
operator of an arbitrary labeled graph groupoid $\Bbb{G}$ in the right graph
von Neumann algebra $\Bbb{M}_{G},$ then the $n$-th $\Bbb{D}_{G}$-valued
moment $E(T_{G}^{n})$ is determined by (5.1).

\strut

Now, we are interested in the more precise formula than (5.1), when we have
a fractaloid $\Bbb{G}$. If we can find the more precise formula, then we not
only explain the $\Bbb{D}_{G}$-valued free distributional data of $T_{G},$
but also show how the fractal property works on graph groupoids (and hence
on graph Hilbert spaces). Moreover, it also shows how the admissibility of
fractaloids works inside the right graph von Neumann algebra $\Bbb{M}_{G}.$
Therefore, the study of the $\Bbb{D}_{G}$-valued free distributional data of
the labeling operator $T_{G}$ of a fractaloid $\Bbb{G}$ is for studying the
groupoidal-and-operator-algebraic fractal properties.

\strut

From now, all graph groupoids in this Section are fractaloids.

\strut

\textbf{Notation} For convenience, we will denote the map $\widetilde{pr}$ $%
\circ $ $\omega $ simply by $\widetilde{\omega }.$ i.e.,

\strut

\begin{center}
$\widetilde{\omega }$ $:$ $\Bbb{F}^{+}(G^{\symbol{94}})$ $\rightarrow $ $%
(\pm X_{0}^{*})$
\end{center}

defined by

\begin{center}
$\widetilde{\omega }(w)$ $=$ $\left\{ 
\begin{array}{ll}
\begin{array}{l}
x_{0} \\ 
\\ 
\widetilde{pr}\left( \omega (w)\right)
\end{array}
& 
\begin{array}{l}
\text{if }w\in V(G) \\ 
\\ 
\text{if }w\in E(G^{\symbol{94}})
\end{array}
\\ 
\begin{array}{l}
\\ 
\left( \widetilde{\omega }(e_{1}),\text{ ..., }\widetilde{\omega }%
(e_{k})\right)
\end{array}
& 
\begin{array}{l}
\\ 
\text{if }w=e_{1}...e_{k}\in FP(G^{\symbol{94}}),\text{ } \\ 
\qquad \qquad \text{for }k\geq 2
\end{array}
\\ 
\begin{array}{l}
\\ 
\emptyset _{*}
\end{array}
& 
\begin{array}{l}
\\ 
\text{if }w=\emptyset ,
\end{array}
\end{array}
\right. $
\end{center}

\strut

for all $w$ $\in $ $\Bbb{F}^{+}(G^{\symbol{94}}).$ Recall the operations $%
\theta $ and $\omega _{+}$ defined in Section 3.2. $\square $

\strut

By Section 4.1, we have the following characterization of fractaloids: the
labeled graph groupoid $\Bbb{G}$ is a fractaloid if and only if the automata
actions $\{\mathcal{A}_{w}$ $:$ $w$ $\in $ $FP(G^{\symbol{94}})\}$ of the
graph automaton $\mathcal{A}_{G}$ act fully on the $2N$-regular tree $%
\mathcal{T}_{2N}.$ So, if $T_{G}$ is the labeling operator of a fractaloid $%
\Bbb{G}$ in the right graph von Neumann algebra $\Bbb{M}_{G},$ then the
formula (5.1) is re-formulated as follows:

\strut

$\quad E(T_{G}^{n})$ $=$ $\underset{(i_{1},...,i_{n})\in \{\pm 1,\text{ ..., 
}\pm N\}^{n}}{\sum }$ $\left( \underset{e_{1}...e_{n}\in FP(G^{\symbol{94}%
}),\,\widetilde{\omega }(e_{j})=t_{i_{j}}}{\sum }\text{ }%
E(R_{e_{1}...e_{n}})\right) $

\strut

by (5.1)

\strut

$\qquad =$ $\underset{(i_{1},...,i_{n})\in \{\pm 1,\text{ ..., }\pm N\}^{n}}{%
\sum }$ $\left( \underset{e_{1}...e_{n}\in FP(G^{\symbol{94}%
}),\,e_{1}...e_{n}\in V(G)\subset \Bbb{G},\text{ }\widetilde{\omega }%
(e_{j})=t_{i_{j}}}{\sum }\text{ }E(R_{e_{1}...e_{n}})\right) $

\strut

$\qquad =$ $\underset{(i_{1},...,i_{n})\in \{\pm 1,\text{ ..., }\pm N\}^{n}}{%
\sum }$ $\left( \underset{v\in V(G)\subset \Bbb{G},\text{ }%
v=\,e_{1}...e_{n}\in FP(G^{\symbol{94}}),\,\widetilde{\omega }%
(e_{j})=t_{i_{j}}}{\sum }\text{ }R_{v}\right) $

\strut

since $E(R_{v})$ $=$ $R_{v},$ for all $v$ $\in $ $V(G)$ $\subset $ $\Bbb{G}$

\strut

\strut (5.2)

$\qquad =$ $\underset{(i_{1},...,i_{n})\in \{\pm 1,\text{ ..., }\pm N\}^{n}}{%
\sum }$ $\left( \underset{e_{1}...e_{n}\in FP(G^{\symbol{94}}),\,\widetilde{%
\omega }(e_{j})=t_{i_{j}},\,\sum_{k=1}^{n}t_{i_{k}}=x_{0}\,}{\sum }\text{ }%
R_{e_{1}...e_{n}}\right) ,$

\strut

by Section 3.2: Remember that a finite path $e_{1}$ ... $e_{n}$ $=$ $v$ $%
e_{1}$ ... $e_{n}$ $v$ in $FP(G^{\symbol{94}}),$ with $v$ $\in $ $V(G^{%
\symbol{94}}),$ is the vertex $v$ in the labeled graph groupoid $\Bbb{G}$
(under the reduction), where $e_{j}$ $\in $ $E(G^{\symbol{94}})$ with $%
\widetilde{\omega }(e_{j})$ $=$ $t_{i_{j}},$ if and only if $\sum_{j=1}^{n}$ 
$t_{i_{j}}$ $=$ $x_{0}$ in $\Bbb{C}^{\oplus \,N},$ in general.

\strut

Fix now a vertex $v$ $\in $ $V(G).$ Define a subset $\mathcal{F}_{v,\,n}$ of
the finite path set $FP(G^{\symbol{94}})$ of the shadowed graph $G^{\symbol{%
94}}$ by

\strut

\begin{center}
$\mathcal{F}_{v,\,n}$ $\overset{def}{=}$ $\left\{ e_{1}\text{ ... }e_{n}\in
FP(G^{\symbol{94}})\left| 
\begin{array}{c}
e_{1}=v\,e_{1}\text{ and} \\ 
\sum_{j=1}^{n}\widetilde{\omega }(e_{j})=x_{0}
\end{array}
\right. \right\} ,$
\end{center}

\strut

for all $n$ $\in $ $\Bbb{N}.$ i.e., the sets $\mathcal{F}_{v,\,n}$ are the
collection of all length-$n$ finite paths of $FP(G^{\symbol{94}})$ which are
all identified with the vertex $v$ in the graph groupoid $\Bbb{G},$ under
the reduction (RR), for all $n$ $\in $ $\Bbb{N}.$ We will say that the set $%
\mathcal{F}_{v,\,n}$ is \emph{the length-}$n$\emph{-}$v$\emph{-subset of} $%
FP(G^{\symbol{94}}),$ for all $n$ $\in $ $\Bbb{N}.$ Remark that the
weighting (or the labeling) process makes us understand the admissibility
and the reduction on $\mathcal{F}_{v,\text{ }n}$'s (inside $\Bbb{G}$), more
in detail, for all $v$ $\in $ $V(G)$ and $n$ $\in $ $\Bbb{N}$. Again, by the
labeling process, we can re-express the sets $\mathcal{F}_{v,\,n}$'s as
follows,

\strut

\begin{center}
$\mathcal{F}_{v,\,n}$ $=$ $\underset{(i_{1},...,i_{n})\in \{\pm 1,\text{
..., }\pm N\}^{n}}{\cup }\left\{ e_{1}\text{ ... }e_{n}\in FP(G^{\symbol{94}%
})\left| 
\begin{array}{c}
e_{1}=v\text{ }e_{1}\in E(G^{\symbol{94}}) \\ 
\widetilde{\omega }(e_{j})=t_{i_{j}},\text{ for }j, \\ 
\sum_{j=1}^{n}t_{i_{j}}=x_{0}
\end{array}
\right. \right\} ,$
\end{center}

\strut

for all $n$ $\in $ $\Bbb{N}.$ Let $v_{1}$ and $v_{2}$ be arbitrarily chosen
vertices in $V(G).$ Then we can have the corresponding sets $\{\mathcal{F}%
_{v_{1},\,n}\}_{n=1}^{\infty }$ and $\{\mathcal{F}_{v_{2},\,n}\}_{n=1}^{%
\infty }.$ Then, for each $n$ $\in $ $\Bbb{N},$ the cardinalities $\left| 
\mathcal{F}_{v_{1},\,n}\right| $ and $\left| \mathcal{F}_{v_{2},\,n}\right| $
are identical, whenever $\Bbb{G}$ is a fractaloid.

\strut

\begin{lemma}
Let $\Bbb{G}$ be a given fractaloid, and let $\mathcal{F}_{v,\,n}$ be the
length-$n$-$v$-subsets of $FP(G^{\symbol{94}}),$ for all $v$ $\in $ $V(G)$
and $n$ $\in $ $\Bbb{N}.$ Then, for any $n$ $\in $ $\Bbb{N},$ the
cardinalities $\left| \mathcal{F}_{v,\,n}\right| $ are identical, for all $v$
$\in $ $V(G).$
\end{lemma}

\begin{proof}
Let $v_{1}$ $\neq $ $v_{2}$ be the distinct vertices in $V(G),$ and let $%
\mathcal{F}_{v_{1},\,n}$ and $\mathcal{F}_{v_{2},\,n}$ be the length-$n$%
-vertex-subsets of $v_{1}$ and $v_{2}$ in $FP(G^{\symbol{94}}),$
respectively, for $n$ $\in $ $\Bbb{N}.$ By definition,

\strut

\begin{center}
$\mathcal{F}_{v_{k},\,n}$ $=$ $\underset{(i_{1},...,i_{n})\in \{\pm 1,\text{
..., }\pm N\}^{n}}{\cup }\left\{ e_{1}\text{ ... }e_{n}\in FP(G^{\symbol{94}%
})\left| 
\begin{array}{c}
e_{1}=v_{k}\text{ }e_{1}\in E(G^{\symbol{94}}) \\ 
\widetilde{\omega }(e_{j})=t_{i_{j}},\text{ for }j, \\ 
\sum_{j=1}^{n}t_{i_{j}}=x_{0}
\end{array}
\right. \right\} ,$
\end{center}

\strut

for $k$ $=$ $1,$ $2.$ Define now the map

\strut

\begin{center}
$\mathcal{E}$ $:$ $E(G^{\symbol{94}})$ $\rightarrow $ $\{\mathcal{A}_{e}$ $:$
$e$ $\in $ $E(G^{\symbol{94}})\}$
\end{center}

by

\begin{center}
$\mathcal{E}(e)$ $=$ $\mathcal{A}_{e},$ for all $e$ $\in $ $E(G^{\symbol{94}%
}).$
\end{center}

\strut

Clearly, this map $\mathcal{E}$ is bijective and hence it is extendable to
the bijective map, also denoted by $\mathcal{E}$, from $FP(G^{\symbol{94}})$
onto $\{\mathcal{A}_{w}$ $:$ $w$ $\in $ $FP(G^{\symbol{94}})\}.$ The
existence of this bijection $\mathcal{E}$ shows that we can regard the
elements in $\mathcal{F}_{v_{k},\,n}$ as the $\mathcal{E}$-corresponding
elements on $FP(\mathcal{T}_{2N}),$ where $\mathcal{T}_{2N}$ is the $2N$%
-regular tree where the automata actions $\mathcal{A}_{w}$'s act. (Remember
that the $\mathcal{A}_{G}$-tree $\mathcal{T}_{G}$ is identical to $\mathcal{T%
}_{2N},$ whenever $\Bbb{G}$ is a fractaloid!) Since the graph groupoid $\Bbb{%
G}$ is a fractaloid, $\mathcal{E}\left( \mathcal{F}_{v_{k},\,n}\right) $
create the full-subgraphs $G_{v_{k},\,n}$, having its finite path set $%
FP(G_{v_{k},\,n})$ $=$ $\mathcal{E}\left( \mathcal{F}_{v_{k},\,n}\right) $
in $\mathcal{T}_{2N}$, and they are graph-isomorphic, for $k$ $=$ $1,$ $2,$
and for $n$ $\in $ $\Bbb{N}.$ Therefore, the generating sets $\mathcal{F}%
_{v_{k},\,n},$ which are regarded as the edge sets of the full-subgraphs $%
G_{v_{k},\,n},$ should have the same cardinality.. i.e.,

\strut

\begin{center}
$\left| \mathcal{E}\left( \mathcal{F}_{v_{1},\,n}\right) \right| $ $=$ $%
\left| \mathcal{E}\left( \mathcal{F}_{v_{2},\,n}\right) \right| $ $%
\Longleftrightarrow $ $\left| \mathcal{F}_{v_{1},\,n}\right| $ $=$ $\left| 
\mathcal{F}_{v_{2},\,n}\right| ,$
\end{center}

\strut

by the bijectivity of $\mathcal{E},$ for all $n$ $\in $ $\Bbb{N}.$ Since $%
v_{1}$ and $v_{2}$ are arbitrary, for any $n$ $\in $ $\Bbb{N},$ the
cardinalities $\left| \mathcal{F}_{v,\,n}\right| $ of the length-$n$-$v$%
-subsets of $FP(G^{\symbol{94}})$ are identical, for all $v$ $\in $ $V(G).$
\end{proof}

\strut

By the previous lemma, we can re-formulate the formula (5.1), where the
given labeled graph groupoid $\Bbb{G}$ is a fractaloid.

\strut

\begin{theorem}
Let $G$ be a canonical weighted graph with its labeled graph groupoid $\Bbb{G%
},$ and assume that $\Bbb{G}$ is a fractaloid. Then the labeling operator $%
T_{G}$ of $\Bbb{G}$ in $\Bbb{M}_{G}$ satisfies that

\strut

(5.3)$\qquad \qquad \qquad \qquad \qquad E\left( T_{G}^{n}\right) $ $=$ $%
\left| \sum_{n}^{(N)}\right| \cdot 1_{\Bbb{D}_{G}},$

where

(5.4)$\qquad \qquad \sum_{n}^{(N)}$ $=$ $\left\{ (t_{i_{1}},\text{ ..., }%
t_{i_{n}})\in \left( \pm X\right) ^{n}\left| \sum_{j=1}^{n}\text{ }%
t_{i_{j}}=x_{0}\right. \right\} ,$

\strut

for all $n$ $\in $ $\Bbb{N}.$
\end{theorem}

\strut

\begin{proof}
By (5.2) and by the previous lemma, we can compute, for any fixed $n$ $\in $ 
$\Bbb{N},$

\strut

\begin{center}
$
\begin{array}{ll}
E(T_{G}^{n}) & =\underset{(i_{1},...,i_{n})\in \{\pm 1,\text{ ..., }\pm
N\}^{n}}{\sum }\left( \underset{e_{1}...e_{n}\in FP(G^{\symbol{94}}),\,%
\widetilde{\omega }(e_{j})=t_{i_{j}},\,\sum_{k=1}^{n}t_{i_{k}}=x_{0}\,}{\sum 
}\text{ }R_{e_{1}...e_{n}}\right) \\ 
&  \\ 
& =\underset{v\in V(G)}{\sum }\left| \mathcal{F}_{v,\,n}\right| R_{v}\text{ }%
=\text{ }\underset{v\in V(G)}{\sum }\text{ }\eta _{n}\text{ }R_{v}\text{ }=%
\text{ }\eta _{n}\text{ }\left( \underset{v\in V(G)}{\sum }\text{ }%
R_{v}\right) .
\end{array}
$
\end{center}

\strut

The last two equalities hold, by the previous lemma. i.e., for the fixed $n$ 
$\in $ $\Bbb{N},$ the cardinalities $\left| \mathcal{F}_{v,\,n}\right| $ of
the length-$n$-$v$-subsets $\mathcal{F}_{v,\,n}$ are identical to the number 
$\eta _{n}$, for all $v$ $\in $ $V(G),$ whenever $\Bbb{G}$ is a fractaloid.
Recall that $\underset{v\in V(G)}{\sum }$ $R_{v}$ is the identity element $%
1_{\Bbb{D}_{G}}$ $=$ $1_{\Bbb{M}_{G}}$ of $\Bbb{D}_{G}$ $\subseteq $ $\Bbb{M}%
_{G}$. Therefore, we have that

\strut

\begin{center}
$E(T_{G}^{n})$ $=$ $\eta _{n}$ $\cdot $ $1_{\Bbb{D}_{G}}$ $\in $ $\Bbb{D}%
_{G},$ for all $n$ $\in $ $\Bbb{N},$
\end{center}

\strut

for some scalar-values $\{\eta _{n}\}_{n=1}^{\infty }.$ Notice that

\strut

\begin{center}
$\eta _{n}$ $=$ $\left| \mathcal{F}_{v,\,n}\right| ,$ for all $v$ $\in $ $%
V(G),$
\end{center}

\strut

for any $n$ $\in $ $\Bbb{N}.$ By the definition of length-$n$-vertex-subsets 
$\mathcal{F}_{v,\text{ }n}$'s of $FP(G^{\symbol{94}}),$ we can construct
sets $\sum_{n}^{(N)}$ (independent from the choice of vertices) by

\strut

\begin{center}
$\sum_{n}^{(N)}$ $\overset{def}{=}$ $\left\{ (t_{i_{1}},\text{ ..., }%
t_{i_{n}})\in \left( \pm X\right) ^{n}\left| \sum_{j=1}^{n}\text{ }%
t_{i_{j}}=x_{0}\right. \right\} ,$
\end{center}

\strut

\strut for all $n$ $\in $ $\Bbb{N}.$ Then the cardinality $\left|
\sum_{n}^{(N)}\right| $ of $\sum_{n}^{(N)}$ is identical to the value $\eta
_{n},$ for all $n$ $\in $ $\Bbb{N}.$ \strut
\end{proof}

\strut

\begin{remark}
In general, without the assumption that $\Bbb{G}$ is a fractaloid, the
equalities

\strut

\begin{center}
$\eta _{n}$ $=$ $\left| \mathcal{F}_{v,\,n}\right| $\quad and\quad $\eta
_{n} $ $=$ $\left| \sum_{n}^{(N)}\right| ,$
\end{center}

\strut

in the proof of the previous theorem, do not hold. (See [40].)
\end{remark}

\strut \strut

The above theorem shows that the $\Bbb{D}_{G}$-valued (or the
operator-valued) moments $E(T_{G}^{n})$ of the labeling operator $T_{G}$ are
completely determined by the scalar-values $\left| \sum_{n}^{(N)}\right| $
which are the cardinalities of certain sets $\sum_{n}^{(N)},$ for all $n$ $%
\in $ $\Bbb{N}$, whenever the labeled graph groupoid $\Bbb{G}$ is a
fractaloid. Also, the above theorem shows how the labeling (in terms of the
admissibility of $\Bbb{G}$) works (in $\Bbb{G}$) and $\Bbb{M}_{G},$ whenever 
$\Bbb{G}$ is a fractaloid.

\strut

\strut

\strut

\subsection{Refinements of (5.3)}

\strut

\strut

In Section 5.3, we showed that the $\Bbb{D}_{G}$-valued free moments $%
\{E(T_{G}^{n})\}_{n=1}^{\infty }$ of the labeling operator $T_{G}$ of a
fractaloid $\Bbb{G}$ in the right graph von Neumann algebra $\Bbb{M}_{G}$ is
completely determined by the cardinalities $\{\left| \sum_{n}^{(N)}\right|
\}_{n=1}^{\infty }$ of certain sets $\{\sum_{n}^{(N)}\}_{n=1}^{\infty },$
where

\strut

\begin{center}
$\sum_{n}^{(N)}$ $=$ $\{(t_{i_{1}},$ ..., $t_{i_{n}})$ $\in $ $(\pm X)^{n}$ $%
:$ $\sum_{j=1}^{n}$ $t_{i_{j}}$ $=$ $x_{0}\},$
\end{center}

\strut where

\begin{center}
$t_{i_{j}}$ $\overset{def}{=}$ $\left\{ 
\begin{array}{ll}
x_{i_{j}} & \text{if }t_{i_{j}}\in X \\ 
-x_{i_{j}} & \text{if }t_{i_{j}}\in -X,
\end{array}
\right. $
\end{center}

\strut

for all $j$ $=$ $1,$ ..., $n,$ for all $n$ $\in $ $\Bbb{N},$ where $X$ $=$ $%
\{x_{1},$ ..., $x_{N}\}$ $\subset $ $\Bbb{C}^{\oplus \,N}$ is the labeling
set of $G$ consisting of the orthonormal vectors $x_{j}$ $=$ $(0,$ ..., $0,$ 
$\underset{j\text{-th}}{1},$ $0,$ ..., $0),$ for $j$ $=$ $1,$ ..., $n,$ and
where

\strut

\begin{center}
$x_{0}$ $=$ $\left( \underset{N\text{-times}}{\underbrace{0,\text{ ........, 
}0}}\right) $ $\in $ $\Bbb{C}^{\oplus \,N},$
\end{center}

where

\begin{center}
$N$ $=$ $\max \{\deg _{out}(v)$ $:$ $v$ $\in $ $V(G)\}.$
\end{center}

\strut

By using the operation $\theta ,$ defined in Section 3.2, we can re-write $%
\sum_{n}^{(N)}$'s by

\strut

\begin{center}
$\sum_{n}^{(N)}$ $=$ $\{w$ $\in $ $(\pm X)^{n}$ $:$ $\theta (w)$ $=$ $%
x_{0}\}.$
\end{center}

\strut

Construct now the lattice paths in $\Bbb{R}^{2}.$

\strut

\begin{definition}
Let $Y$ $=$ $\{1,$ ..., $N\},$ for $N$ $\in $ $\Bbb{N}.$ For the given index
set $Y,$ define the vectors $l_{\pm k}$ by the vectors in $\Bbb{R}^{2}$
satisfying that

\strut

\begin{center}
$l_{k}$ $=$ $(1,$ $e^{k})$ and $l_{-k}$ $=$ $(1,$ $-e^{k}),$ for all $k$ $=$ 
$1,$ ..., $N,$
\end{center}

\strut

where $e$ $\in $ $\Bbb{R}$ is the natural (exponential) number. i.e., $l_{k}$
(resp. $l_{-k}$) is a vector starting from the point $(0,$ $0)$ and ending
at the point $(1,$ $e^{k})$ (resp. the point $(1,$ $-e^{k})$), for all $k$ $%
= $ $1,$ ..., $N.$ Define now the binary operation, denoted simply by ($%
\cdot $), on $\{l_{\pm k}$ $:$ $k$ $=$ $1,$ ..., $N\}$, by

\strut

\begin{center}
$l_{i}$ $l_{j}$ $\overset{def}{=}$ the path in $\Bbb{R}^{2}$ connecting the
vector $l_{i}$ and $l_{j},$ by identifying the ending point $(1,$ $e^{i})$
of $l_{i}$ and the starting point $(0,$ $0)$ of $l_{j},$
\end{center}

\strut

where $i,$ $j$ $\in $ $\{\pm 1,$ ..., $\pm N\}.$ Inductively, we can
determine the paths $l_{i_{1}}$ $l_{i_{2}}$ ... $l_{i_{n}},$ where $(i_{1},$
..., $i_{n})$ $\in $ $\{\pm 1,$ ..., $\pm N\}^{n},$ for all $n$ $\in $ $\Bbb{%
N}.$ Such paths in $\Bbb{R}^{2},$ generated by the vectors $l_{\pm 1},$ ..., 
$l_{\pm N},$ are called the lattice paths generated by $Y.$ We denote the
collection of all lattice paths generated by the vectors $l_{\pm 1},$ ..., $%
l_{\pm N}$, by $\mathcal{L}_{N}.$ And we call $\mathcal{L}_{N},$ the lattice
path set induced by $Y.$

\strut

Let $l$ $=$ $l_{i_{1}}$ ... $l_{i_{n}}$ $\in $ $\mathcal{L}_{N}$ be a
lattice path. Define the length of $l$ by the number $n.$ And we denote the
length of $l$ by $\left| l\right| .$
\end{definition}

\strut \strut \strut \strut \strut \strut

Notice that every element in the lattice path set $\mathcal{L}_{N}$ is a
(non-reduced and nonempty) words in $\{l_{\pm 1},$ ..., $l_{\pm N}\}.$ We
can decompose $\mathcal{L}_{N}$ by

\strut

\begin{center}
$\mathcal{L}_{N}$ $=$ $\underset{k=1}{\overset{\infty }{\sqcup }}$ $\mathcal{%
L}_{N}(k)$
\end{center}

with

\begin{center}
$\mathcal{L}_{N}(k)$ $\overset{def}{=}$ $\{l$ $\in $ $\mathcal{L}_{N}$ $:$ $%
\left| l\right| $ $=$ $k\},$ for all $k$ $\in $ $\Bbb{N}$,
\end{center}

\strut

where ``$\sqcup $'' means the disjoint union.\strut

\strut \strut

\begin{remark}
Remark that, we define the lattices $l_{\pm k}$ by the vectors $(1,$ $\pm
e^{k}),$ for $k$ $=$ $1,$ ..., $N.$ The choice of the $y$-coordinates $\pm
e^{k}$ (by regarding the space $\Bbb{R}^{2}$ as the usual $(x,$ $y)$%
-coordinate system) is for our purpose, for $k$ $=$ $1,$ ..., $N$. i.e., we
want to make the sets $\mathcal{L}_{N}(k)$ be equipotent (or bijective) to
the sets $\sum_{k}^{(N)},$ for each $k$ $=$ $1,$ ..., $N,$ where $%
\sum_{k}^{(N)}$ are the given subsets of $\Bbb{C}^{\oplus \,N},$ in the
previous section. It is okay to take $\pi $ instead of $e.$
\end{remark}

\strut \strut

Consider a certain subset $\mathcal{L}_{N}^{o}(k)$ of $\mathcal{L}_{N}(k),$
for each $k$ $\in $ $\Bbb{N}$: Define the subset $\mathcal{L}_{N}^{o}(k)$ of 
$\mathcal{L}_{N}(k)$ by

\strut

\begin{center}
$\mathcal{L}_{N}^{o}(k)$ $\overset{def}{=}$ $\left\{ l\in \mathcal{L}%
_{N}(k):\left| 
\begin{array}{c}
l\text{ ends at the point on } \\ 
\text{the horizontal axis} \\ 
\text{(or the }x\text{-axis) of }\Bbb{R}^{2}
\end{array}
\right. \right\} .$
\end{center}

\strut

i.e., if $l$ $=$ $l_{i_{1}}$ ... $l_{i_{n}}$ $\in $ $\mathcal{L}_{N}^{o}(k),$
then it is a lattice path starting from $(0,$ $0),$ ending at $(k,$ $0).$

\strut

\begin{definition}
If a lattice path $l$ is contained in $\mathcal{L}_{N}^{o}(k),$ for some $k$ 
$\in $ $\Bbb{N},$ then we say that the length-$k$ lattice path $l$ \emph{%
satisfies the }(\emph{horizontal- or }$x$\emph{-})\emph{axis property}.
\end{definition}

\strut

Similarly, we can define the subset $\mathcal{L}_{N}^{o}$ of the lattice
path set $\mathcal{L}_{N}$ by

\strut

\begin{center}
$\mathcal{L}_{N}^{o}$ $\overset{def}{=}$ $\underset{k=1}{\overset{\infty }{%
\cup }}$ $\mathcal{L}_{N}^{o}(k).$
\end{center}

\strut

i.e., the subset $\mathcal{L}_{N}^{o}$ is the collection of all lattice
paths satisfying the axis property.

\strut

By defining a map $\Theta $ $:$ $\sum_{n}^{(N)}$ $\rightarrow $ $\mathcal{L}%
_{N}^{o}(n),$ we can realize the relation between the sets $\sum_{n}^{(N)}$
and $\mathcal{L}_{N}^{o}(n),$ for all $n$ $\in $ $\Bbb{N}.$ Define the map $%
\Theta $ by

\strut

\begin{center}
$(t_{i_{1}},$ ..., $t_{i_{n}})$ $\in $ $\sum_{n}^{(N)}$ $\overset{\Theta }{%
\longrightarrow }$ $l_{k_{1}}$ ... $l_{k_{n}}$ $\in $ $\mathcal{L}%
_{N}^{o}(n),$
\end{center}

where

\begin{center}
$k_{j}$ $=$ $\left\{ 
\begin{array}{lll}
i_{j} &  & \text{if }t_{i_{j}}=x_{i_{j}}\in X \\ 
&  &  \\ 
-i_{j} &  & \text{if }t_{i_{j}}=-x_{i_{j}}\in -X,
\end{array}
\right. $
\end{center}

\strut

for all $j$ $=$ $1,$ ..., $n,$ for all $n$ $\in $ $\Bbb{N}.$ By the very
definition of the map $\Theta ,$ this map is bijective. Therefore, we can
realize that the property

\strut

\begin{center}
$\theta \left( (t_{i_{1}},\text{ ..., }t_{i_{n}})\right) $ $=$ $%
\sum_{j=1}^{n}$ $t_{i_{j}}$ $=$ $x_{0}$ $\in $ $\Bbb{C}^{\oplus \,N}$
\end{center}

\strut

is equivalent to the axis property: the lattice path

\strut

\begin{center}
$\Theta \left( (t_{i_{1}},\text{ ..., }t_{i_{n}})\right) $ $=$ $l_{k_{1}}$
... $l_{k_{n}}$ $\in $ $\mathcal{L}_{N}^{o}(n)$
\end{center}

\strut

satisfies the axis property. By the existence of the above bijection $\Theta
,$ we have the following proposition.

\strut

\begin{proposition}
The sets $\sum_{n}^{(N)}$ and $\mathcal{L}_{N}^{o}(n)$ are equipotent (or
bijective), for all $n$ $\in $ $\Bbb{N}.$ i.e.,

\strut

\begin{center}
$\left| \sum_{n}^{(N)}\right| $ $=$ $\left| \mathcal{L}_{N}^{o}(n)\right| ,$
for all $n$ $\in $ $\Bbb{N}.$
\end{center}

$\square $
\end{proposition}

\strut

The following corollary is the direct consequence of the previous
proposition.

\strut

\begin{corollary}
Let $G$ be a canonical weighted graph with its labeled graph groupoid $\Bbb{G%
},$ and assume that $\Bbb{G}$ is a fractaloid. Then the labeling operator $%
T_{G}$ of $\Bbb{G}$ in $\Bbb{M}_{G}$ has its $\Bbb{D}_{G}$-valued moments,

\strut

\begin{center}
$E(T_{G}^{n})$ $=$ $\left| \mathcal{L}_{N}^{o}(n)\right| $ $\cdot $ $1_{\Bbb{%
D}_{G}},$
\end{center}

\strut

for all $n$ $\in $ $\Bbb{N}.$ $\square $
\end{corollary}

\strut

In fact, the length-$n$ lattice path set $\mathcal{L}_{N}^{o}(n)$ satisfying
the axis property is introduced for visualizing the set $\sum_{n}^{(N)},$
for all $n$ $\in $ $\Bbb{N}.$ By identifying (bijectively) the sets $%
\{\sum_{n}^{(N)}\}_{n=1}^{\infty }$ and $\{\mathcal{L}_{N}^{o}(n)\}_{n=1}^{%
\infty },$ we can get the following proposition.

\strut

\begin{proposition}
Let $T_{G}$ be the labeling operator of a fractaloid $\Bbb{G}$ in $\Bbb{M}%
_{G}.$ Then all odd $\Bbb{D}_{G}$-valued moments vanish. i.e., $E(T_{G}^{n})$
$=$ $0_{\Bbb{D}_{G}},$ whenever $n$ is odd.
\end{proposition}

\strut

\begin{proof}
We can easily verify that if $n$ is odd, then the set $\mathcal{L}%
_{N}^{o}(n) $ is empty. So, if $n$ is odd, then $\left| \mathcal{L}%
_{N}^{o}(n)\right| $ $=$ $0.$ So, we have that

\strut

\begin{center}
$E(T_{G}^{n})$ $=$ $\left| \sum_{n}^{(N)}\right| \cdot 1_{\Bbb{D}_{G}}$ $=$ $%
\left| \mathcal{L}_{N}^{o}(n)\right| \cdot 1_{\Bbb{D}_{G}}$ $=$ $0_{\Bbb{D}%
_{G}},$
\end{center}

\strut

whenever $n$ is odd in $\Bbb{N}.$\strut
\end{proof}

\strut \strut

By the previous proposition, we can get the more detailed $\Bbb{D}_{G}$-free
distributional data of the labeling operator $T_{G}$ of the fractaloid $\Bbb{%
G}$:

\strut

\begin{corollary}
Let $T_{G}$ be the labeling operator of a fractaloid $\Bbb{G}$ in the right
graph von Neumann algebra $\Bbb{M}_{G}.$ Then

\strut

\strut (5.5)

\begin{center}
$E(T_{G}^{n})$ $=$ $\left\{ 
\begin{array}{ll}
\left| \sum_{n}^{(N)}\right| \cdot 1_{\Bbb{D}_{G}}=\left| \mathcal{L}%
_{N}^{o}(n)\right| \cdot 1_{\Bbb{D}_{G}} & \text{if }n\in 2\Bbb{N} \\ 
&  \\ 
0_{\Bbb{D}_{G}} & \text{if }n\in 2\Bbb{N}-1,
\end{array}
\right. $
\end{center}

\strut

for all $n$ $\in $ $\Bbb{N}.$ $\square $
\end{corollary}

\strut

Now, we concentrate on finding the cardinality $\left| \mathcal{L}%
_{N}^{o}(2n)\right| ,$ for $n$ $\in $ $\Bbb{N}.$ To do that, we consider the
results in [38] and [39]. In [38], we compute the powers of multinomials in
a commutative unital algebra. We showed that the computations are determined
by the well-known famous Pascal's triangle, recursively.

\strut

Let $A$ be a commutative unital algebra with its identity $1_{A},$ and let $%
a_{1},$ ..., $a_{n}$ be arbitrary distinct elements in $A,$ for $n$ $\in $ $%
\Bbb{N}.$ Then we can construct a new element $a$ $=$ $\sum_{k=1}^{n}$ $%
a_{j} $ in $A,$ and the powers $a^{m}$ of $a$ in $A,$ for all $m$ $\in $ $%
\Bbb{N}.$ We will call the element $a$ the \emph{multinomial induced by} $%
a_{1},$ ..., $a_{n},$ and the elements $a^{m},$ the $m$-\emph{th powers} 
\emph{of} $a,$ for all $m$ $\in $ $\Bbb{N}.$ Then, by the commutativity on $%
A,$ we can get that

\strut

\begin{center}
$
\begin{array}{ll}
a^{m} & =\underset{(i_{1},\text{ ..., }i_{m})\in \{1,...,n\}^{m}}{\sum }%
\text{ }a_{i_{1}}\text{ ... }a_{i_{m}} \\ 
&  \\ 
& =\underset{(i_{1},...,i_{m})\in \{1,...,n\}^{m},\,i_{1}\leq \,i_{2}\,\leq 
\text{ ... }\leq \text{ }i_{m}}{\sum }\text{ }c_{i_{1},...,i_{m}}^{(n)}\text{
}a_{i_{1}}...a_{i_{m}},
\end{array}
$
\end{center}

\strut

where $c_{i_{1},...,i_{m}}^{(n)}$ $\in $ $\Bbb{N}$ are called the $(i_{1},$
..., $i_{m})$-\emph{th coefficients of }$a_{i_{1}}$ ... $a_{i_{m}}$ \emph{in}
$a^{m},$ for $(i_{1},$ ..., $i_{m})$ $\in $ $\{1,$ ..., $n\}^{m},$ where $%
i_{1}$ $\leq $ $i_{2}$ $\leq $ ... $\leq $ $i_{m}.$ The existence of the
coefficients of $a^{m}$ is guaranteed by the commutativity of $A,$ for all $%
m $ $\in $ $\Bbb{N}.$ It is natural that if we know how to compute all
coefficients of $a^{m},$ then we can compute $a^{m}.$

\strut

For convenience, let's denote $1_{A}$ be $a_{0}.$ Fix $a_{j}$ $\in $ $%
\{a_{1},$ ..., $a_{n}\}.$ Then we can consider the binomial $y_{j}$ $=$ $%
1_{A}$ $+$ $a_{j}$ $=$ $a_{0}$ $+$ $a_{j}$ in $A,$ for all $j$ $=$ $1,$ ..., 
$n.$ Then the powers $y_{j}^{m}$ have their coefficients completely
determined by the Pascal's triangle, for all $j$ $=$ $1,$ ..., $n,$ and $m$ $%
\in $ $\Bbb{N}.$ i.e., we can write

\strut

\begin{center}
$y_{j}^{m}$ $=$ $\underset{(i_{1},...,i_{m})\in \{0,\text{ }j\}^{m},\text{ }%
i_{1}\leq \text{ }i_{2}\text{ }\leq \text{ ... }\leq \text{ }i_{m}}{\sum }$ $%
\varepsilon _{j:\,i_{1},...,i_{m}}$ $a_{i_{1}}$ ... $a_{i_{m}},$
\end{center}

\strut

where $\varepsilon _{j:i_{1},...,i_{m}}$'s are the $(i_{1},$ ..., $i_{m})$%
-th coefficients of $y_{j}^{m},$ for all $j$ $=$ $1,$ ..., $n,$ and $m$ $\in 
$ $\Bbb{N}$, and they are the entries of the Pascal's triangle,

\strut

\begin{center}
$
\begin{array}{lllllllllll}
&  &  &  &  & 1 &  &  &  &  &  \\ 
&  &  &  & 1 &  & 1 &  &  &  &  \\ 
&  &  & 1 &  & 2 &  & 1 &  &  &  \\ 
&  & 1 &  & 3 &  & 3 &  & 1 &  &  \\ 
& 1 &  & 4 &  & 6 &  & 4 &  & 1 &  \\ 
1 &  & 5 &  & 10 &  & 10 &  & 5 &  & 1 \\ 
&  &  &  &  &  &  &  &  &  &  \\ 
\vdots &  & \vdots &  & \vdots &  & \vdots &  & \vdots &  & \vdots
\end{array}
.$
\end{center}

\strut

In other words, the entry $\varepsilon _{j:\,\underset{k\text{-times}}{%
\underbrace{0,........,0}}\,,\,\,\underset{(m-k)\text{-times}}{\underbrace{%
j,..........,j}}}$ are the $(m$ $-$ $k)$-th entry (from the left) of the $m$
-th level of the Pascal's triangle. Equivalently, the above triangle can be
re-written as follows, in terms of the coefficients of $y_{j}^{m}$:

\strut

\begin{center}
$
\begin{array}{lllllllll}
&  &  &  & \varepsilon _{j:0} &  &  &  &  \\ 
&  &  & \varepsilon _{j:0} &  & \varepsilon _{j:j} &  &  &  \\ 
&  & \varepsilon _{j:0,0} &  & \varepsilon _{j:0,\,j} &  & \varepsilon
_{j:\,j,\,j} &  &  \\ 
& \varepsilon _{j:0,0,0} &  & \varepsilon _{j:0,0,j} &  & \varepsilon
_{j:0,j,j} &  & \varepsilon _{j:j,\,j,\,j} &  \\ 
\varepsilon _{j:0,0,0,0} &  & \varepsilon _{j:0,0,0,j} &  & \varepsilon
_{j:0,0,j,j} &  & \varepsilon _{j:0,j,j,j} &  & \varepsilon _{j:j,j,j,\,j}
\\ 
&  &  &  &  &  &  &  &  \\ 
\vdots &  & \vdots &  & \vdots &  & \vdots &  & \vdots
\end{array}
,$
\end{center}

\strut

for all $j$ $=$ $1,$ ..., $n,$ and for $m$ $\in $ $\Bbb{N}.$ The following
proposition is the main result of [38]:

\strut

\begin{proposition}
(See [38]) Let $a$ $=$ $\sum_{j=1}^{n}$ $a_{j}$ be a multinomial induced by
the distinct elements $a_{1},$ ..., $a_{n}$ in a commutative unital algebra $%
A,$ for $n$ $\in $ $\Bbb{N}.$ Then the powers $a^{m}$ of $a$ has the
coefficients $c_{i_{1},...,i_{m}}^{(n)},$ satisfying that

\strut \strut

\begin{center}
$c_{i_{1},...,i_{m}}^{(n)}$ $=$ $\left( c_{i_{1},...,i_{m-k}}^{(n-1)}\right)
\left( \varepsilon _{m:\text{ }\underset{k\text{-times}}{\underbrace{%
0,........,0}}\text{ },\text{ }i_{m-k+1},...,i_{m}}\right) ,$
\end{center}

\strut

for $k$ $\in $ $\{1,$ ..., $m\}$, for all $m$ $\in $ $\Bbb{N},$ whenever

\strut

\begin{center}
$i_{m-k+1}$ $=$ $i_{m-k+2}$ $=$ ... $=$ $i_{m}$ in $\{1,$ ..., $n\}.$
\end{center}

$\square $
\end{proposition}

\strut

The proof of the above proposition is straightforward. The above proposition
shows that the coefficients of the $m$-th power $a^{m}$ of a multinomial $a$ 
$=$ $\sum_{j=1}^{n}$ $a_{j}$ are determined by the entries of the Pascal's
triangle, recursively. The following three examples are helpful to
understand the use of the above proposition.\strut

\strut

\begin{example}
The coefficient $c_{1,1,1,1,3,3,4}^{(4)}$ of $a_{1}^{4}$ $a_{3}^{2}$ $a_{4}$
in $(\sum_{k=1}^{4}$ $a_{k})^{7}$ is

\strut

\begin{center}
$
\begin{array}{ll}
c_{1,1,1,1,3,3,4}^{(4)} & =c_{1,1,1,1,3,3}^{(3)}\text{ }\varepsilon
_{4:0,0,0,0,0,0,4}=7\cdot c_{1,1,1,1,3,3}^{(3)} \\ 
&  \\ 
& =7\cdot c_{1,1,1,1}^{(2)}\text{ }\varepsilon _{3:0,0,0,0,3,3}=7\cdot
15\cdot 1=105
\end{array}
$
\end{center}
\end{example}

\strut

\begin{example}
The coefficient $c_{2,2,3,4}^{(5)}$ of $a_{2}^{2}$ $a_{3}$ $a_{4}$ in $%
(\sum_{k=1}^{5}$ $a_{k})^{4}$ is

\strut

\begin{center}
$
\begin{array}{ll}
c_{2,2,3,4}^{(5)} & =c_{2,2,3}^{(4)}\cdot \varepsilon _{5:0,0,0,4}=4\cdot
c_{2,2,3}^{(4)} \\ 
&  \\ 
& =4\cdot c_{2,2}^{(3)}\cdot \varepsilon _{4:0,0,3}=4\cdot 1\cdot 3=12.
\end{array}
$
\end{center}

\strut
\end{example}

\strut

\begin{example}
The coefficient $c_{1,2,3,4}^{(4)}$ of $a_{1}$ $a_{2}$ $a_{3}$ $a_{4}$ in $%
(\sum_{k=1}^{4}$ $a_{k})^{4}$ is

\strut

\begin{center}
$
\begin{array}{ll}
c_{1,2,3,4}^{(4)} & =c_{1,2,3}^{(3)}\cdot \varepsilon _{4:0,0,0,4}=4\cdot
c_{1,2,3}^{(3)} \\ 
&  \\ 
& =4\cdot c_{1,2}^{(2)}\cdot \varepsilon _{3:0,0,3}=4\cdot 3\cdot
c_{1,2}^{(2)}=4\cdot 3\cdot 2=24.
\end{array}
$
\end{center}

\strut
\end{example}

\strut

Now, consider the lattice paths in $\mathcal{L}_{N}(2n)$. Every lattice path 
$l$ contained in $\mathcal{L}_{N}(2n)$ is the length-$2n$ (non-reduced and
nonempty) word in $\{l_{\pm 1},$ ..., $l_{\pm N}\},$ say $l$ $=$ $l_{j_{1}}$ 
$l_{j_{2}}$ ... $l_{j_{2n}}.$ As we observed before, we can realize that
this lattice $l$ satisfies the axis property if and only if

\strut

\begin{center}
$\sum_{k=1}^{2n}$ $j_{k}$ $=$ $0.$
\end{center}

\strut

Assume that we have distinct elements $a_{1},$ $a_{-1},$ $a_{2},$ $a_{-2},$
..., $a_{N},$ $a_{-N}$ in a commutative unital algebra $A.$ Then we can have
an element

$\strut $

\begin{center}
$a_{\pm }$ $=$ $\sum_{k=-N}^{-1}$ $a_{k}$ $+$ $\sum_{i=1}^{N}$ $a_{i},$ in $%
A.$
\end{center}

\strut

Then the coefficients, denoted by $c_{i_{1},...,i_{2n}}^{(\pm N)}$ of the $%
2n $-th power $a_{\pm }^{2n}$ of $a_{\pm }$ are also determined similarly by
the above proposition, for all $(i_{1},$ ..., $i_{2n})$ $\in $ $\{\pm 1,$
..., $\pm N\}^{2n},$ for $n$ $\in $ $\Bbb{N}.$ By the previous observation,
we can realize that the addition of all coefficients

\strut

\begin{center}
$c_{i_{1},...,i_{2n}}^{(\pm N)}$ satisfying that $\sum_{k=1}^{2n}$ $i_{k}$ $%
= $ $0,$
\end{center}

\strut

represents the cardinalities of the length-$2n$ lattice paths satisfying the
axis property. Therefore, we can get the following theorem.

\strut

\begin{theorem}
Let $\mathcal{L}_{N}^{o}(2n)$ be the set of all length-$2n$ lattice paths
induced by $l_{\pm 1},$ ..., $l_{\pm N}$, for any $N$, $n$ $\in $ $\Bbb{N}.$
And assume that $c_{i_{1},...,i_{2n}}^{(\pm N)}$ are the coefficients of the 
$2n$-th power $a_{\pm }^{2n},$ for $n$ $\in $ $\Bbb{N},$ where $a_{\pm }$ is
defined as in the previous paragraph. Then

\strut

\begin{center}
$\left| \mathcal{L}_{N}^{o}(2n)\right| $ $=$ $\underset{(j_{1},...,j_{2n})%
\in \{\pm 1,\text{ ..., }\pm N\}^{2n},\,\sum_{k=1}^{2n}\,j_{k}=0}{\sum }$ $%
c_{j_{1},...,j_{2n}}^{(\pm N)},$
\end{center}

\strut

for all $n$ $\in $ $\Bbb{N}.$ $\square $
\end{theorem}

\strut

The above theorem shows that the amalgamated moments of the labeling
operator of fractaloids are also determined by certain entries of the
Pascal's triangle, recursively.

\strut

\begin{corollary}
Let $c_{i_{1},..,i_{n}}^{(\pm N)}$ be the coefficients of the $2n$-power $%
a_{\pm }^{2n}$ of $a_{\pm }$, defined as above. And let $T_{G}$ be the
labeling operator of a fractaloid $\Bbb{G}$ in the right graph von Neumann
algebra $\Bbb{M}_{G},$ where $N$ $=$ $\max \{$ $\deg _{out}(v)$ $:$ $v$ $\in 
$ $V(G)\}.$ Then the\textbf{\ nonzero} $\Bbb{D}_{G}$-valued moments $%
\{E(T_{G}^{2n})\}_{n=1}^{\infty }$ of $T_{G}$ are determined by

\strut

\begin{center}
$E\left( T_{G}^{2n}\right) $ $=$ $\left( \underset{(j_{1},...,j_{2n})\in
\{\pm 1,\text{ ..., }\pm N\}^{2n},\,\sum_{k=1}^{2n}\,j_{k}=0}{\sum }%
c_{j_{1},...,j_{2n}}^{(\pm N)}\right) \cdot 1_{\Bbb{D}_{G}},$
\end{center}

\strut

for all $n$ $\in $ $\Bbb{N}.$ $\square $
\end{corollary}

\strut

\strut \strut \strut

\strut

\subsection{An Example: Computing$\left| \mathcal{L}_{1}^{o}(2n)\right| $ $=$
$\left| \sum_{2n}^{(1)}\right| ,$ for $n$ $\in $ $\Bbb{N}$}

\strut

\strut

In this section, we will provide an example to use the lattice path model $\{%
\mathcal{L}_{N}^{o}(n)\}_{n=1}^{\infty }$ to compute $\{\left|
\sum_{n}^{(N)}\right| \}_{n=1}^{\infty }.$ Again, remark that

\strut

\begin{center}
$E(T_{G}^{n})$ $=$ $\left| \sum_{n}^{(N)}\right| \cdot 1_{D_{G}},$ for all $%
n $ $\in $ $\Bbb{N},$
\end{center}

\strut

and the sets $\sum_{n}^{(N)}$ are visualized by the length-$n$ lattice paths
contained in $\mathcal{L}_{N}^{o}(n)$ satisfying the axis property, for all $%
n$ $\in $ $\Bbb{N}$, \strut where

\strut

\begin{center}
$N$ $=$ $\max \{\deg _{out}(v)$ $:$ $v$ $\in $ $V(G)\}$.
\end{center}

\strut \strut

In the rest of this section, as an example of Section 5.4, we will compute

$\strut $

\begin{center}
$\left| \mathcal{L}_{1}^{o}(2n)\right| $ $=$ $\left| \sum_{2n}^{(1)}\right|
, $ for all $n$ $\in $ $\Bbb{N}.$
\end{center}

\strut

(\textbf{Case:} $n$ $=$ $2$) If $n$ $=$ $2,$ then we have the following
elements of $\mathcal{L}_{1}^{o}(2)$:

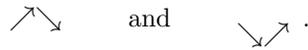
\strut

\begin{center}
$
\begin{array}{l}
\nearrow \searrow
\end{array}
$\qquad and\qquad $
\begin{array}{l}
\\ 
\searrow \nearrow
\end{array}
.$
\end{center}

\strut And hence

\begin{center}
$\left| \mathcal{L}_{1}^{o}(2)\right| $ $=$ $\left| \sum_{2}^{(1)}\right| $ $%
=$ $2.$
\end{center}

\strut

(\textbf{Case}: $n$ $=$ $4$) If $n$ $=$ $4,$ then $\mathcal{L}_{1}^{o}(4)$
has the following elements:

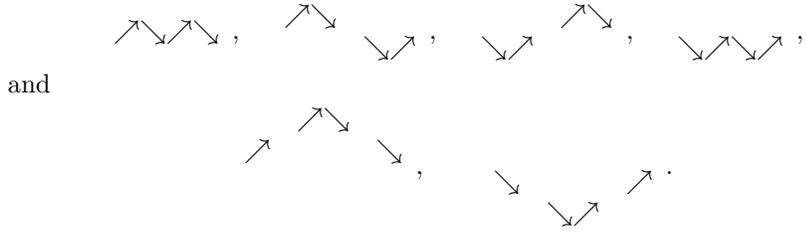
\strut

\begin{center}
$
\begin{array}{l}
\nearrow \searrow \nearrow \searrow
\end{array}
,\quad 
\begin{array}{ll}
\nearrow \searrow &  \\ 
& \searrow \nearrow
\end{array}
,\quad 
\begin{array}{ll}
& \nearrow \searrow \\ 
\searrow \nearrow & 
\end{array}
,\quad 
\begin{array}{l}
\\ 
\searrow \nearrow \searrow \nearrow
\end{array}
,$
\end{center}

and

\begin{center}
$
\begin{array}{lll}
& \nearrow \searrow &  \\ 
\nearrow &  & \searrow \\ 
&  &  \\ 
&  & 
\end{array}
,\qquad 
\begin{array}{lll}
&  &  \\ 
&  &  \\ 
\searrow &  & \nearrow \\ 
& \searrow \nearrow & 
\end{array}
.$
\end{center}

\strut \strut So, we have

\begin{center}
\strut $\left| \mathcal{L}_{1}^{o}(4)\right| $ $=$ $\left|
\sum_{4}^{(1)}\right| $ $=$ $6.$
\end{center}

\strut

(\textbf{Case}: $n$ $=$ $6$) If $n$ $=$ $6,$ then $\mathcal{L}_{1}^{o}(6)$
has the following elements:

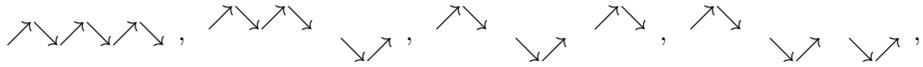
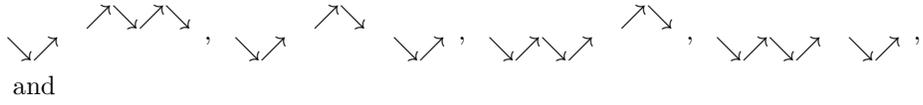
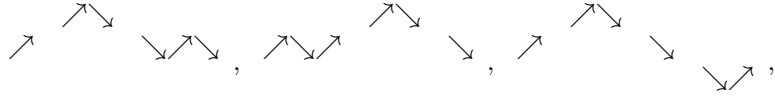
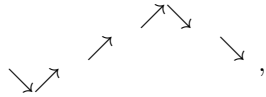
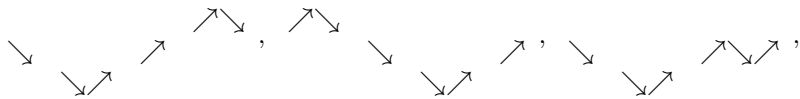
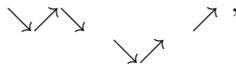
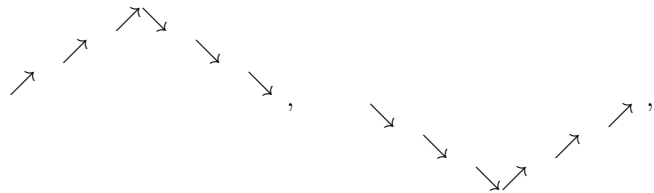
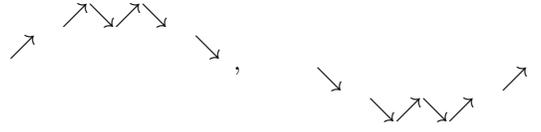
\strut

\begin{center}
$
\begin{array}{l}
\nearrow \searrow \nearrow \searrow \nearrow \searrow
\end{array}
,$ $
\begin{array}{ll}
\nearrow \searrow \nearrow \searrow &  \\ 
& \searrow \nearrow
\end{array}
,$ $
\begin{array}{lll}
\nearrow \searrow &  & \nearrow \searrow \\ 
& \searrow \nearrow & 
\end{array}
,$ $
\begin{array}{lll}
\nearrow \searrow &  &  \\ 
& \searrow \nearrow & \searrow \nearrow
\end{array}
,$
\end{center}

\strut

\begin{center}
$
\begin{array}{ll}
& \nearrow \searrow \nearrow \searrow \\ 
\searrow \nearrow & 
\end{array}
,$ $
\begin{array}{lll}
& \nearrow \searrow &  \\ 
\searrow \nearrow &  & \searrow \nearrow
\end{array}
,$ $
\begin{array}{ll}
& \nearrow \searrow \\ 
\searrow \nearrow \searrow \nearrow & 
\end{array}
,$ $
\begin{array}{ll}
&  \\ 
\searrow \nearrow \searrow \nearrow & \searrow \nearrow
\end{array}
,$
\end{center}

and

\begin{center}
$
\begin{array}{lll}
& \nearrow \searrow &  \\ 
\nearrow &  & \searrow \nearrow \searrow \\ 
&  &  \\ 
&  & 
\end{array}
,$ $
\begin{array}{lll}
& \nearrow \searrow &  \\ 
\nearrow \searrow \nearrow &  & \searrow \\ 
&  &  \\ 
&  & 
\end{array}
,$ $
\begin{array}{llll}
& \nearrow \searrow &  &  \\ 
\nearrow &  & \searrow &  \\ 
&  &  & \searrow \nearrow \\ 
&  &  & 
\end{array}
,$ $
\begin{array}{llll}
&  & \nearrow \searrow &  \\ 
& \nearrow &  & \searrow \\ 
\searrow \nearrow &  &  &  \\ 
&  &  & 
\end{array}
,$
\end{center}

\strut

\begin{center}
$
\begin{array}{llll}
&  &  &  \\ 
&  &  & \nearrow \searrow \\ 
\searrow &  & \nearrow &  \\ 
& \searrow \nearrow &  & 
\end{array}
,$ $
\begin{array}{llll}
&  &  &  \\ 
\nearrow \searrow &  &  &  \\ 
& \searrow &  & \nearrow \\ 
&  & \searrow \nearrow & 
\end{array}
,$ $
\begin{array}{lll}
&  &  \\ 
&  &  \\ 
\searrow &  & \nearrow \searrow \nearrow \\ 
& \searrow \nearrow & 
\end{array}
,$ $
\begin{array}{lll}
&  &  \\ 
&  &  \\ 
\searrow \nearrow \searrow &  & \nearrow \\ 
& \searrow \nearrow & 
\end{array}
,$
\end{center}

and

\begin{center}
$
\begin{array}{lllll}
&  & \nearrow \searrow &  &  \\ 
& \nearrow &  & \searrow &  \\ 
\nearrow &  &  &  & \searrow \\ 
&  &  &  &  \\ 
&  &  &  &  \\ 
&  &  &  & 
\end{array}
,$ $\qquad 
\begin{array}{lllll}
&  &  &  &  \\ 
&  &  &  &  \\ 
&  &  &  &  \\ 
\searrow &  &  &  & \nearrow \\ 
& \searrow &  & \nearrow &  \\ 
&  & \searrow \nearrow &  & 
\end{array}
,$
\end{center}

\strut

\begin{center}
$
\begin{array}{lll}
&  &  \\ 
& \nearrow \searrow \nearrow \searrow &  \\ 
\nearrow &  & \searrow \\ 
&  &  \\ 
&  &  \\ 
&  & 
\end{array}
,$ $\qquad 
\begin{array}{lll}
&  &  \\ 
&  &  \\ 
&  &  \\ 
\searrow &  & \nearrow \\ 
& \searrow \nearrow \searrow \nearrow &  \\ 
&  & 
\end{array}
.$
\end{center}

\strut

Therefore, we can get that

\strut

\begin{center}
$\left| \mathcal{L}_{1}^{o}(6)\right| $ $=$ $\left| \sum_{6}^{(1)}\right| $ $%
=$ $20,$
\end{center}

etc.

\strut

The above visual observations for $n$ $=$ $2,$ $4,$ $6,$ provide a way how
to compute $\left| \mathcal{L}_{1}^{o}(n)\right| ,$ for $n$ $=$ $2,$ $4,$ $%
6. $ However, in general, if $n$ is an even number greater than $6,$ the
above observation is extremely hard. But, by using the technique determined
recursively by the Pascal's triangle, we can compute $\left| \mathcal{L}%
_{1}^{o}(2n)\right| ,$ for all $n$ $\in $ $\Bbb{N}.$

\strut

As we have seen, we can realize that the cardinalities of the set $\mathcal{L%
}_{1}^{o}(2k),$ for $k$ $=$ $1,$ $2,$ $3,$ are indeed determined by the
certain entries of the Pascal's triangle,

\strut

\begin{center}
$
\begin{array}{lllllllllllll}
&  &  &  &  &  & \,1 &  &  &  &  &  &  \\ 
&  &  &  &  & 1 &  & 1 &  &  &  &  &  \\ 
&  &  &  & 1 &  & \,\frame{2} &  & 1 &  &  &  &  \\ 
&  &  & 1 &  & 3 &  & 3 &  & 1 &  &  &  \\ 
&  & 1 &  & 4 &  & \,\frame{6} &  & 4 &  & 1 &  &  \\ 
& 1 &  & 5 &  & 10 &  & 10 &  & 5 &  & 1 &  \\ 
1 &  & 6 &  & 15 &  & \frame{20} &  & 15 &  & 6 &  & 1 \\ 
&  &  &  &  &  &  &  &  &  &  &  &  \\ 
\vdots &  & \vdots &  &  &  & \vdots &  & \vdots &  & \vdots &  & \vdots
\end{array}
.$
\end{center}

\strut

Recall that, by Section 5.4, we can have that:

\strut

\begin{corollary}
For any $n$ $\in $ $\Bbb{N},$

\strut

\begin{center}
$\left| \mathcal{L}_{1}^{o}(2n)\right| $ $=\quad $ $c_{\underset{n\text{%
-times}}{\underbrace{-1,.........,-1}},\,\,\underset{n\text{-times}}{%
\underbrace{1,...........,1}}}^{(\pm 1)}$ $\quad =$ $\left|
\sum_{2n}^{(1)}\right| ,$
\end{center}

\strut

where $c_{-1,...,-1,1,...,1}^{(\pm 1)}$ are the $(-1,$ ..., $-1,$ $1,$ ..., $%
1)$-th coefficients of the $2n$-power $(a_{-1}$ $+$ $a_{1})^{2n}$ of the
binomial $a_{-1}$ $+$ $a_{1},$ for all $n$ $\in $ $\Bbb{N}.$ Equivalently,
we have that

\strut

\begin{center}
$\left| \mathcal{L}_{1}^{o}(2n)\right| $ $=$ $_{2n}C_{n}$ $=$ $\left|
\sum_{2n}^{(1)}\right| ,$ for all $n$ $\in $ $\Bbb{N},$
\end{center}

\strut

where $_{m}C_{k}$ $\overset{def}{=}$ $\frac{m!}{k!(m-k)!},$ for all $k$ $%
\leq $ $n$ in $\Bbb{N}.$ $\square $
\end{corollary}

\strut

Indeed, the coefficients in the above corollary represents the middle
entries of $2n$-levels of the Pascal's triangle, for all $n$ $\in $ $\Bbb{N}%
. $

\strut

\begin{example}
$\left| \mathcal{L}_{1}^{o}(8)\right| $ $=$ $70$ $=$ $\left|
\sum_{8}^{(1)}\right| ,\quad \left| \mathcal{L}_{1}^{o}(10)\right| $ $=$ $%
252 $ $=$ $\left| \sum_{10}^{(1)}\right| ,$

\strut

$\left| \mathcal{L}_{1}^{o}(12)\right| $ $=$ $924$ $=$ $\left|
\sum_{12}^{(1)}\right| ,\quad \left| \mathcal{L}_{1}^{o}(14)\right| $ $=$ $%
3432$ $=$ $\left| \sum_{14}^{(1)}\right| ,$

\strut

$\left| \mathcal{L}_{1}^{o}(16)\right| $ $=$ $12870$ $=$ $\left|
\sum_{16}^{(1)}\right| ,$ $\left| \mathcal{L}_{1}^{o}(18)\right| $ $=$ $%
48620 $ $=$ $\left| \sum_{18}^{(1)}\right| ,$

\strut \strut

$\left| \mathcal{L}_{1}^{o}(20)\right| $ $=$ $184756$ $=$ $\left|
\sum_{20}^{(1)}\right| ,$ etc.
\end{example}

\strut \strut

\strut \strut

\strut \strut

\strut

\strut

\strut

\strut \textbf{References}

\strut

\begin{quote}
\strut

{\small [1] \ \ A. G. Myasnikov and V. Shapilrain (editors), Group Theory,
Statistics and Cryptography, Contemporary Math, 360, (2003) AMS.}

{\small [2] \ \ A. Gibbons and L. Novak, Hybrid Graph Theory and Network
Analysis, ISBN: 0-521-46117-0, (1999) Cambridge Univ. Press.}

{\small [3]\strut \ \ \ \strut B. Solel, You can see the arrows in a Quiver
Operator Algebras, (2000), preprint.}

{\small [4] \ \ C. W. Marshall, Applied Graph Theory, ISBN: 0-471-57300-0
(1971) John Wiley \& Sons}

{\small [5] \ \ \strut D.Voiculescu, K. Dykemma and A. Nica, Free Random
Variables, CRM Monograph Series Vol 1 (1992).\strut }

{\small [6] \ \ D.W. Kribs and M.T. Jury, Ideal Structure in Free
Semigroupoid Algebras from Directed Graphs, preprint.}

{\small [7] \ \ D.W. Kribs, Quantum Causal Histories and the Directed Graph
Operator Framework, arXiv:math.OA/0501087v1 (2005), Preprint.}

{\small [8] \ \ F. Balacheff, Volume Entropy, Systole and Stable Norm on
Graphs, arXiv:math.MG/0411578v1, (2004) Preprint.}

{\small [9] \ \ G. C. Bell, Growth of the Asymptotic Dimension Function for
Groups, (2005) Preprint.}

{\small [10]\ I. Cho, Graph von Neumann algebras, ACTA. Appl. Math, 95,
(2007) 95 - 135.}

{\small [11]\ I. Cho, Characterization of Free Blocks of a right graph von
Neumann algebra, Compl. An. \& Op. theo (2007) To be Appeared.}

{\small [12]\ I. Cho, Direct Producted }$W^{*}${\small -Probability Spaces
and Corresponding Free Stochastic Integration, B. of KMS, 44, No. 1, (2007),
131 - 150.}

{\small [13] I. Cho, Vertex-Compressed Algebras of a Graph von Neumann
Algebra, (2007) Submitted to ACTA. Appl. Math.}

{\small [14] I. Cho, Group-Freeness and Certain Amalgamated Freeness, J. of
KMS, 45, no. 3, (2008) 597 - 609.}

{\small [15]\ I. Cho and Palle Jorgensen, }$C^{*}${\small -Algebras
Generated by Partial Isometries, JAMC, (2008) To Appear.}

{\small [16]\ I. Raeburn, Graph Algebras, CBMS no 3, AMS (2005).}

{\small [17]\ P. D. Mitchener, }$C^{*}${\small -Categories, Groupoid
Actions, Equivalent KK-Theory, and the Baum-Connes Conjecture,
arXiv:math.KT/0204291v1, (2005), Preprint.}

{\small [18] R. Scapellato and J. Lauri, Topics in Graph Automorphisms and
Reconstruction, London Math. Soc., Student Text 54, (2003) Cambridge Univ.
Press.}

{\small [19] R. Exel, A new Look at the Crossed-Product of a }$C^{*}${\small %
-algebra by a Semigroup of Endomorphisms, (2005) Preprint.}

{\small [20] R. Gliman, V. Shpilrain and A. G. Myasnikov (editors),
Computational and Statistical Group Theory, Contemporary Math, 298, (2001)
AMS.}

{\small [21] R. Speicher, Combinatorial Theory of the Free Product with
Amalgamation and Operator-Valued Free Probability Theory, AMS Mem, Vol 132 ,
Num 627 , (1998).}

{\small [22] S. H. Weintraub, Representation Theory of Finite Groups:
Algebra and Arithmetic, Grad. Studies in Math, vo. 59, (2003) AMS.}

{\small [23] V. Vega, Finite Directed Graphs and }$W^{*}${\small %
-Correspondences, (2007) Ph. D thesis, Univ. of Iowa.}

{\small [24] W. Dicks and E. Ventura, The Group Fixed by a Family of
Injective Endomorphisms of a Free Group, Contemp. Math 195, AMS.}

{\small [25] F. Radulescu, Random Matrices, Amalgamated Free Products and
Subfactors of the von Neumann Algebra of a Free Group, of Noninteger Index,
Invent. Math., 115, (1994) 347 - 389.}

{\small [26] D. A. Lind, Entropies of Automorphisms of a Topological Markov
Shift, Proc. AMS, vo 99, no 3, (1987) 589 - 595.}

{\small [27] D. A. Lind and B. Marcus, An Introduction to Symbolic Dynamics
and Coding, (1995) Cambridge Univ. Press.}

{\small [28] D. A. Lind and S. Tuncel, A Spanning Tree Invariant for Markov
Shifts, IMA Vol. Math. Appl., vo 123, (2001), 487 - 497.}

{\small [29] D. A. Lind and K. Schmidt, Symbolic and Algebraic Dynamical
Systems, Handbook of Dynamical System, Vol.\TEXTsymbol{\backslash}1A, (2002)
765 - 812.}

{\small [30] R. V. Kadison and J. R. Ringrose, Fundamentals of the Theory of
Operator Algebra, Grad. Stud. Math., vo. 15, (1997) AMS.}

{\small [31] D. E. Dutkay and P. E. T. Jorgensen, Iterated Function Systems,
Ruelle Operators and Invariant Projective Measures,
arXiv:math.DS/0501077/v3, (2005) Preprint.}

{\small [32] P. E. T. Jorgensen, Use of Operator Algebras in the Analysis of
Measures from Wavelets and Iterated Function Systems, (2005) Preprint.}

{\small [33] D. Guido, T. Isola and M. L. Lapidus, A Trace on Fractal Graphs
and the Ihara Zeta Function, arXiv:math.OA/0608060v1, (2006) Preprint.}

{\small [34] P. Potgieter, Nonstandard Analysis, Fractal Properties and
Brownian Motion, arXiv:math.FA/0701649v1, (2007) Preprint.}

{\small [35] L. Bartholdi, R. Grigorchuk, and V. Nekrashevych, From Fractal
Groups to Fractal Sets, arXiv:math.GR/0202001v4, (2002) Preprint.}

{\small [36] I. Cho, The Moments of Certain Perturbed Operators of the
Radial Operator of the Free Group Factor }$L(F_{N})${\small , JAA, 5, no. 3,
(2007) 137 - 165.}

{\small [37] I. Cho and P. E. T. Jorgensen, }$C^{*}${\small -Subalgebras
Generated by Partial Isometries in }$B(H)${\small , (2007) Submitted to JMP.}

{\small [38] S. Thompson and I. Cho, Powers of Mutinomials in Commutative
Algebras, (2008) (Undergraduate Research) Submitted to PMEJ.}

{\small [39] S. Thompson, C. M. Mendoza, and A. J. Kwiatkowski, and I. Cho,
Lattice Paths Satisfying the Axis Property, (2008) (Undergraduate Research)
Preprint.}

{\small [40] I. Cho, Labeling Operators of Graph Groupoids, (2008) Preprint.}

{\small [41] R. T. Powers, Heisenberg Model and a Random Walk on the
Permutation Group, Lett. Math. Phys., 1, no. 2, (1975) 125 - 130.}

{\small [42] R. T. Powers, Resistance Inequalities for }$KMS${\small -states
of the isotropic Heisenberg Model, Comm. Math. Phys., 51, no. 2, (1976) 151
- 156.}

{\small [43] R. T. Powers, Registance Inequalities for the Isotropic
Heisenberg Ferromagnet, JMP, 17, no. 10, (1976) 1910 - 1918.}

{\small [44] E. P. Wigner, Characteristic Vectors of Bordered Matrices with
Infinite Dimensions, Ann. of Math. (2), 62, (1955) 548 - 564.}

{\small [45] D. Voiculescu, Symmetries of Some Reduced Free Product }$C^{*}$%
{\small -Algebras, Lect. Notes in Math., 1132, Springer, (1985) 556 - 588.}

{\small [46] T. Shirai, The Spectrum of Infinite Regular Line Graphs, Trans.
AMS., 352, no 1., (2000) 115 - 132.}

{\small [47] J. Kigami, R. S. Strichartz, and K. C. Walker, Constructing a
Laplacian on the Diamond Fractal, Experiment. Math., 10, no. 3, (2001) 437 -
448.}

{\small [48] I. V. Kucherenko, On the Structurization of a Class of
Reversible Cellular Automata, Diskret. Mat., 19, no. 3, (2007) 102 - 121.}

{\small [49] J. L. Schiff, Cellular Automata, Discrete View of the World,
Wiley-Interscience Series in Disc. Math .\& Optimazation, ISBN:
978-0-470-16879-0, (2008) John Wiley \& Sons Press.}

{\small [50] P. E. T. Jorgensen, and M. Song, Entropy Encoding, Hilbert
Spaces, and Kahunen-Loeve Transforms, JMP, 48, no. 10, (2007)}

{\small [51] P. E. T. Jorgensen, L. M. Schmitt, and R. F. Werner, }$q$%
{\small -Canonical Commutation Relations and Stability of the Cuntz Algebra,
Pac. J. of Math., 165, no. 1, (1994) 131 - 151. }

{\small [52] A. Gill, Introduction to the Theory of Finite-State Machines,
MR0209083 (34\TEXTsymbol{\backslash}\#8891), (1962) McGraw-Hill Book Co.}

{\small [53] J. E. Hopcroft, and J. D. Ullman, Introduction to Automata
Theory, Language, and Computation, ISBN: 0-201-02988-X, (1979)
Addision-Wesley Publication Co.}

{\small [54] M. Fannes, B. Nachtergaele, and R. F. Werner, Ground States of }%
$VSB${\small -Models on Cayley Trees, J. of Statist. Phys., 66, (1992) 939 -
973.}

{\small [55] M. Fannes, B. Nachtergaele, and R. F. Werner, Finitely
Correlated States on Quantum Spin Chains, Comm. Math. Phys., 144, no. 3,
(1992) 443 - 490.}

{\small [56] M. Fannes, B. Nachtergaele, and R. F. Werner, Entropy Estimates
for Finitely Correlated States, Ann. Inst. H. Poincare. Phys. Theor., 57, no
3, (1992) 259 - 277.}

{\small [57] M. Fannes, B. Nachtergaele, and R. F. Werner, Finitely
Correlated Pure States, J. of Funt. Anal., 120, no 2, (1994) 511 - 534.}

{\small [58] J. Renault, A Groupoid Approach to }$C^{*}${\small -Algebras,
Lect. Notes in Math., 793, ISBN: 3-540-09977-8, (1980) Springer.}

{\small [59] S. Sakai, }$C^{*}${\small -Algebras and }$W^{*}${\small %
-Algebras, MR number: MR0442701, (1971) Springer-Verlag.}
\end{quote}

\end{document}